\documentclass[a4paper,12pt]{article}
\usepackage[english]{babel}
\usepackage{amsmath,amsfonts,amssymb,amscd,color,epsfig,amsthm,enumerate,
calc,comment}
\usepackage{graphicx}

\usepackage[all]{xy}
\usepackage{todonotes}

 \definecolor{ao}{rgb}{0.0, 0.5, 0.0}
\def\AA{{\cal A}}

\def\DD{{\cal D}}
\def\EE{{\cal E}}

\def\GG{{\cal G}}

\def\GP{{\cal G\!P}}
\def\GU{{\cal G\!U}}
\def\GY{{\cal G\!Y}}

\def\LL{{\cal L}}
\def\MM{{\cal M}}
\def\NN{{\cal N}}
\def\OO{{\cal O}}
\def\PP{{\cal P}}
\def\PPP{{\cal PP}}

\newcommand{\dY}{\ensuremath{\operatorname{Y}}} %Yoccoz dynamical puzzle piece
\newcommand{\pY}{\ensuremath{\operatorname{YY}}} %Yoccoz parameter puzzle piece
\newcommand{\dP}{\ensuremath{\operatorname{P}}} %Parabolic dynamical puzzle piece
\newcommand{\pP}{\ensuremath{\operatorname{PP}}} %Parabolic parameter puzzle piece
\newcommand{\pPP}{\ensuremath{\operatorname{\cal PP}}} %Parabolic parameter puzzle
\newcommand{\dN}{\ensuremath{\operatorname{N}}} %dynamical Yoccoz nest
\newcommand{\pN}{\ensuremath{\operatorname{\cal N}}} %Parameter Yoccoz nest

\newcommand{\RMone}[1]{{{\cal R}_{{#1}}^\Mone}}
\def\RR{{\cal R}}

\newcommand\TTT{{\cal T}}

\def\UY{{\cal UY}}

\def\W{{\mathcal W}}

\def\YY{{\cal Y}}
\def\ZZ{{\cal Z}}

\def\al{\alpha}
\def\be{\beta}
\def\ga{\gamma}

\newcommand{\bd}{\ensuremath{\operatorname{\partial}}}
\newcommand{\bftheta}{\ensuremath{\operatorname{\bf\theta}}}

\newcommand{\Bc}{\ensuremath{{B^\star}}}

\newcommand{\Bla}{\ensuremath{{\mathcal{B}l}}}

\newcommand{\Bstar}{\ensuremath{{B_*}}}
\newcommand{\cc}{\ensuremath{{c^\star}}}

\newcommand{\ckg}{\ensuremath{{\check g}}}
\newcommand{\ckga}{\ensuremath{{\check \ga}}}

\newcommand{\ckP}{\ensuremath{\widetilde P}}
\newcommand{\ckS}{\ensuremath{\widetilde S}}
\newcommand{\ckX}{\ensuremath{\widetilde X}}
\newcommand{\ckY}{\ensuremath{\widetilde Y}}

\newcommand{\cvE}[1]{\ensuremath{{E'_{#1}}}}

\newcommand{\cvG}[1]{\ensuremath{{G'_{#1}}}}

\newcommand{\cz}{\ensuremath{{c_0}}}
\newcommand{\cE}[1]{\ensuremath{{E^*_{#1}}}}
\newcommand{\cG}[1]{\ensuremath{{G^*_{#1}}}}
\newcommand{\C}{{\ensuremath{\mathbb  C}}}
\newcommand{\Card}{\ensuremath{{Card}}}
\newcommand{\Cbar}{\ensuremath{\overline{\mathbb  C}}}

\newcommand{\Chat}{\ensuremath{\widehat{\mathbb  C}}}

\newcommand{\diam}{\ensuremath{{\operatorname{diam}}}}

\newcommand{\D}{\ensuremath{{\mathbb   D}}}
\newcommand{\Dbar}{\ensuremath{{\overline{\mathbb  D}}}}

\newcommand{\e}{\ensuremath{{\operatorname{e}}}}
\newcommand{\eps}{{\epsilon}}
\newcommand{\End}{\ensuremath{\operatorname{End}}}

\newcommand{\Ga}{\ensuremath{{\Gamma}}}

\renewcommand{\H}{\ensuremath{\mathbb  H}}

\newcommand{\Hbz}{\ensuremath{{\bf{\ov H}_0}}}

\newcommand{\Hpb}{\ensuremath{{{\overline{\mathbb  H}}_+}}}
\newcommand{\Hplus}{\ensuremath{{{\mathbb  H}_+}}}

\newcommand{\Hz}{\ensuremath{{\bf H_0}}}
\newcommand{\itemref}[1]{\ref{#1}}

\newcommand{\mapfromto}[3]{\hbox{\ensuremath{#1 : #2 \longrightarrow #3}}}
\newcommand{\mtwo}{\ensuremath{{\operatorname{m_2}}}}

\newcommand{\Mbrot}{\ensuremath{{\bf M}}}

\newcommand{\Mod}{\ensuremath{{\operatorname{mod}}}}
\newcommand{\Mone}{\ensuremath{{\bf M_1}}}

\newcommand{\N}{\ensuremath{\mathbb  N}}

\newcommand{\ov}[1]{{\overline{#1}}}

\newcommand{\Peroo}{\ensuremath{{Per_1(1)}}}

\newcommand{\R}{\ensuremath{\mathbb  R}}

\newcommand{\Rminus}{\ensuremath{{{\mathbb  R}_-}}}

\newcommand{\siminfstar}{\ensuremath{\sim_\infty^\bigstar}}
\newcommand{\siminf}{{\sim_\infty}}
\newcommand{\siminfB}{{\sim_\infty^B}}
\newcommand{\siminfc}{{\sim_\infty^c}}

\newcommand{\siminfT}{{\sim_\infty^T}}

\newcommand{\Sm}{\ensuremath{{\setminus}}}
\newcommand{\Sen}{\ensuremath{{{\mathbb  S}^1}}}

\newcommand{\T}{\ensuremath{\mathbb T}}

\newcommand{\veps}{\ensuremath{{\underline{\epsilon}}}}
\newcommand{\vepsilon}{\veps}
\newcommand{\vepso}{\ensuremath{{\underline{\epsilon}^1}}}
\newcommand{\vepst}{\ensuremath{{\underline{\epsilon}^2}}}

\newcommand{\whP}{\ensuremath{{\widehat P}}}

\newcommand{\whtheta}{\ensuremath{{\widehat\theta}}}

\newcommand{\whz}{\ensuremath{\widehat z}}

\newcommand{\wH}{\ensuremath{\widehat H}}

\newcommand{\wtOm}{\ensuremath{\widetilde\Omega}}

\newcommand{\WMbrotpq}{\ensuremath{{\W^\Mbrot(p/q)}}}
\newcommand{\WMonepq}{\ensuremath{{\W^\Mone(p/q)}}}

\def\al{\alpha}
\def\be{\beta}

\def\ga{\gamma}

\def\si{\sigma}
\def\Si{\Sigma}

\def\eps{\epsilon}
\def\la{\lambda}
\def\La{\Lambda}

\def\de{\delta}
\def\om{\omega}
\def\Om{\Omega}

\def\R{\mbox{$\mathbb R$}}
\def\H{\mbox{$\mathbb H$}}

\def\C{\mbox{$\mathbb C$}}
\def\TT{\mbox{$\mathbb T$}}
\def\Z{\mbox{$\mathbb Z$}}

\def\N{\mbox{$\mathbb N$}}

\def\D{\mbox{$\mathbb D$}}

\def\Bottcher{B{\"o}ttcher}

\def\Mobius{M{\"o}bius} 

\newtheorem{newthm}{Theorem}
\newtheorem{theorem}{Theorem}[section]
\newtheorem{lemma}[theorem]{Lemma}
\newtheorem{proposition}[theorem]{Proposition}
\newtheorem{corollary}[theorem]{Corollary}
\newtheorem{Remark}[theorem]{Remark}
\newtheorem{defthm}[theorem]{Definition and Theorem}
\newtheorem{definition}[theorem]{Definition}
\newtheorem{example}[theorem]{Example}

\newcommand{\ENUM}{\begin{enumerate}}
\newcommand{\ENUMa}{\begin{enumerate}[a.]}
\newcommand{\ENUMA}{\begin{enumerate}[A.]}
\newcommand{\ENUMD}{\begin{enumerate}[D1.]}
\newcommand{\ENUMDp}{\begin{enumerate}[D1'.]}
\newcommand{\ENUMDB}{\begin{enumerate}[DB1.]}
\newcommand{\ENUMDBp}{\begin{enumerate}[DB1'.]}
\newcommand{\ENUMi}{\begin{enumerate}[i)]}
\newcommand{\ENDENUM}{\end{enumerate}}
\newcommand{\ITMZ}{\begin{itemize}}
\newcommand{\ENDITMZ}{\end{itemize}}
\newcommand{\REFEQN}[1] { \begin{equation}\label{#1} }
\newcommand{\ENDEQN}{\end{equation}}
\newcommand{\THM}{\begin{theorem}}
\newcommand{\REFEXA}[1] { \begin{example}\label{#1} }
\newcommand{\ENDEXA}{\end{example}}
\newcommand{\REFTHM}[1] { \begin{theorem}\label{#1} }
\newcommand{\RREFTHM}[2] { \begin{theorem}[#1]\label{#2} }
\newcommand{\ENDTHM}{\end{theorem}}
\newcommand{\REFNTH}[1] { \begin{newthm}\label{#1} }
\newcommand{\ENDNTH}{\end{newthm}}
\newcommand{\REFPROP}[1]{\begin{proposition}\label{#1} }
\newcommand{\PROP}{\begin{proposition}}
\newcommand{\ENDPROP}{\end{proposition} }
\newcommand{\REFDEF}[1]{\begin{definition}\label{#1} }
\newcommand{\DEF}{\begin{definition}}
\newcommand{\ENDDEF}{\end{definition} }
\newcommand{\REFLEM}[1]{\begin{lemma}\label{#1} }
\newcommand{\RREFLEM}[2]{\begin{lemma}[#1]\label{#2} }
\newcommand{\LEM}{\begin{lemma}}
\newcommand{\ENDLEM}{\end{lemma} }
\newcommand{\REFCOR}[1]{\begin{corollary}\label{#1} }
\newcommand{\COR}{\begin{corollary}}
\newcommand{\ENDCOR}{\end{corollary} }
\newcommand{\CONJ}{\begin{conjecture}}
\newcommand{\ENDCONJ}{\end{conjecture} }
\newcommand{\REFREM}[1]{\begin{Remark}\label{#1} }
\newcommand{\REM}{\begin{Remark}}
\newcommand{\ENDREM}{\end{Remark} }
\newcommand{\REFDEFTHM}[1] { \begin{defthm}\label{#1} }
\newcommand{\ENDDEFTHM}{\end{defthm}}

\newcommand{\corref}[1]{Corollary~\ref{#1}}
\newcommand{\remref}[1]{Remark~\ref{#1}}
\newcommand{\defref}[1]{Definition~\ref{#1}}

\newcommand{\lemref}[1]{Lemma~\ref{#1}}
\newcommand{\thmref}[1]{Theorem~\ref{#1}}
\newcommand{\propref}[1]{Proposition~\ref{#1}}
\newcommand{\secref}[1]{Section~\ref{#1}}
\newcommand{\subsecref}[1]{Subsection~\ref{#1}}

\newcommand{\figref}[1]{Fig.~\ref{#1}}
\newcommand{\Itemref}[1]{{\itemref{#1}.}}
\newcommand{\ItemDref}[1]{{{\rm{D\itemref{#1}.}}}}
\newcommand{\ItemDpref}[1]{{{\rm{D\itemref{#1}'.}}}}
\newcommand{\ItemDBref}[1]{{{\rm{DB\itemref{#1}.}}}}
\newcommand{\ItemDBpref}[1]{{{\rm{DB\itemref{#1}'.}}}}
\newcommand{\Itemiref}[1]{{\emph{\itemref{#1})}}}

\newcommand{\PROOF}{\begin{proof}}
\newcommand{\ENDPROOF}{\end{proof}}

\def\GPP{{\cal GP\!P}}
\newtheorem*{conjecture}{Conjecture}

\newtheorem*{thmA}{Theorem A}
\newtheorem*{thmB}{Theorem B}
\newtheorem*{thmi}{Theorem}

\title{The Parabolic Mandelbrot Set}
\author{Carsten Lunde Petersen and Pascale Roesch}
\date{\today}
\begin{document}

\maketitle

\begin{abstract}
We solve the longstanding conjecture by Milnor (1993) concerning the connectedness locus  
$\Mone$ of the family of quadratic rational maps tangent to the identity at $\infty$. 
We prove that this locus in homeomorphic to the Mandelbrot set $\Mbrot$ 
and that the homeomorphism is unique, provided it identifies maps that are 
"hybridly" conjugate on their filled-in Julia set.  
Moreover this homeomorphism from $\Mbrot$ to $\Mone$ is nowhere H{\"o}lder 
on the boundary and so can not have even locally a quasi-conformal extension to complements. 
\end{abstract}

\section{Introduction}
Dynamical systems given by iteration of holomorphic maps have attracted a lot 
of attention over the past 45 years. The simplest non trivial case being that of
iteration of quadratic polynomials conveniently normalized as $Q_c(z) = z^2 +c$, 
where $c\in\C$ is a parameter. The Julia set $J_f$ of a holomorphic map $f$ is the chaotic locus, which can be characterized in several way e.g. the minimal invariant set containing 
at least three points or the set of non-normality in the sense of Montel for the family of iterates. 
For a thorough introduction to holomorphic dynamics see e.g. \cite{Milnor_Dynamics}, 
\cite{Carleson-Gamelin} or \cite{Steinmetz}.

For the quadratic polynomial $Q_c$ we denote by $J_c$ its Julia set. 
In the quadratic family there is a natural dichotomy given by connectedness of the Julia set. 
The Mandelbrot set is the connectedness locus of this family:
$$
\Mbrot = \{c\in\C : J_c\textrm{ is connected}\,\}.
$$
In many cases computer generated images of the parameter space of 
a holomorphic family of holomorphic maps contains objects that looks like the 
Mandelbrot set. This led in the 1980-ties Douady and Hubbard to develop a theory of polynomial like maps, \cite{polylikemaps}, i.e.~proper holomorphic maps 
{\mapfromto f U {U'}}, where $U \subset\subset U'$ are Riemann surfaces isomorphic to $\D$. 
In this framework Douady and Hubbard were able to show that 
under a certain natural hypothesis on the family of holomorphic maps, 
denoted \emph{Mandelbrot-like family}
those objects are connectedness loci of the family and are 
in fact homeomorphic copies of the Mandelbrot set by 
a canonical dynamics preserving homeomorphism. 
They conjectured that the copies satisfying their hypothesis are 
moreover quasi-conformally homeomorphic to $\Mbrot$. 
McMullen took a step further and showed that the Mandelbrot set is universal, 
\cite{universalityofM}. 
And Lyubich, \cite{Lyubich} who developed a refined theory 
of polynomial-like maps of degree $2$ proved the Douady-Hubbard conjecture.  
In more detail Lyubich considered normalized polynomial-like maps with 
$U\subset\subset U' \subset\subset\C$ and called such maps Quadratic-like maps. 
He showed that within the space of germs of Quadratic-like maps, 
the connectedness locus of a Mandelbrot-like family of Quadratic-like maps 
is connected to the Mandelbrot set through a finite number of holomorphic motions. 

As example, let $\MM_2:=Rat_2/PSL(2,\C)$ 
denote the moduli space of rational maps of degree two. 
Milnor in~\cite{Milnorquad} introduced natural biholomorphic coordinates, 
i.e.~a complex structure {\mapfromto {(\si_1, \si_2)} {\MM_2} {\C^2}} on $\MM_2$, 
where $\si_1, \si_2$ are the first and second elementary symmetric functions 
in the three fixed point multipliers.
He then consider the curves 
$$
Per_1(\mu)=\{[f]\in \MM_2\mid f \hbox{ has a fixed point with multipier } \mu\}
$$ 
for $\mu\in\C$ and shows that each $Per_1(\mu)$ is a straight line in the above 
complex structure. 

For $\mu\in\D\cup\{1\}$ denote by $\Mbrot_\mu$ 
the connectedness locus in $Per_1(\mu)$
$$
\Mbrot_{\mu}:=\{[f]\in Per_1(\mu)\mid J_f \hbox{ is connected\}}. 
$$

It follows from Lyubich's theory of Quadratic-like maps 
that $\Mbrot_{\mu}$ is quasi-conformally homeomorphic to
$\Mbrot$ for $\mu\in\D$. 

\begin{figure}\label{fhomeoM}
\begin{center}\includegraphics[height=2 in]{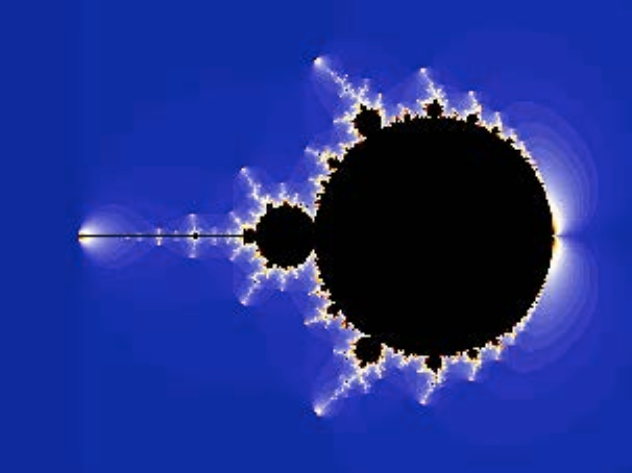}
\includegraphics[height=2 in]{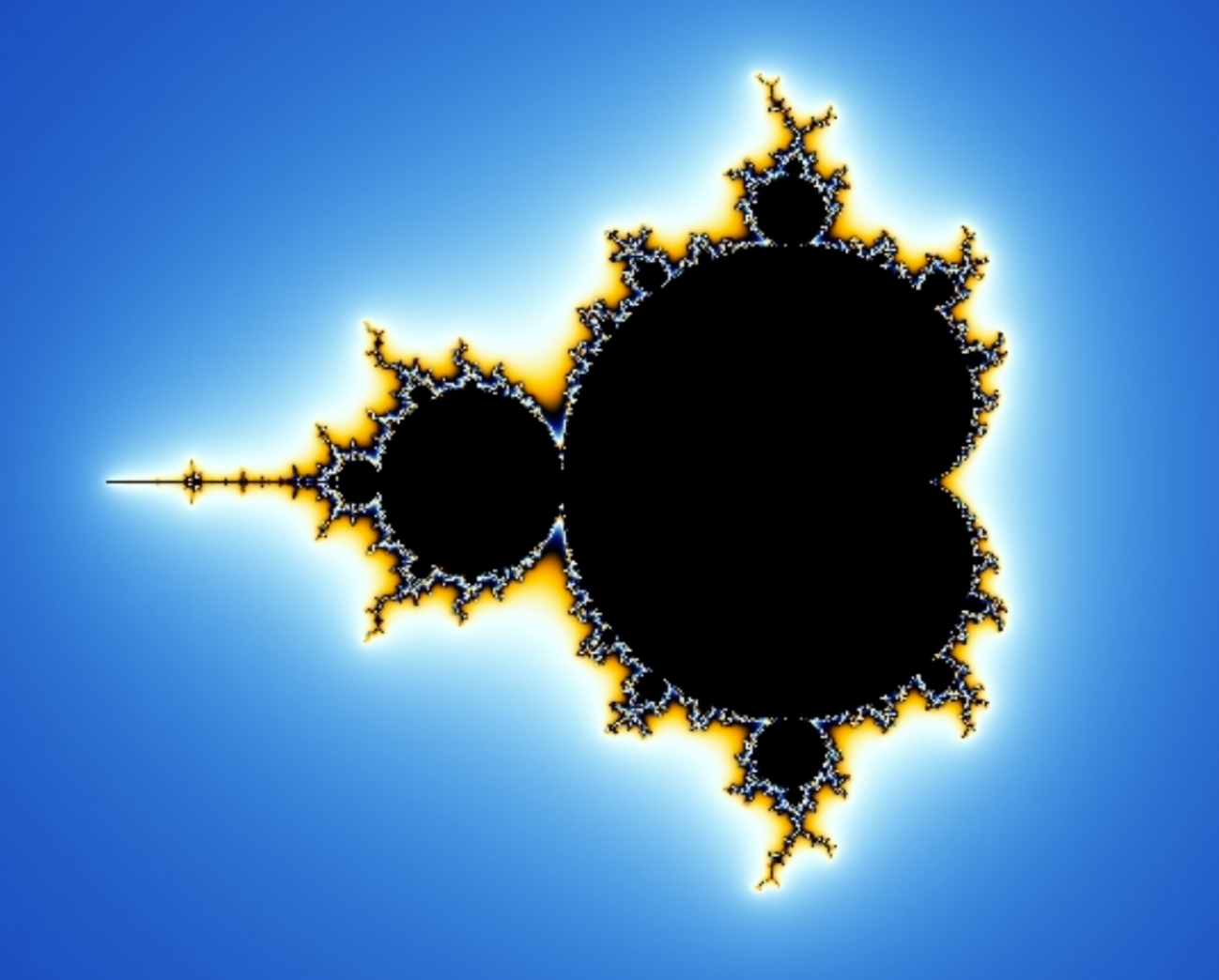}  
\end{center}
\caption{The parabolic connectedness locus $\Mone$ (left) and the $\Mbrot$ (right)}
\end{figure}

In the afore mentioned paper \cite[1993]{Milnorquad} Milnor showed pictures of the connectedness locus $\Mone$ for $\mu=1$, see \figref{f:Mone} and \figref{fhomeoM} 
for illustrations of $\Mone$. Moreover he proposed the following conjecture.
\CONJ[Milnor, 1983]
The connectedness locus $\Mone$ is homeomorphic to the Mandelbrot set $\Mbrot$.
\ENDCONJ
After this the set $\Mone$ was referred to as the Parabolic Mandelbrot set. 
However for $\mu=1$ maps are not polynomial-like, 
i.e.~the maps $f$, $[f]\in Per_1(1)$ do not posses a polynomial like restriction. 

In this paper we prove a stronger version of Milnor's conjecture.
\begin{thmA} There is a unique dynamics preserving homeomorphism 
{\mapfromto {\Phi^1} \Mbrot \Mone} between the Mandelbrot set $\Mbrot$ 
and the parabolic Mandelbrot set $\Mone$. 
Moreover $\Phi^1$ admits no quasi-conformal extension to any neighbourhod of any 
boundary point of $\Mbrot$.
\end{thmA}

\noindent In this theorem, \emph {dynamics preserving} means:
 
For $c\in \Mbrot$ and $[g]=\Phi^1(c)\in\Mone$ 
there exists a homeomorphism {\mapfromto {\rho_c} \Chat\Chat}, 
which is conformal a.e.~on $K(Q_c)$, 
which conjugates $Q_c$ to $g$ (see below Figure~\ref{fdynamics}) : 
\[\xymatrix{K(Q_c)\ar[r]^{Q_c}\ar[d]_{\rho_c} &K(Q_c)\ar[d]^{\rho_c} \\
K(g)\ar[r]_g&K(g) } \]

The statement of no quasi-conformal extension is a consequence 
of the following stronger result. 

\begin{thmB}
For a dense set of parameters in $\partial\Mbrot$ the homeomorphism $\Phi^1$
is not H{\"o}lder for any exponent. 
\end{thmB}

Copies of $\Mone$ appear in several works as a consequence of Lomonaco's theory of parabolic-like maps,
\cite{parabolic_like}, \cite{parabolic_like_families}. 
As examples the result has already found applications in the work by Bullet and Lomonaco on
algebraic correspondences \cite{BL}. 
Further examples are the lemniscate copies in the slices $Per_1(\mu)$, $\mu$ 
a root of unity, in the space 
$\MM_2$ \cite{Lomonaco-Uhre }, and  in the moduli space of cubic polynomials \cite{parabolic_like_families}, \cite{Runze1}, \cite{Runze2}. 
 
In this and many other cases the maps do not posess a 
polynomial-like restriction, because it is not possible to choose $U, U'$ with 
$U\subset\subset U'$. 
They are what is called pinched quadratic-like maps 
and such families do not form Mandelbrot-like families. 
Thus in order to prove Milnor's conjecture and more generally prove 
that more of the observed Mandelbrot-look-alikes are homeomorphic 
to the actual Mandelbrot set, 
one needs to device completely new strategies and tools. 
For a more thorough discussion of this see also \secref{Background_etc}. 

Our theorem is the first instance of a proof that a connectedness locus 
in a family of pinched quadratic-like maps is homeomorphic to $\Mbrot$. 
As a consequence of the proof we also obtain new examples of maps in the family $Per_1(1)$ with positive area Julia sets.

%{\bf \underline{The strategy of the  proof of the Main Theorem}}
%%%%%%%%%%%%%%%%%%%%%%%%%%%%%%%%%%%

The basic idea of our proof, which was laid out in the papers \cite{PR1} and \cite{PR2},  
is to develop a theory of Yoccoz puzzles for parabolic maps and 
to use the combinatorics, which is encoded in the Yoccoz puzzles, 
as a vehicle to define a homeomorphism. 
Such ideas has since then been applied in several other settings e.g. 
\cite{Dudko3}, \cite{DudkoSchleicher}. 

The structure of the paper is as follows. 
In section \ref{s:DHM1} we introduce a natural parametrization of the 
complement of $\Mone$. 
In section \ref{s:parabolicrays} we introduce the notion of parabolic rays 
originally defined in \cite{PR2} and port these to parabolic rays in parameter space. 
As for the Mandelbrot set this leads to defining wakes and limbs of $\Mone$ and 
it leads to a parabolic Yoccoz inequality. 
In section \ref{s:ParabolicPuzzles} we recall from \cite{PR2} the 
construction of parabolic puzzles, similar to Yoccoz puzzles. 
These are defined via parabolic rays in both 
dynamical space and parameter space. 
In section \ref{s:Tower} we recall the notion of combinatorial-analytic invariants 
introduced in \cite{PR1}. These allow for defining a dynamics preserving 
bijection between corresponding limbs of $\Mone$ and $\Mbrot$. 

In section \ref{s:Transfer_to_Mone} we develop a dynamical Yoccoz theory for parabolic puzzles. The basic problem compared with the classical Yoccoz theory 
is that the puzzle pieces around the $\beta$-fixed point and its preimage are 
not dynamical in the sense that each puzzle piece does not map onto the puzzle piece one level up.

In section \ref{s:PPP} we use the dynamical parabolic Yoccoz theory to develop 
a parameter parabolic Yoccoz theory. Also here the new difficulty is coming from 
the $\beta$ and $\beta'$ nests not being dynamical. 

In section \ref{Yoccoztheoremformone} we prove local connectivity of $\Mone$. 
Finally in section \ref{s:mainproof} we prove the main theorems.

Our strategy is to compare the representation of $\Mbrot$ given in \cite{PR1} 
in terms of combinatorial data together 
with analytic data with a similar representation of $\Mone$.  
The representation of $\Mbrot$ was obtained with the use of Yoccoz' theorem. 
More precisely, replacing every maximal, i.e.~level one renormalization copy of  $\Mbrot $ 
strictly inside $\Mbrot$ we obtain a tree. 
This tree is faithfully described by a space of stratified equivalence relations/laminations called towers.
The union of the equivalence relations in a tower is an equivalence relation, 
which is forward invariant under $Q_0$ and which corresponds to the co-landing pattern 
of those rays which eventually land on the $\al$-fixed point.

In \cite{PR2} we provided a dynamically defined map {\mapfromto {\Psi^1} \Mone \Mbrot}. 
This was constructed by associating combinatorial and analytic data similar to those for quadratic polynomials to each element of $\Mone$. 
The map $\Psi^1$ takes $g$ in $\Mone$ to $c= \Psi^1(g)$ such that $g$ and $Q_c$ 
have the same combinatorial and analytic data.

In this paper we prove that $\Psi^1$ is a homeomorphism. We do this in two steps.
First we prove an analogue of Yoccoz parameter theorem for $\Mone$. 
From this theorem it easily follows that $\Psi^1$ is a bijection and is continuous, 
except possibly at the boundary of the level one renormalization copies of $\Mbrot$.
The second and last step is to prove the continuity at those remaining parameters. 
The continuity of the inverse {\mapfromto {\Phi^1} \Mbrot \Mone} 
then follows from abstract reasons.

In the course of the proof we prove the aforementioned parabolic Yoccoz parameter theorem, 
which falls in two parts shrinking of limbs along the unit disk and a parameter puzzle theorem for $\Mone$. 

\section{Background and state of the art}\label{Background_etc}
By definition of $Per_1(\mu)$, one of the fixed point multipliers is $\mu$. 
The product of the two remaining fixed point multipliers define an isomorphism 
{\mapfromto {\si^\mu} {Per_1(\mu)} \C} between $Per_1(\mu)$ and $\C$. 
See \cite[Lemma 3.4]{Milnorquad} for details. 
To simplify the notation we shall henceforth identify $Per_1(\mu)$ with $\C$ via this isomorphism.
 
In  $Per_1(\mu) $, the interesting  dynamical systems are located in 
the connectedness locus $\Mbrot_{\mu}$. 
Indeed, outside  $\mathbf M_{\mu}$ the map $f$ is conjugate to the shift map on a Cantor set. 
We shall often identify $f$ and $[f]$, when there is no risk of confusion.
  
For $\mu=0$, $\Mbrot_0 $ corresponds to the classical Mandelbrot set $\Mbrot$. 
There is an extensive knowledge and literature about quadratic polynomials 
and the Mandelbrot set pioneered by Douady and Hubbard (see \cite{orsaynotesI} and \cite{orsaynotesII}). 
The family of quadratic polynomials $Q_c(z) = z^2+c$, $c\in\C$ parametrizes $Per_1(0)$. 
The Julia set $J_c = J(Q_c)$ is the common boundary of the filled Julia set 
$K_c=K(Q_c)=\{z\in \C\mid Q_c^n(z)\hbox{ is bounded} \} $ 
and the basin of infinity $ B_c(\infty):=\{z\in \C\mid Q_c^n(z)\to \infty\} = \C\Sm K_c$. 
The fixed landing point of the external ray of argument $0$ is called $\be_c$, 
the other fixed point is called $\al_c$. We have $\al_c=\be_c$ if and only if $c=\frac{1}{4}$.

The Mandelbrot set $\Mbrot=\{c\in \C\mid J(Q_c) \hbox{ is connected}\}$ 
is also the set of parameters $c\in\C$ such that $c\in K_c$. 
The product $\si^0$ of the multipliers of the two finite fixed points of $Q_c$ equals $4c$, 
so that $\Mbrot_0 = 4\Mbrot$.
The central hyperbolic component $\Card$, i.e.~the connected component 
of the interior of $\Mbrot$ containing $0$ consists of parameters $c$  for which $Q_c$ has an attracting fixed point $\al_c\in \C$. 

 Recall that a holomorphic motion of $L\subset\C$ over a complex analytical manifold 
 $\La$ with base point $\la_0\in\La$ is a map {\mapfromto H {\La\times L} \C} satisfying 
 \ENUM
 \item
 For each fixed $z\in L$ the map $\la\to H_z(\la) := H(\la, z)$ is holomorphic.
 \item
 For each $\la\in\La$ the map $z\to H^\la(z) := H(\la, z)$ is injective.
 \item
The map {\mapfromto {H^{\la_0}} L \C} is the identity.
\ENDENUM
By the celebrated $\la$-lemma each $H^\la$ is (the restriction of) a quasi-conformal homeomorphism.

For $\mu\in\D$ the connectedness locus $\Mbrot_\mu$ is 
quasi-conformally homeomorphic to $\Mbrot$. 
Indeed identifying each $Per_1(\mu)$ with $\C$ via $\si^\mu$ 
one has the following theorem:
 
\begin{thmi}[Lyubich, Uhre, Bassanelli-Berteloot]
There exists a dynamical holomorphic motion {\mapfromto \Phi {\D\times\Mbrot_0} \C} 
of $\Mbrot_0$ over $\D$ with base point $\mu_0 = 0$ such that 
$\Phi^\mu(\Mbrot_0) = \Mbrot_\mu$ for all $\mu\in\D$.
\end{thmi}   
Here dynamical means that $Q_c$ and $f\in\Phi^\mu(4c)$ have 
polynomial-like restrictions which are hybridly equivalent, 
i.e.~are conjugate by a quasi-conformal homeomorphism, 
which is conformal a.e.~on the filled-in Julia sets.

 \begin{figure}[h]
\begin{center}  
\includegraphics[height=1.5 in]{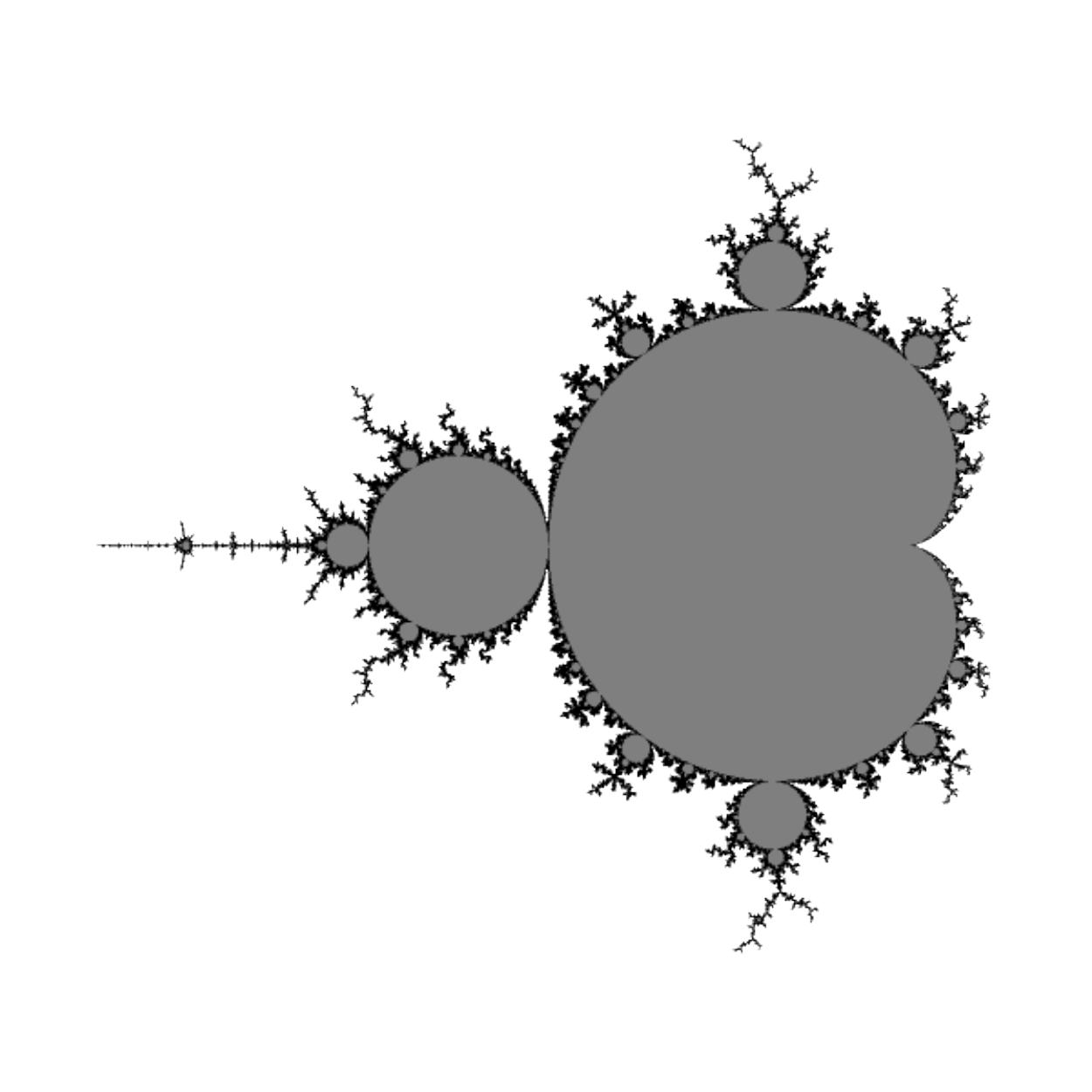} 
\includegraphics[height=1.5 in]{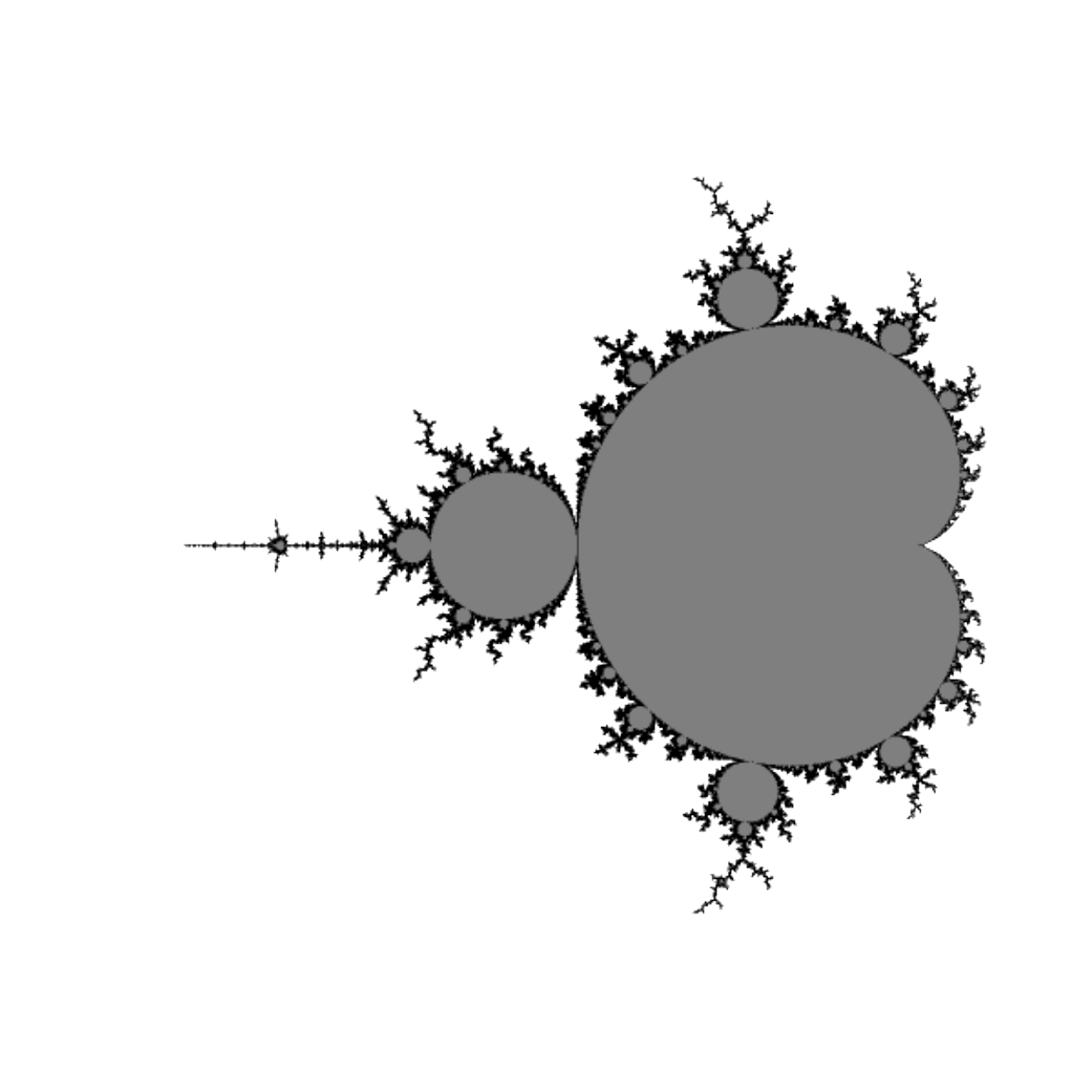}
\includegraphics[height=1.5 in]{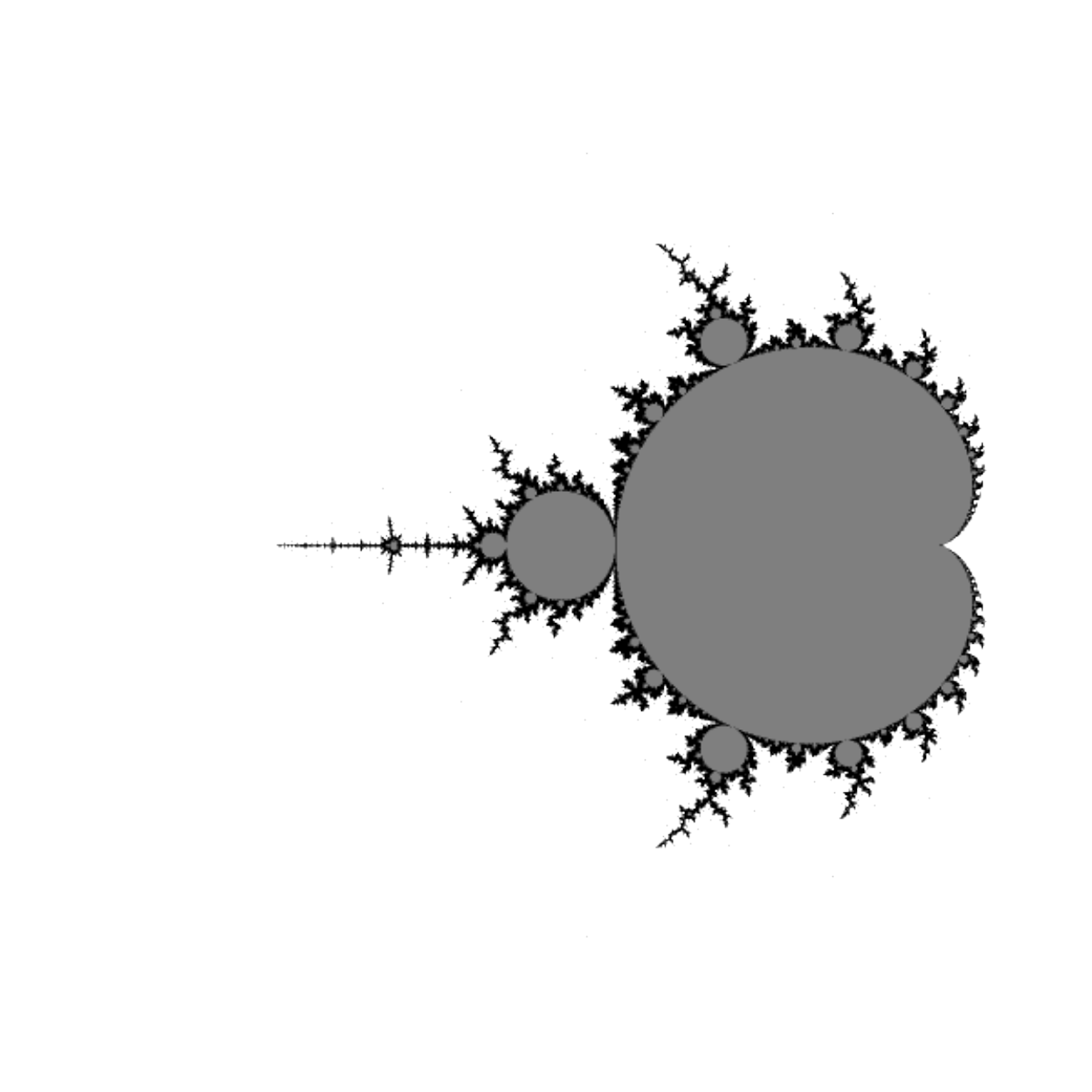}
\includegraphics[height=1.5 in]{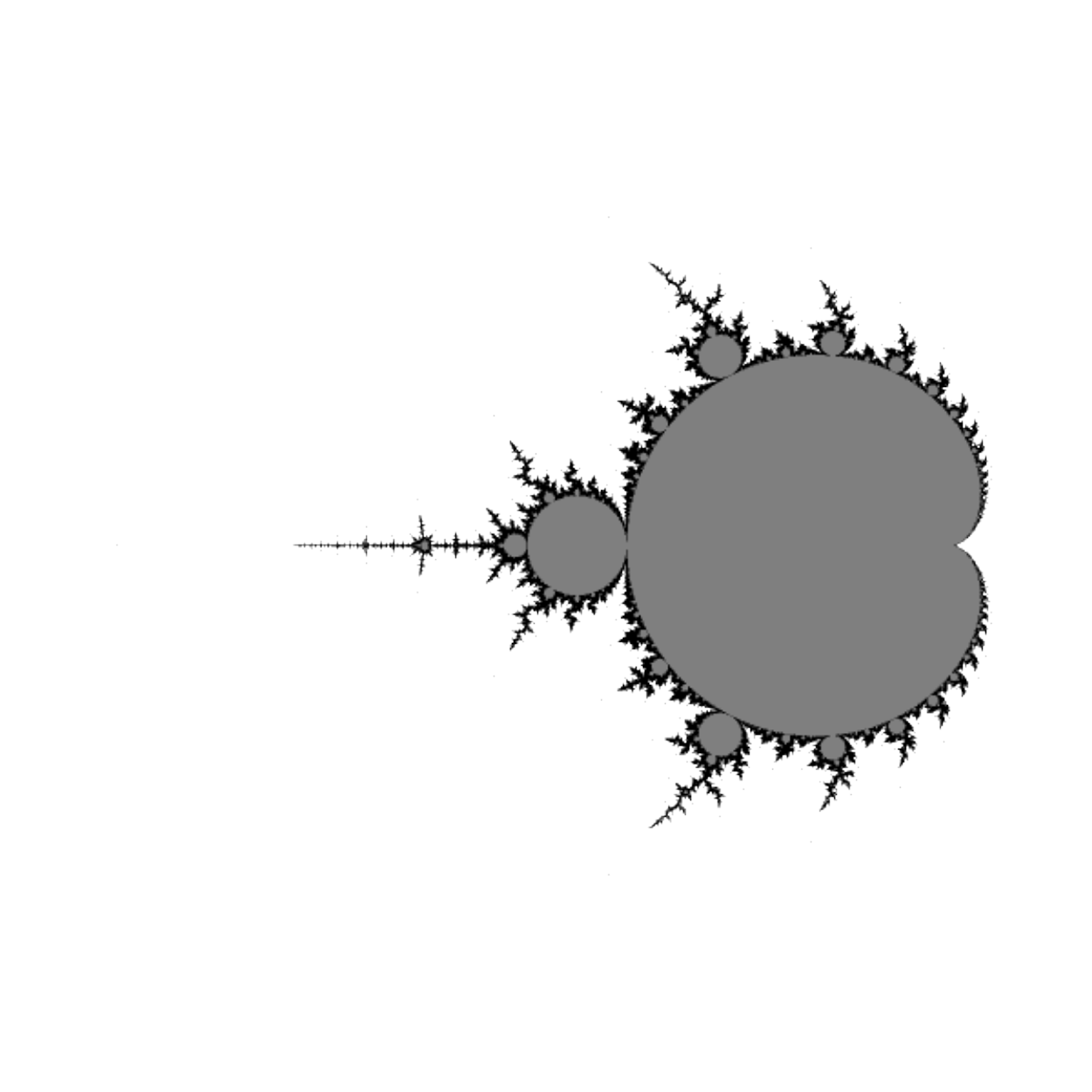}\\
\includegraphics[height=1.5 in]{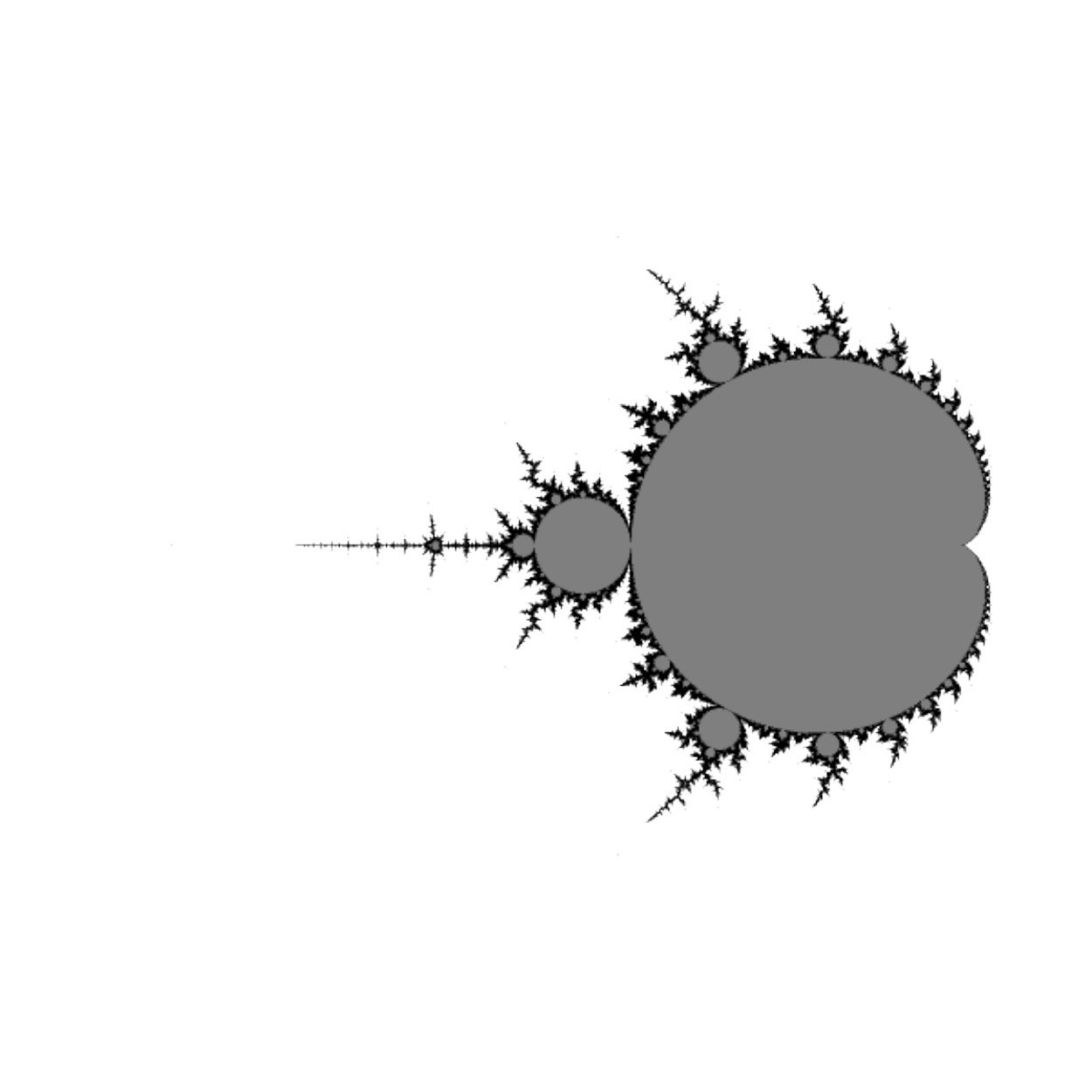}
\includegraphics[height=1.5 in]{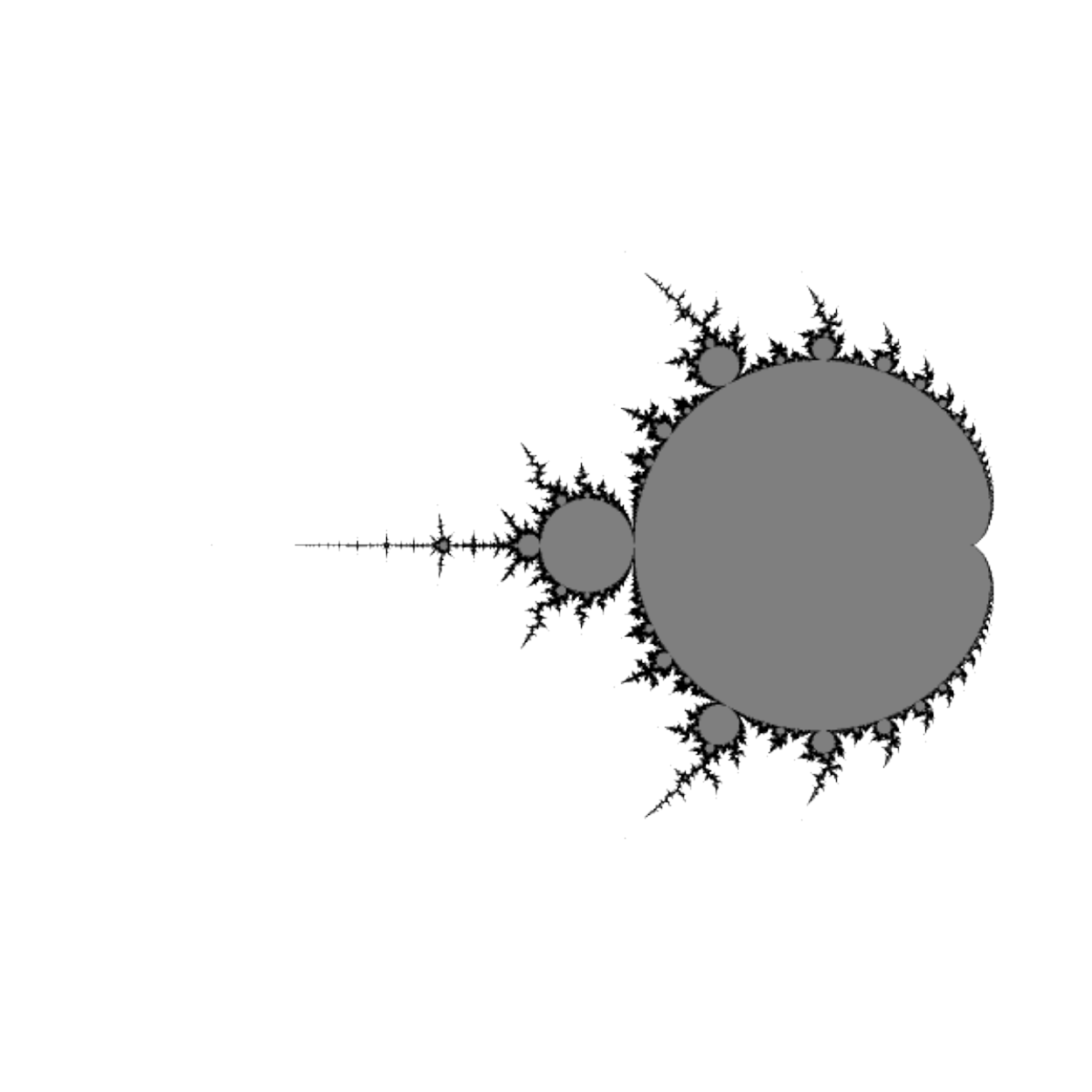}
\includegraphics[height=1.5 in]{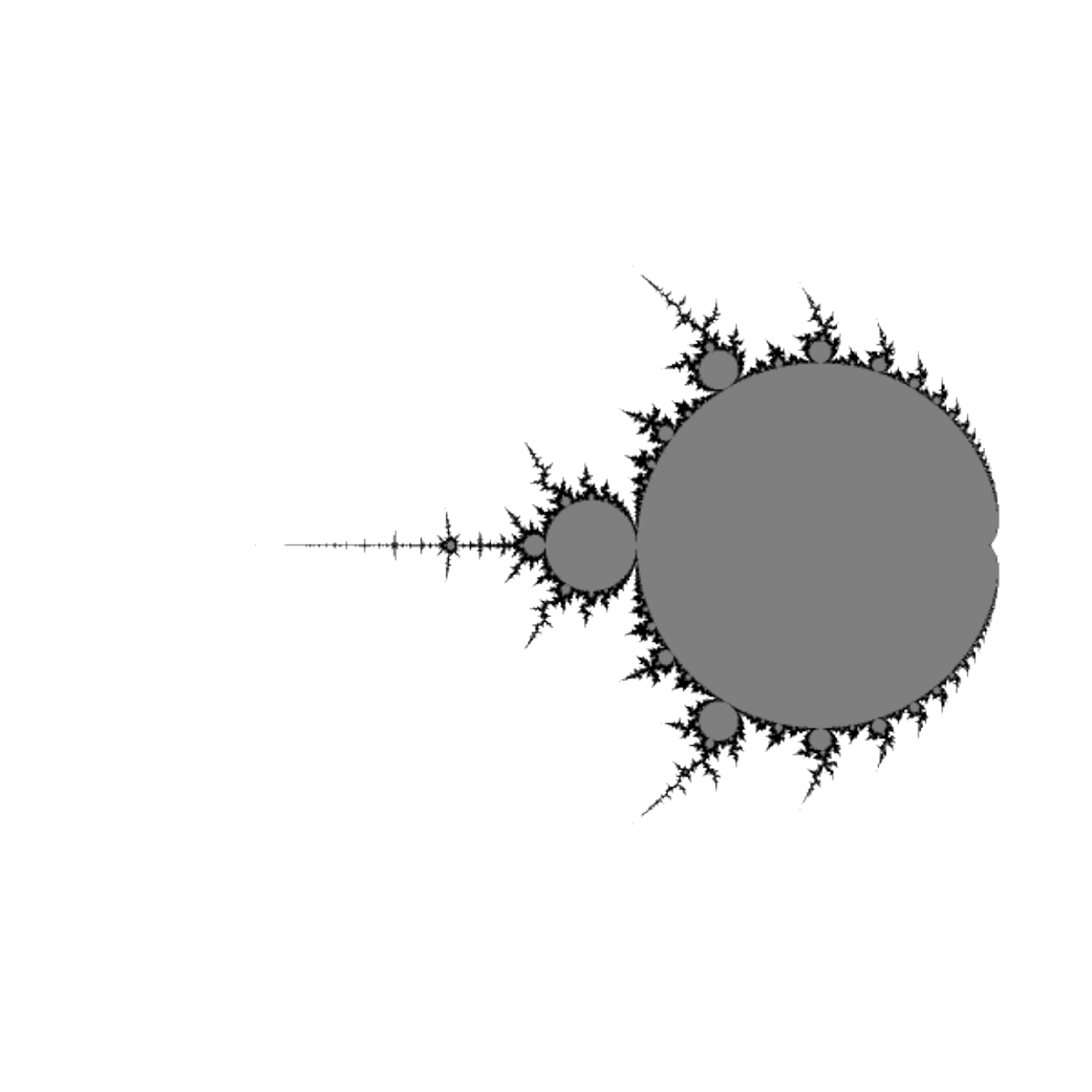}
\includegraphics[height=1.51 in]{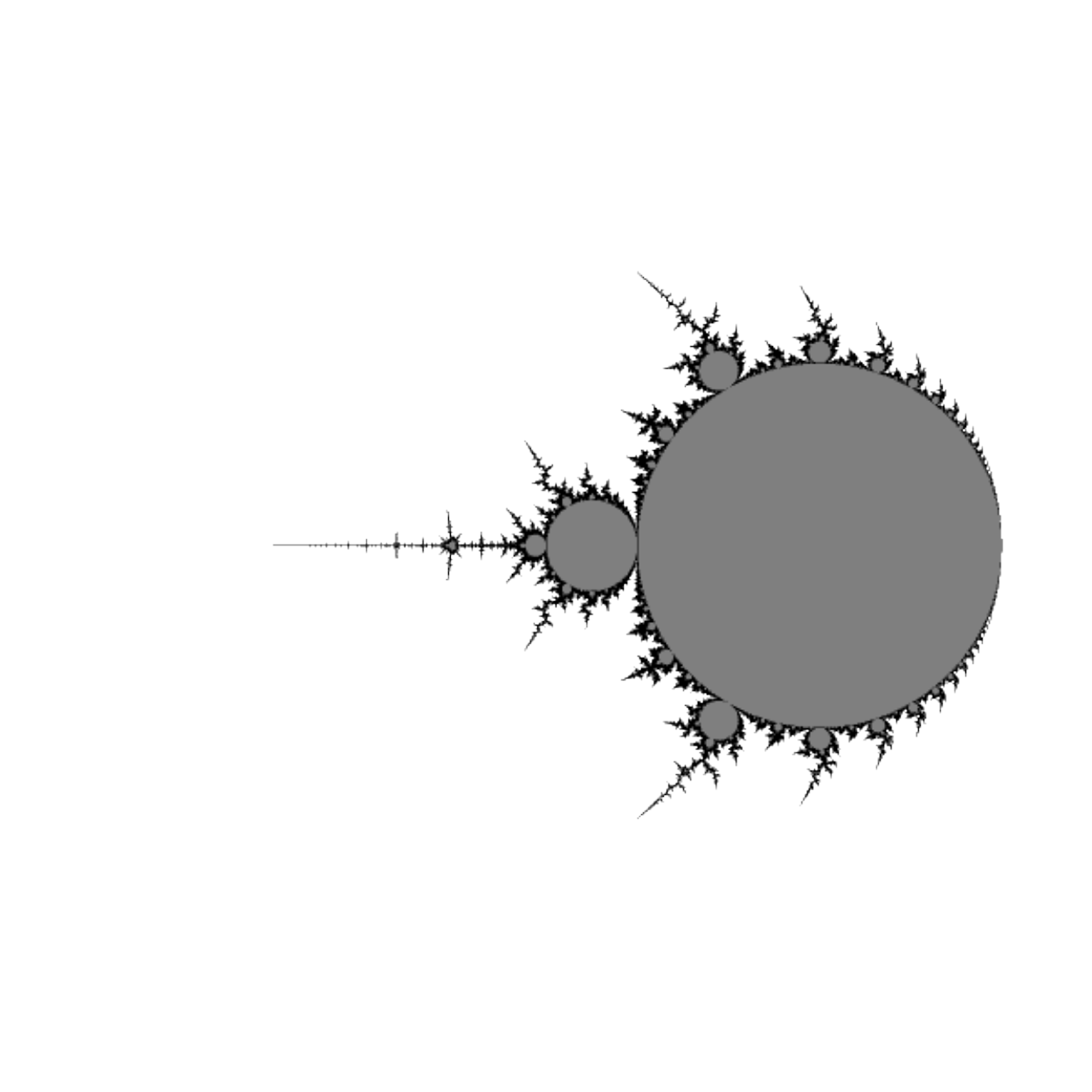}
 \end{center}
\caption{The sets  $\Mbrot_\mu$  for some parameters $\mu \in [0,1]$.}\label{fholoM} 
\end{figure}

\begin{figure}[h]
\begin{center} 
\includegraphics[height=1.45 in]{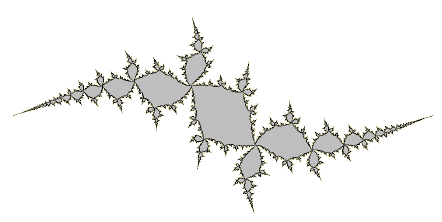}
\includegraphics[height=1.485 in]{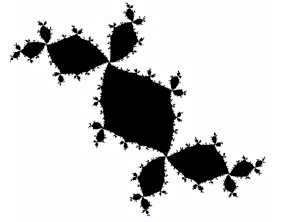}  
\end{center}
\caption{A parabolic Julia set (left) and the corresponding one in $\Mbrot$ (right)}\label{fdynamics} 
\end{figure}

\begin{figure}[h]
\centerline{\includegraphics[height=80mm]{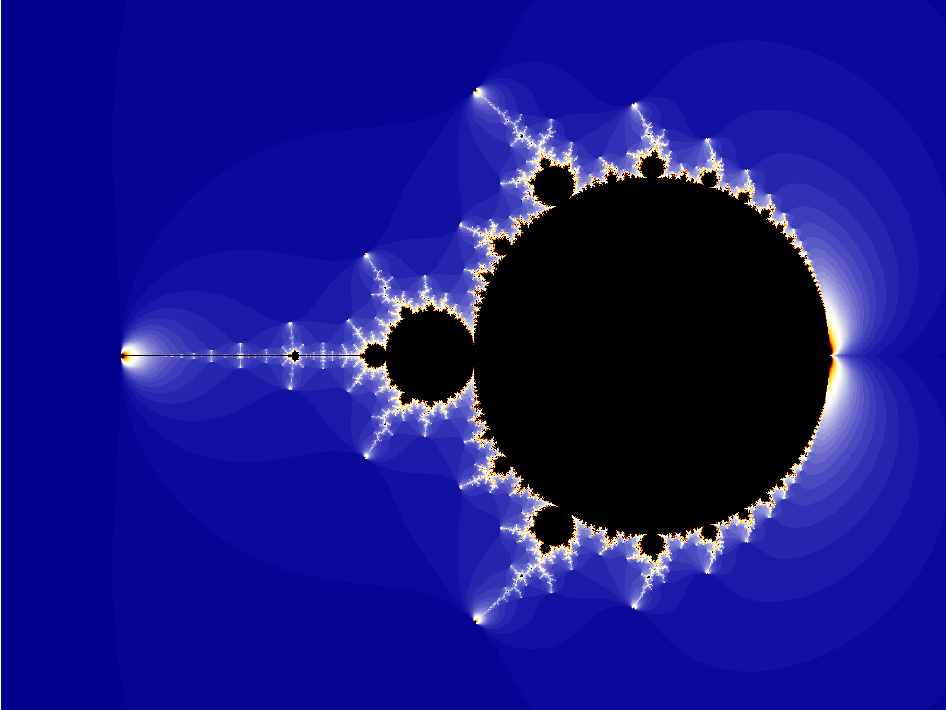}}
\caption{In black, $\Mone$, viewed in the coordinate $\sigma^1$, 
the product of the two remaining fixed point multipliers. 
Maps in $Per_1(1)$ have a degenerate fixed point. 
Hence $\si^1([g])$ is also the multiplier of the unique third fixed point $\al_g$ of $g$. 
In particular the big black disk is the unit disk, it corresponds to $\alpha_g$ being attracting. }
\label{f:Mone}
\end{figure}

The consecutive images of $\Mbrot_\mu$ for $\mu\in [0,1]$ in \figref{fholoM} illustrates why Milnor could be led to conjecture that for $\mu=1$ there is still a homeomorphism, 
see also see~\figref{fhomeoM} and~\figref{fdynamics}. 
However there is no general Theorem on the boundary values of a holomorphic motion on 
say $\D$, which would allow to conclude that there is a limiting homeomorphism 
{\mapfromto {\Phi^1}{\Mbrot_0}{\Mbrot_1}}.

The notion of filled-in Julia set for $[g]\in Per_1(1)$, needs a little clarification. 
For $\si^1(g) = 1$, the parabolic basin for the fixed point of multiplier $1$ 
has two symmetric components, choose one of them to be the 
external basin, in all other cases the parabolic basin is connected 
and so we may unambiguously call this the external basin. 
{\it The filled Julia set } $K(g)$ is the complement of the external basin.

Ha{\"i}ssinsky, \cite{Haissinsky} has developed a trans-quasi-conformal surgery 
on simply connected attracting basins:
\begin{thmi}[Haissinsky]\label{haissinskysurgery} 
Suppose $c\in M$ and that
the critical point $0$ of $Q_c$ is not recurrent to $\beta_c$.
Then there exists a trans-quasi conformal homeomorphism
{\mapfromto {h_c} \Cbar \Cbar} conformal on the interior of $K_c$
and a unique quadratic rational map $g_B(z)= z+1/z + B$ such
that $1\in \La_B(\infty)$ and such that $h_c$ conjugates $Q_c$ to
$g_B$ on $\Cbar\Sm \Delta_c$, where $\Delta_c$ is a forward
invariant topological disk in $\La_c(\infty)$ accessing $\beta_c$
quasi-radially. 
\end{thmi}

Ha{\"i}ssinsky's theorem provides a way to uniquely define $\Phi^1(c)\in\Mone$ 
for any $c\in\Mbrot$ such that the critical point is not recurrent 
to the $\beta$-fixed point $\beta_c$.
However the condition that the critical point is not recurrent to the $\beta$-fixed point 
has zero harmonic measure in $\partial\Mbrot$.
A recent result of Dudko and Lyubich states that every quadratic polynomial with an indifferent fixed point has a maximal hedgehog, 
\cite{Dudko-Lyubich}. 
According to private communication with Dudko, 
this implies that the critical point is not recurrent to the 
$\beta$ fixed point of such a polynomial. 
Hence Ha{\"i}ssinsky's theorem applies to all polynomials on the boundary of 
$\Card$.

To sum up the challenges to be overcome in order to prove Milnor's conjecture include: 
\begin{enumerate}
 \item 
As discussed above there is no straightening  result (like Douady-Hubbard straightening theorem) going from the parabolic world to the hyperbolic world 
and more generally from pinched quadratic-like maps to quadratic-like maps 
(the disks defining the polynomial like map touch at their boundary). 
In the case of $\Mone$ the work of Haissinsky, \cite{Haissinsky} using Guy David's theorem, see also below, does apply to interior points of $\Mbrot$, 
but fails to apply for $\om$-almost all boundary points,  
where $\om$ denotes the equilibrium measure on $\Mbrot$. 
 \item 
 In the case at hand with $\Mone$ there is no existing  theory, 
 which allows us to extend the holomorphic motion $\Phi$ at the boundary point $1$. 
 There is definitely no extension of the holomorphic motion in neighbourhood of $1$, 
 because $\Mbrot$ and $\Mone$ are not quasi-conformally equivalent\,;
\item There is no  complete description of  the boundary $\partial \Mbrot$ that describes the dynamics of the maps for instance by the exterior and that could have been transported. (It is conjectured that $\Mbrot$ is locally connected).
 \item There is no complete description of the dynamics inside of the interior of $\Mbrot$ that would allow us to compare the dynamics\,;
 \item  There could be queer components in $\Mone$ which do not correspond to a queer component in $\Mbrot$, i.e.~components of the interior for which all periodic points are repelling for those maps. 
 The above mentioned work of Haissinsky implies that every queer component in $\Mbrot$ if any
 would correspond to a queer component in $\Mone$ (Fatou's conjecture says that there are none in $\Mbrot$);
 \end{enumerate}
 \vskip 1em
 
\section{Basic notations and description of parameter spaces.}\label{s:DHM1}
\subsection{Basic notions for quadratic polynomials}

An essential tool in the study of the dynamics of quadratic polynomials 
and the Mandelbrot set is the notion of external rays. 
For $c\in\C$ we let $\phi_c$ denote the {\Bottcher}-coordinate conjugating $Q_c$ to $Q_0$ 
in a neighbourhood of $\infty$ normalized by being tangent to identity at $\infty$. 
The Green's function $G_c$ for $B_c(\infty)$ 
is the subharmonic function on $\C$ defined by $G_c(z) = \log|\phi_c(z)|$ on the domain of $\phi_c$, 
the recursive relation $G_c(Q_c(z)) = 2G_c(z)$ and $G_c \equiv 0$ on $K_c$. 
The {\Bottcher}-coordinate has a unique univalent extension {\mapfromto {\phi_c}{\{z|G_c(z)>G_c(0)\}}{\Chat\Sm\Dbar(\e^{G_c(0)})}}. 

Denote by $\psi_c$ the inverse of $\phi_c$. The map $\psi_c$ analytically extends along rays $\RR_\theta$, 
where $\RR_\theta$ is the straight line of angle $\theta$, $\theta\in\R/\Z=:\TT$. 
This extension stops when reaching $\Dbar$ or 
a point whose image by $\psi_c$ is a pre-critical value of $Q_c$. 
In the rest of the paper we shall let $\psi_c$ mean the maximal radial extension.
The external ray of angle $\theta$ is defined by $\RR_\theta^c=\psi_c(\RR_\theta)$. 
An external ray which stops at a pre-critical value is said to bump. 
The dynamics of $Q_c$ on rays is semi-conjugate to angle-doubling $\mtwo(\theta) = 2\theta\mod 1$.
Douady and Hubbard proved that a non-bumping (pre)-periodic ray lands 
at a (pre)-periodic repelling or parabolic orbit with the same pre-period 
and period dividing that of the ray (i.e. that of $\theta$). More precisely

\RREFTHM{Douady-Hubbard}{polyraylanding}
Let $c\in\C$. For any {\rm(}pre-{\rm)}periodic argument $\theta\in\TT$,
i.e.~$\mtwo^k(\mtwo^l(\theta))=\mtwo^l(\theta)$,
if the external ray $\RR=\RR^c_\theta$ does not bump, 
then it converges to a $Q_c$ {\rm(}pre-{\rm)}periodic point
$z\in J_c$ with $Q_c^k(Q_c^l(z))=Q_c^l(z)$.
If  the argument is periodic {\rm(}i.e.~$l=0${\rm)}, let $k'$ denote the exact period of $z$ and let $q=k/k'$.
Then the ray $\RR$ defines the combinatorial rotation number $p/q$, $(p,q)=1$ for $z$.
The periodic point $z$ is repelling or parabolic with multiplier $\e^{i2\pi p/q}$.
Moreover any other external ray landing at $z$ is also $k$-periodic and
defines the same rotation number.
\ENDTHM
 \begin{figure}[h]
\centerline{\includegraphics[height=60mm]{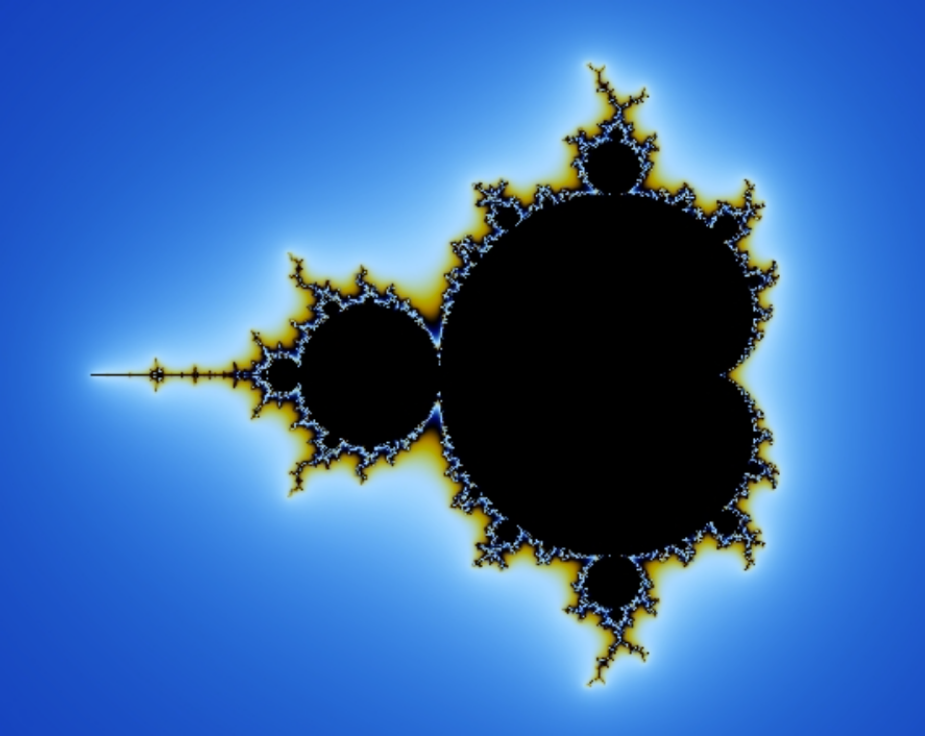}}
\caption{ Mandelbrot set -- the central cardioid with limbs attached}
\label{f:M}
\end{figure} 
 
\RREFTHM{Douady}{landingportrait}
Let $c\in\Mbrot$ and suppose $z$ is a (pre)-periodic point, $Q_c^{k'}(Q_c^l(z)) = Q_c^l(z)$, 
$l\geq 0$ and $k\geq 1$. And suppose that $w=Q_c^l(z)$ is either repelling or parabolic. 
Then $w$ has a combinatorial rotation number in the sense above. 
That is $w$ is the landing point of at least one external ray and all rays landing at $w$ form 
a single cycle under $Q_c^{k'}$. Moreover $z$ is the landing point of the same number of rays, 
each one being a preimage of a ray landing at $w$.
\ENDTHM

\begin{figure}[h]\label{f:Wake}
\centerline{\includegraphics[height=70mm]{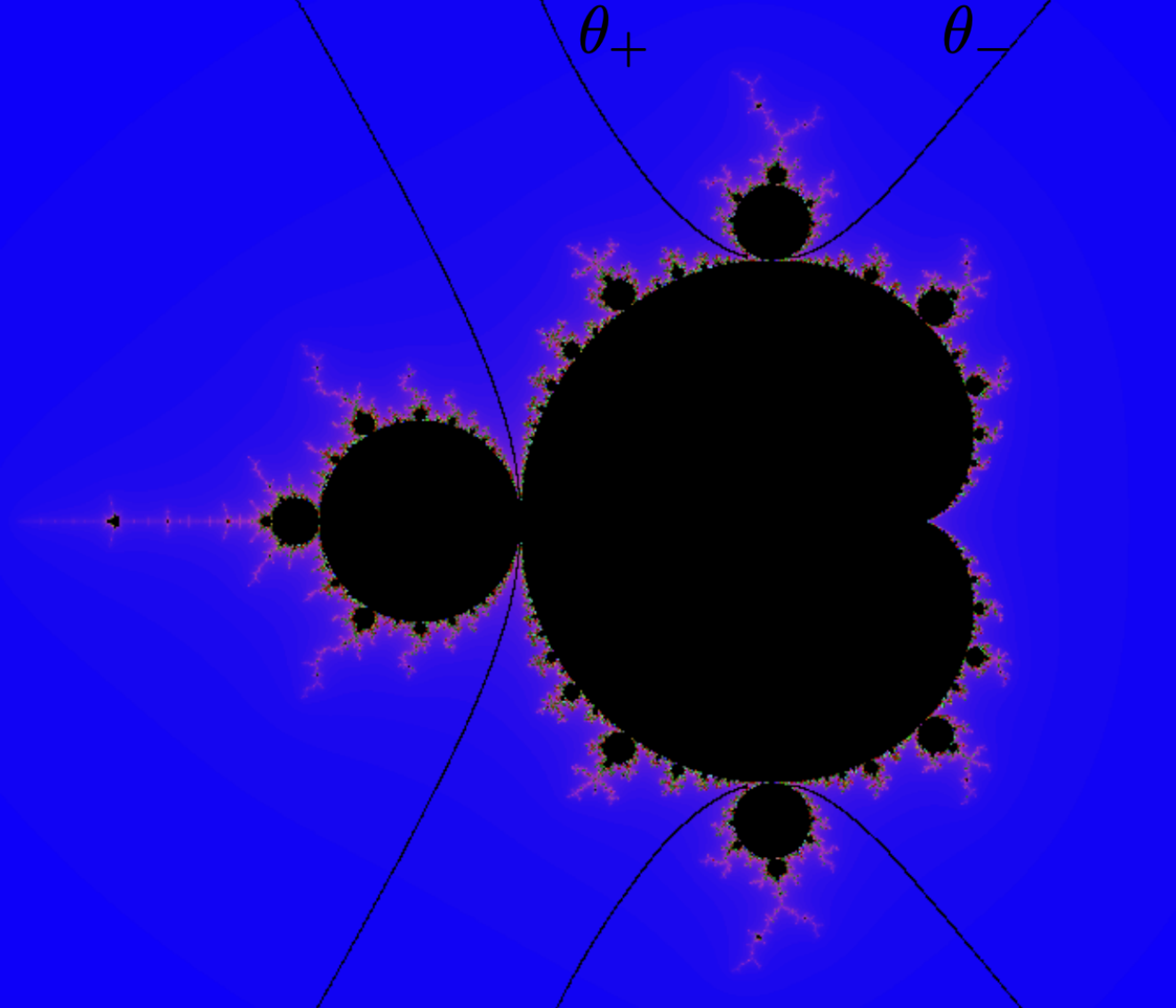}\ \includegraphics[height=70mm]{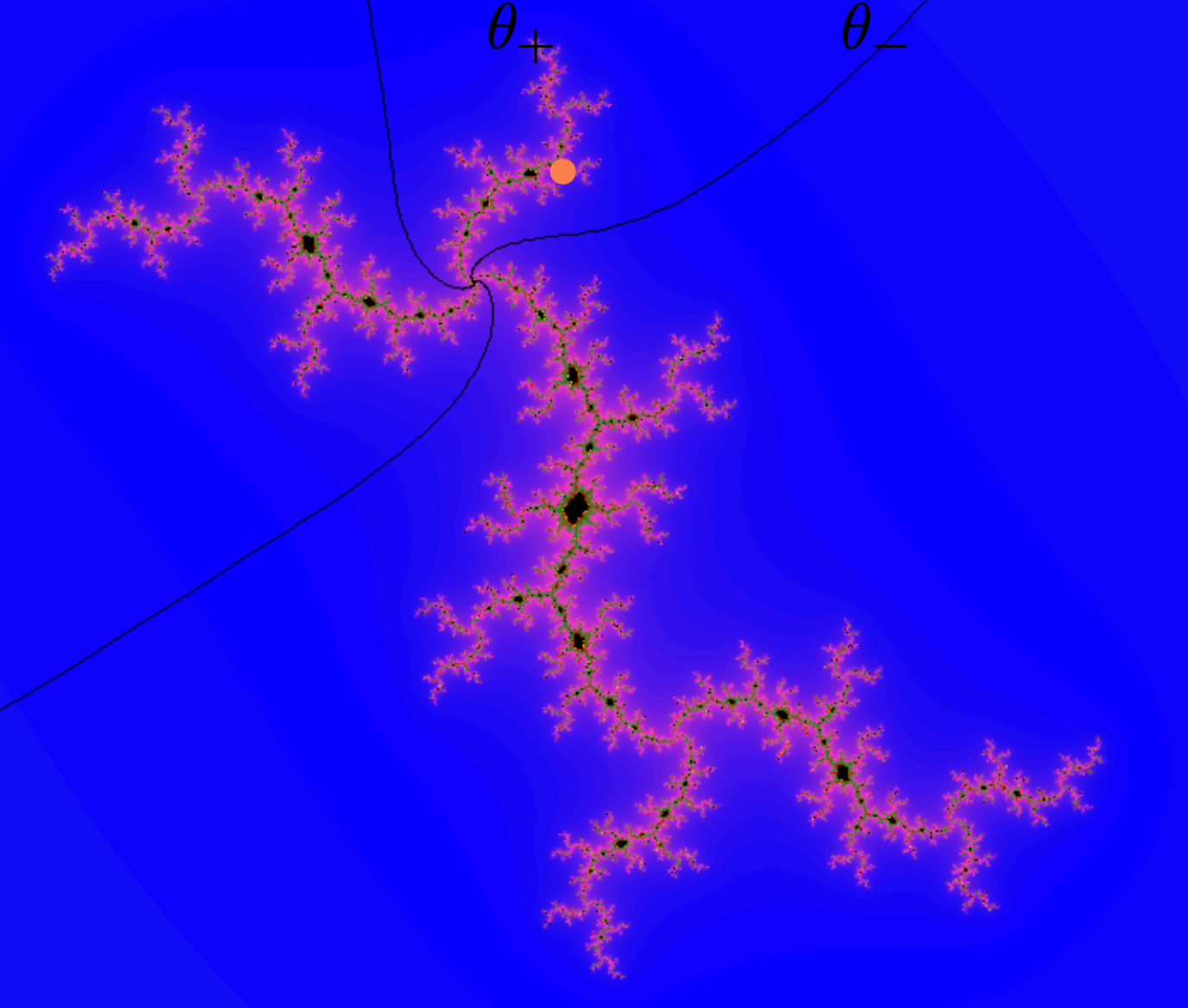}}
\caption{Parameter and dynamical Wakes}
\label{polynomial_wakes} 
\end{figure}

\begin{figure}[h]
\centerline{
\includegraphics[height=70mm]{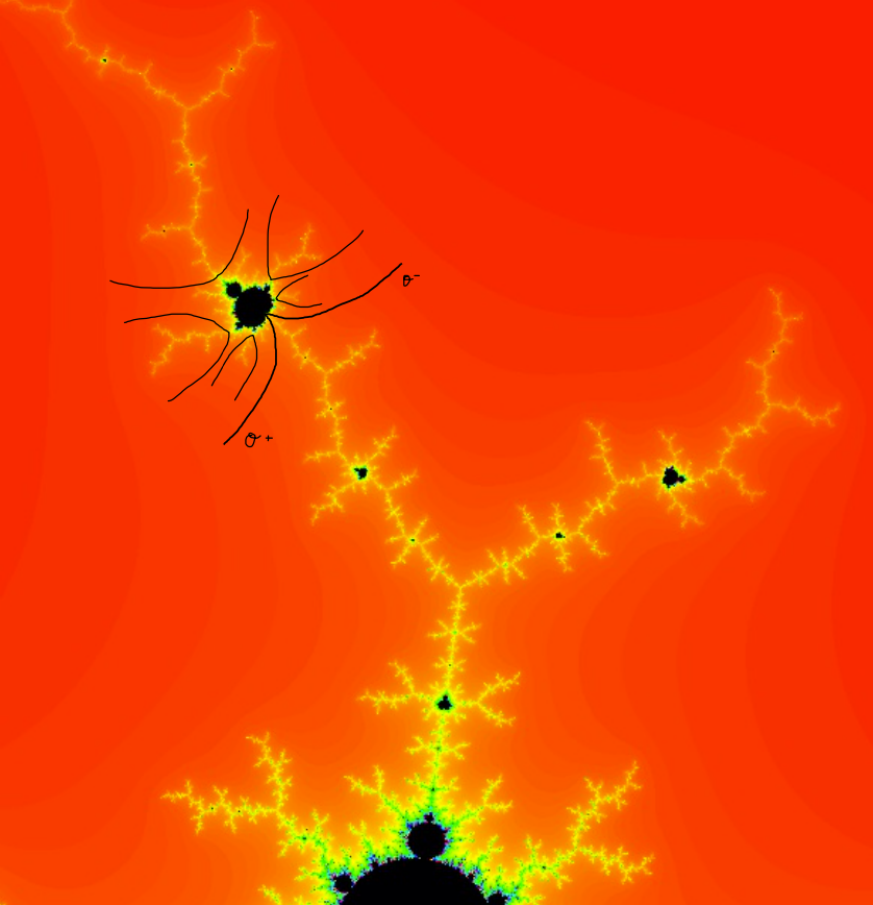}
\includegraphics[height=70mm]{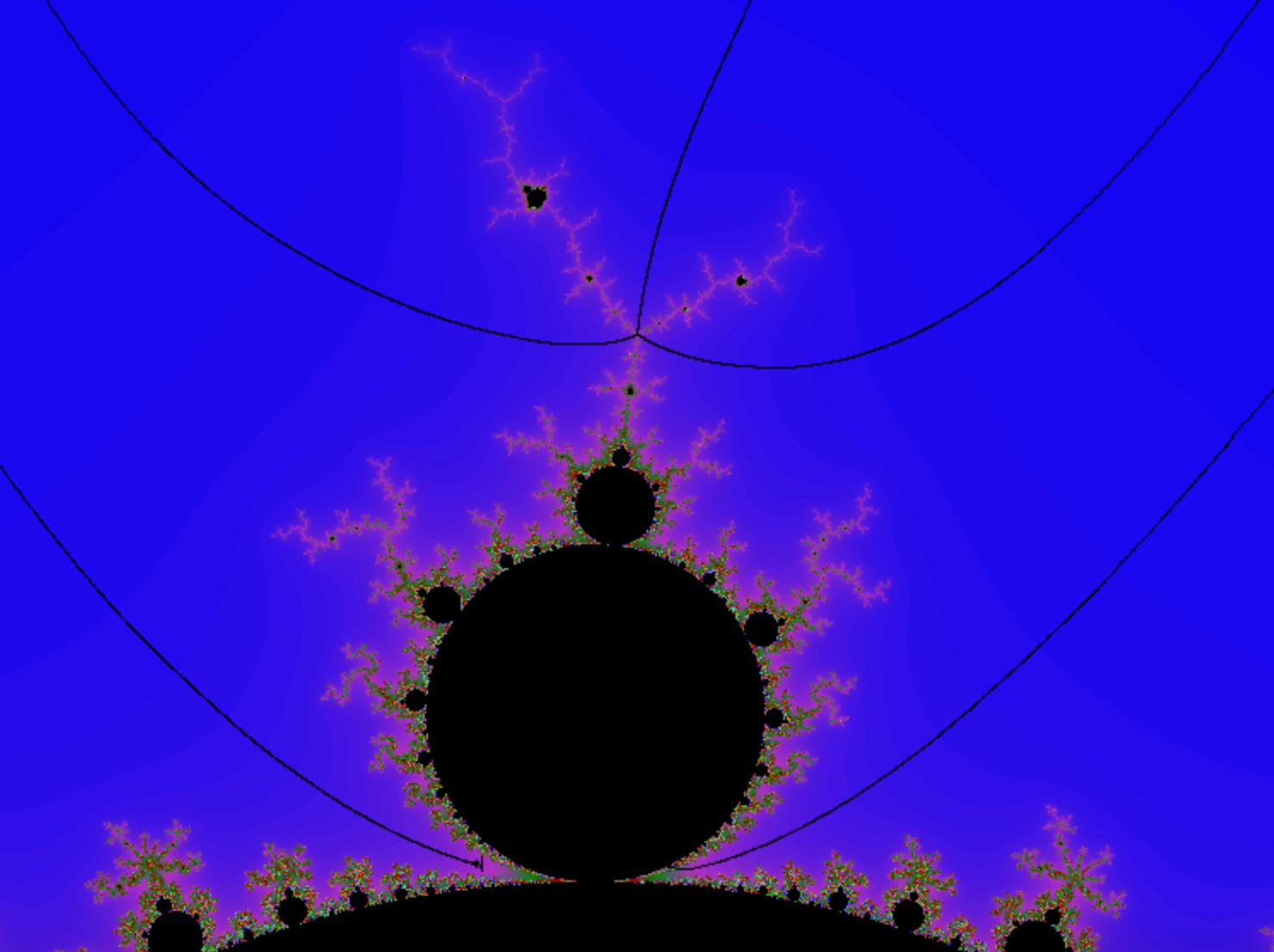}}
\caption{ Dyadic wakes of primitive and satellite Mandelbrot copy.}
\label{f:DiadicWake}
\end{figure}

For $c\in\Mbrot$, the map $\phi_c$ extends to an isomorphism 
between the basin of infinity and $\Chat\Sm\Dbar$, so that no ray bumps.
The two fixed points of $Q_c(z)$ are labelled $\beta_c$ and $\al_c$, with the convention that 
$\be_c$ is the landing point of the unique fixed ray, $\RR_0^c$. 
The other fixed point $\al_c$ can be attracting, neutral or  repelling. 
It is non repelling precisely when $c\in\overline{\Card}$.
Thus by \thmref{landingportrait} $\al_c$ is the landing point of $q>1$ external rays that define a cycle of 
 {\it combinatorial rotation number $p/q$}, and that thus assigns rotation number $p/q$ to $\al_c$. 
 This leads to the following stratification of $\Mbrot$, (see~\cite{Milnorasterisque}).
\begin{theorem}[Douady-Hubbard]~\\
$$
\displaystyle\Mbrot= \overline {\Card} \cup \bigcup_{\frac p q \not= \frac 0 1}L^\star_{p/q},
$$
where the uprooted limb  $L^\star_{p/q}$ consists of those parameters $c\in \Mbrot$
for which the separating fixed point $\al_c$ is repelling and has combinatorial rotation number $p/q$. See also \figref{f:M}, \figref{polynomial_wakes} and \figref{f:DiadicWake} for illustrations.
\end{theorem}

\subsection{Wakes} 
Since the {\Bottcher}-coordinate $\phi_c$ extends univalently to the set of points of potential greater than $G_c(0)$ 
it follows that {\mapfromto {\Phi_\Mbrot}{\C\Sm\Mbrot}{\C\Sm\Dbar}} 
given by $\Phi_\Mbrot(c) := \phi_c(c)$ is well defined for $c\in\C\Sm\Mbrot$. 
Douady and Hubbard proved that $\Phi_M$ is an isomorphism tangent to the identity at infinity.
Let $\Psi_\Mbrot$ denote its inverse. The parameter ray $\RR_\theta^\Mbrot$ of argument $\theta\in\TT$ 
is by definition $\RR_\theta^\Mbrot := \Psi_\Mbrot(\RR_\theta)$.

Rational parameter rays land at special parameters on the boundary of $\Mbrot$, 
see e.g. \cite{orsaynotesI} and \cite{orsaynotesII}.

\RREFTHM{Douady-Hubbard}{paramraylanding}
Let $\theta\in\TT$ be any {\rm(}pre-{\rm)}periodic argument under $\mtwo$,
i.e.~$\mtwo^k(\mtwo^l(\theta))=\mtwo^l(\theta)$, for some $l\ge 0$ $k\ge 1$.
Then $\RR^\Mbrot_\theta$ lands at a parameter $c\in\bd\Mbrot$. 

If  the argument $\theta$ is periodic {\rm(}i.e.~$l=0${\rm)}, then $Q_c$ admits a parabolic orbit 
of period $k'$ and multiplier $\e^{i2\pi p/q}$ such that $k'q = k$. 
Moreover the dynamical ray $\RR_\theta^c$ lands at a point in this parabolic orbit.
If furthermore $k>1$ then $c$ is the landing point of  $\RR^\Mbrot_{\theta}$ 
and of precisely one more parameter ray $\RR^\Mbrot_{\theta'}$. 
The angles  $\theta$ and $\theta'$ belong to the same cycle if $q>1$ --- 
the satellite case --- and to two different cycles if $q=1$ --- the primitive case. 
In either case the critical value $c$ is contained in the wake $\W_c(\theta,\theta')$, 
i.e. the domain bounded by the closure of the co-landing rays 
$\RR^c_\theta, \RR_{\theta'}^c$ and not containing $0$.

If $l>0$ then $\RR^\Mbrot_\theta$ lands at a Misiurewicz parameter $c$, 
such that $Q_c^l(c)$ belongs to a repelling cycle of exact period $k$. 
Moreover for any $\theta'\in\TT$ such that  the dynamical ray $\RR_{\theta'}^c$ lands on $c$, 
the parameter ray $\RR_{\theta'}^\Mbrot$ also lands on $c$.
\ENDTHM

Assume at first that $l=0$ and  $k>1$ (periodic arguments), then the parameter rays $\RR^\Mbrot_\theta, \RR^\Mbrot_{\theta'}$ co-land on a parabolic parameter $c$ as in the theorem above. 
Define the parameter wake $\W^\Mbrot(\theta,\theta')$ 
as the domain bounded by the closure of these rays and not containing the origin. 
It is easy to see that for any $c'\in\W^\Mbrot(\theta,\theta')$ 
the dynamical rays $\RR^{c'}_\theta, \RR_{\theta'}^{c'}$ move holomorphically with $c'$ 
and co-land on a repelling periodic point $w_{c'}$, 
which becomes the parabolic periodic point of $Q_c$ with the same rays landing, 
when $c'$ converges to the root $c$ of the wake. 
And that for $c'\notin\W^\Mbrot(\theta,\theta')$ the two rays $\RR^{c'}_\theta, \RR_{\theta'}^{c'}$ 
either land on different periodic orbits or at least one of them bump. 
For $c'\in\W^\Mbrot(\theta,\theta')$ we may thus also define the dynamical wake $\W_{c'}(\theta,\theta')$ 
by the same description. 

Without loss of generality we can assume $0<\theta<\theta'<1$. 
There exists $\whtheta, \whtheta'$, $0<\theta <\whtheta'<\whtheta<\theta'<1$ such that 
$\mtwo^k$ maps each of the intervals $[\theta, \whtheta']$ and $[\whtheta,\theta']$ 
diffeomorphically onto the full interval $[\theta,\theta']$. 
Let $C(\theta,\theta')$ denote the dyadic Cantor set consisting of those points which never escapes 
$[\theta, \whtheta']\cup[\whtheta,\theta']$ under iteration by $\mtwo^k$. 
Then $C(\theta,\theta')$ naturally corresponds to the classical middle third Cantor set.
As described by the Douady tuning algorithm the points of $C(\theta,\theta')$ 
are almost all of the arguments of external rays, which accumulates the filled-in Julia set $K_c'$ with 
$w_{c'}\in K_c'\subset \ov{\W_{c'}(\theta,\theta')}$ of a polynomial-like restriction of degree $2$ of $Q_c^k$ 
(see also \lemref{q_renormalization}) In the case of $q>2$ each gap of the Cantor set contains $q-2$ 
additional arguments of $K_{c'}$. 
It is a deep theorem that the set of such parameters for which $K'_c$ is connected, 
is a (derooted in the satelite case) copy $M^\Mbrot_{\theta,\theta'}$ of the Mandelbrot set. 
And the set of arguments for parameter rays accumulating $M^\Mbrot_{\theta,\theta'}$ contains $C(\theta,\theta')$.

\label{Renormalization_Cantor_set_of_arguments}
The gaps of the Cantor set $C(\theta,\theta')$ naturally corresponds to the dyadic numbers 
$r/2^s$, $0<r,s$, $r$ odd and $r<2^s$ with the initial gap $]\whtheta',\whtheta[$ corresponding to $1/2$. 
The endpoints of the gap corresponding to $r/2^s$ map under $\mtwo^{ks}$ to the $\mtwo^k$ fixed points 
$\theta, \theta'$ and so the corresponding parameter rays co-land on a Misiurewicz parameter 
$c'' = c(\theta,\theta',r,s)$ such that $Q_c^{sk}(c'') = w_{c''}$. 
We denote by $\W^\Mbrot(\theta,\theta',r,s)$ the wake bounded by the closure of these rays. 
For $c'\in\W^\Mbrot(\theta,\theta',r,s)$ the corresponding dynamical rays co-land on a 
pre-image of $w_{c'}$ under $Q_{c'}^{ks}$ and so bound a corresponding dynamical wake 
$\W_c(\theta,\theta',r,s)$. 
We define the $r/2^s$-dyadic limb of $M^\Mbrot_{\theta,\theta'}$ as the intersection 
$$
L^\Mbrot(\theta,\theta',r,s) := \W_c(\theta,\theta',r,s)\cap\Mbrot.
$$
For a similar discussion see \cite[Section 1]{PRcarrots}.

In the special case where the $k$-periodic parameter rays $\RR^\Mbrot_\theta, \RR^\Mbrot_{\theta'}$ 
co-land on a parabolic parameter $c$ for which the parabolic periodic point of period $k'=1$, i.e. equals $\al_c$, so that $\theta,\theta'$ belong to a $p/q$-cycle with $q=k$, we use the standard short hand $\W^\Mbrot(p/q)$ for the parameter wake $\W^\Mbrot(\theta,\theta')$, 
$\W_{c'}(p/q)$ for the dynamical wake $\W_{c'}(\theta,\theta')$ when $c'\in\W^\Mbrot(p/q)$, $\Mbrot_{p/q}$ for 
$\Mbrot_{\theta,\theta'}$, $\W^\Mbrot(p/q,r,s)$ for the dyadic parameter wake $\W^\Mbrot(\theta,\theta',r,s)$ 
and $\W_{c'}(p/q,r,s)$ for corresponding dyadic dynamical wakes $\W_{c'}(\theta,\theta',r,s)$.

Note that by definition the uprooted limb $L^\star_{p/q} = \Mbrot\cap\W^\Mbrot(p/q)$.

\subsection{Basic notation for maps in $Per_1(1)$.} 
\begin{figure}[h]\label{f:Moneparam}
\centerline{\includegraphics[height=70mm]{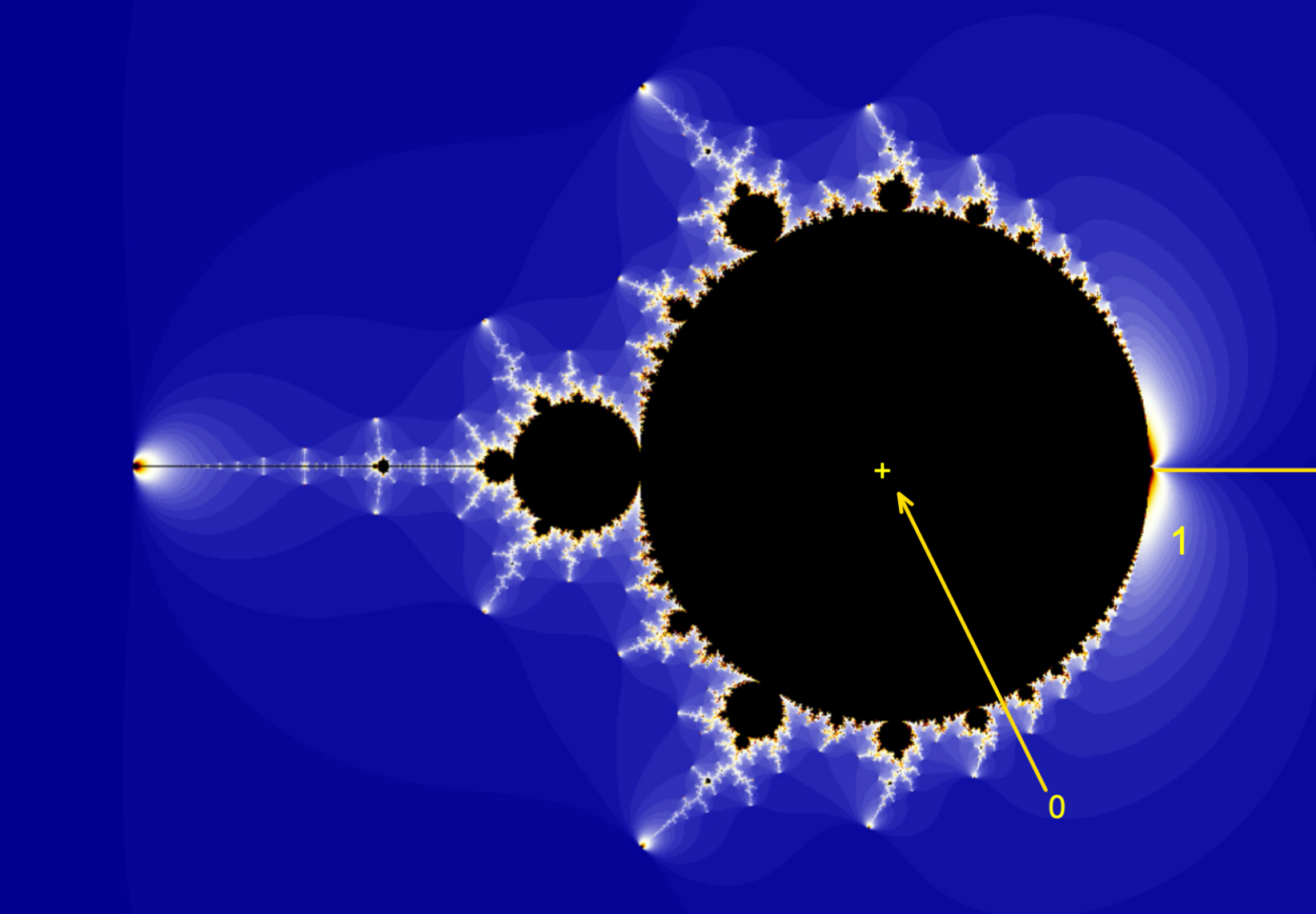}}
\caption{The $A$ plane ($A=1-B^2$ is the multiplier of the $\al$-fixed point of $g_B$)}
\label{f:Moneparam}
\end{figure}

The slice $\Peroo$ does not admit a  normal form which univalently parameterizes it. 
As a consequence there is not a universal choice of parametrization. 
We shall henceforth use several parametrizations interchangeably. 
First of all we shall write $[g]$ for the element in $Per_1(1)$ represented by $g$. 
Such a map $g$ has a parabolic fixed point of multiplier $1$ and one more fixed point 
of multiplier $A\in\C$. This fixed point coincides with the two others precisely when $A=1$. 
Thus $\sigma_1([g])=A$, because the parabolic fixed point of multiplier $1$  is multiple.
We invite the reader to think of $g$ as taking the form $g(z) = g_B(z) = z + B + 1/z$ 
for $B\in\C$, in which case the fixed point of multiplier $1$ is at infinity 
and the critical points are located at $\pm 1$ and the corresponding critical values are 
$B\pm2$.For $B=0$ the three fixed points coincide at $\infty$ and otherwise $g_B$ has 
a finite fixed point at $-1/B$ with multiplier $A = 1 - B^2 = \sigma_1([g_B])$.
As a consequence we have
\REFREM{peroneparametrizations}The correspondance $B\in \C\mapsto [g_B]\in \Peroo$ is  a $2$ to $1$ branched covering.\\
 \noindent We shall use interchangeably the notations $g=g_B$, $[g]$, $B$ and $A$. 
In particular we shall use $g\in\Mone$, $[g]\in\Mone$, $A\in\Mone$ and $B\in\Mone$ in the obvious meaning see also \figref{f:Moneparam}. 

However we shall mainly be interested in $A\in\C\Sm[1,\infty[$ which biholomorphically 
corresponds to $B\in\Hplus := \{x+iy| x>0\}$. 
For this reason our preferred representation of $Per_1(1)$ shall be via the parameter 
$B \in \Hplus\cup\{0\}$ or $B\in\Hpb := \{B' \mid \Re(B')\geq 0\}$. 
\ENDREM

 \noindent   
For a map $g = g_B$  denote by $\Lambda_B$  the parabolic basin {\it i.e.} the maximal open subset of points converging to the parabolic fixed point $\infty$ of multiplier $1$. 
It is completely invariant~: $g^{-1}(\Lambda_B)=\Lambda_B$. 
If $B=0$, the parabolic basin $\Lambda_0$ of $g_0$ has two connected components $\pm  \Hpb$, each  containing a critical point. Indeed,   the Julia set $J(g_0) $ is simply the imaginary axis $  i\R\cup\{\infty\}$. 

For $B\not=0$ the parabolic basin $\Lambda_B$ is connected.  
Denote by filled Julia set its complement $K(g_B)  = K_B := \C\setminus\Lambda_B$. 
The Julia set $J(g_B) = J_B$ is then the  common boundary 
$J_B=\partial K_B=\partial \Lambda_B$. 
The Julia set  is either connected and $B\in \Mone$ or it is a Cantor set (see~\cite{Milnorquad}). In the  Cantor  case $\Lambda_B$ is connected, infinitely connected, contains both critical points and the dynamics on the Julia set  is conjugate to the one-side shift map on two symbols. 

\begin{figure}
\includegraphics[height= 2.73in]{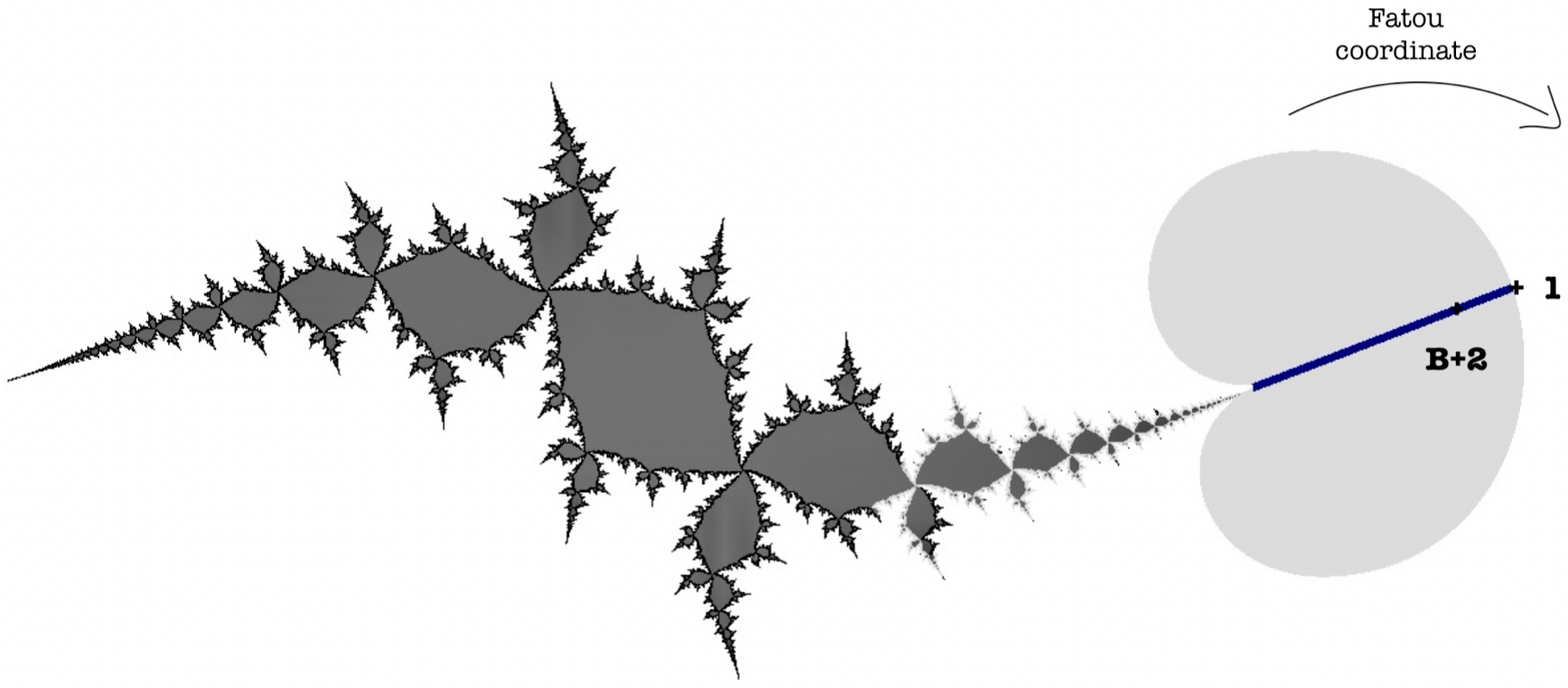}\qquad 
\includegraphics[height=2in]{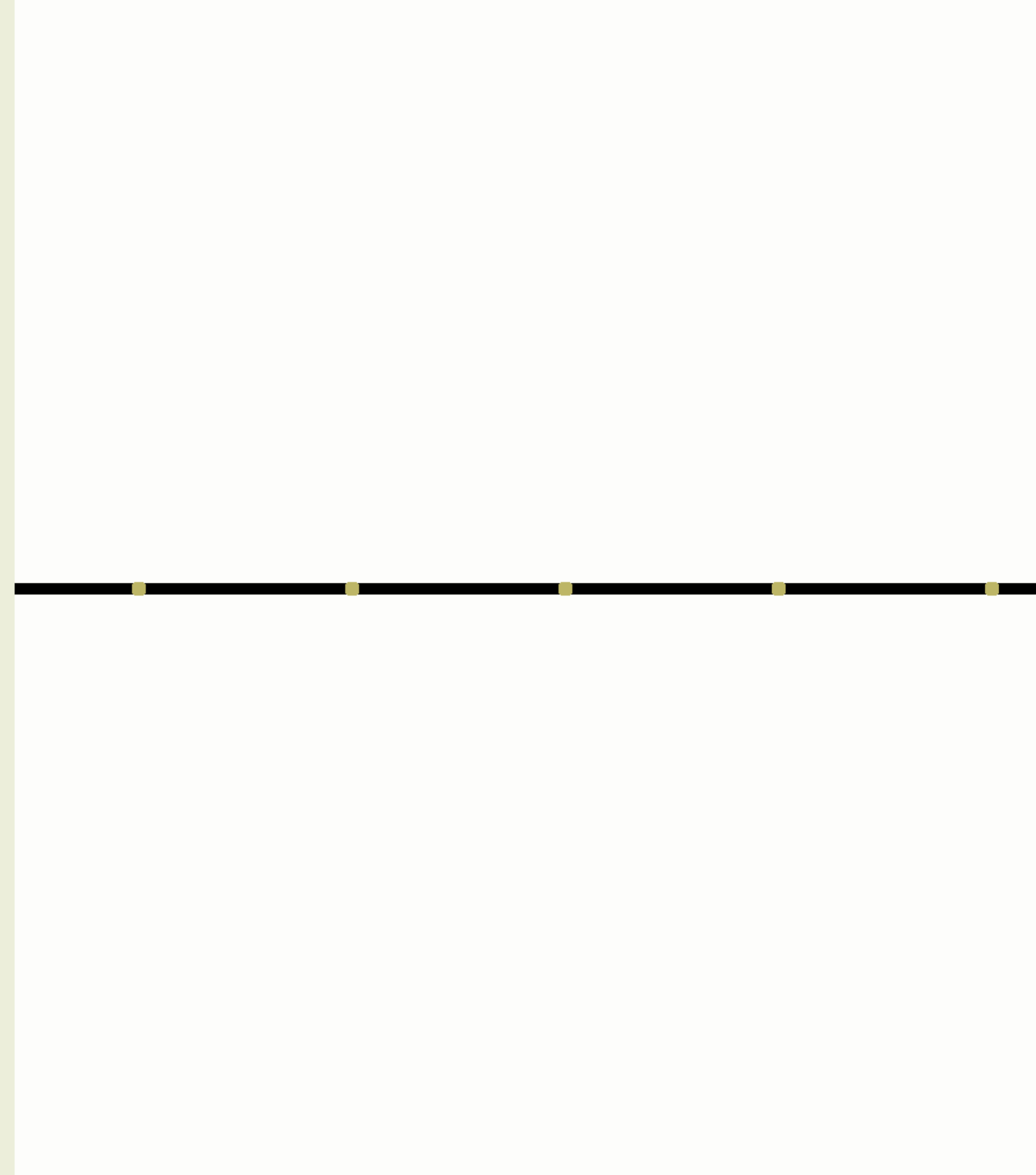}
\caption{Fatou coordinate $\phi_B$ and $\Omega^B$}
\label{Fatou-coordinate-illustration}
\end{figure}

For $B\not=0$, $g_B$ admits an attracting Fatou coordinate 
{\mapfromto {\phi_B} {\Lambda_B}\C}, 
which semi-conjugates $g_B$ to translation by $1$. 
The map $\phi_B$ is unique up to post-composition by a translation. 
See \figref{Fatou-coordinate-illustration} for an illustration.

\REFDEF{Omega-B}
There is a unique maximal forward invariant domain $\Omega^B\subset \Lambda_B$, 
whose boundary contains at least one critical point 
and which is mapped univalentely onto a right half plane by $\phi_B$. 
Adjusting $\phi_B$ we can assume the half plane is $\Hplus$ 
$\phi_B$ sends a critical point to $0$. 
This critical point is denoted the \emph{fastest critical point}.
\ENDDEF

For $B=0$ there are two such coordinates one on each connected component 
of the basin $\Lambda_{0}$.

In the next section we shall describe 
another representative $\Bla$ of $[g_0]$, which is more suitable for comparison 
with the polynomial $z^2$.

\REFLEM{premier}For $B\in\Hplus$ the critical point $1$ is the fastest escaping critical point. 
\ENDLEM
\proof We prove that the critical points escape at the same rate if and only if the parameter belongs to the imaginary axis.  
Thus in $\H_+$ either $+1$ is always the fastest escaping critical point or $-1$ is. 
Note that for $B=1$ we have $g_1(-1) = -1$ so that the critical point $-1$ is fixed and thus $+1$ is 
the only escaping and hence fastest escaping critical point.

Now consider first a map   $g_{ir}$, whith $r\in\R$.  It commutes with the reflection in the imaginary axis, $z\mapsto -\overline{z}$.
Hence, the critical points $\pm 1$ escape at equal rates for any $r$ 
(for $r=0$ the critical points however are in distinct components of the parabolic basin).

Suppose that for parameters  $B\neq B'$, $g_B$ and $g_{B'}$ both have two first attracted critical points, i.e.  critical points $c_1$, $c_2$ for $g_B$ and  $c_1'$, $c_2'$ for $g_{B'}$ 
respectively satisfy  for  well  chosen Fatou coordinates $\phi_B, \phi_{B'}$ that 
$\phi_B(c_1) = \phi_{B'}(c_1') = 0$ and $\phi_B(c_2) = \phi_{B'}(c_2') = iy$ for some non-zero real $y$. Then $\eta := \phi_{B'}^{-1}\circ \phi_B$ defines a biholomorphic conjugacy on a petal $\Omega$ --- which is  mapped biholomophically to the right half plane by $\phi_B$. Moreover $\eta$ maps critical values to critical values and so extends as a biholomorphic conjugacy between the parabolic basins by iterated lifting. Then it extends to a global topological conjugacy, since Julia sets are Cantor sets. Finally this conjugacy is holomorphic between $g_B$ and $g_{B’}$, because the Julia set is holomorphically removable. The map $\eta$  fixes $\infty$ and hence $0$, pole of the maps,  so $\eta$ is linear and since it maps critical points to critical points  ($\pm 1$) it is either the identity or $z\mapsto -z$.  Hence $B=\pm B'$. Applying, this to $B\in i\R$ and $B'\in\C$ we get that $B'\in i\R$. 
\endproof

\REM
For $B\in\Hpb$ we normalize the Fatou-coordinate $\phi_B$ by $\phi_B(1) = 0$. 
By \lemref{premier} there is no option for $B\in\Hplus$. 
But for $B \in i\R = \partial\Hplus$, this is a choice. 
For $B = 0$ this corresponds to considering 
$\ov{\H}_-$ as $K_0$ and $\La_0$ as $\Hplus$. In any case
the restriction 
$$
{\mapfromto {\phi_B}{\Omega^B}{\Hplus}}
$$
is a biholomorphic conjugacy and $\Om^B$ depends continuously on $B\in\Hpb$.

For $B\in\Hpb$ the critical value $B+2 = g_B(1) \in \Om^B$, this implicitly defines $\Om^0$. 
We denote by $v_B = -2+B = g_B(-1)$ the other critical value. 
For $ B\in \Mone$,  the points $-1$ and $v_B$ both belong to $K_B$, 
$\Lambda_B$ is isomorphic to $\D$ and contains only the critical point $1$. 
\ENDREM

In case $B\in i\R$ both critical points are on the boundary of $\Omega_B$ 
and we could have chosen to normalise $\phi_B$ using the critical point $-1$ in this case.  
The two choices for $\phi_B$ thus differ by a purely imaginary translation and 
the values of $\phi_B(c_1)$ for the two different choices
are purely imaginary and complex conjugate. 
This fact is used by Shishikura, when constructing a natural isomorphism between 
$Per_1(1)\Sm\Mone$ and $\Chat\Sm\Dbar$, see e.g. \cite{Milnorbicrit}. 
We shall in order to ease the notation not use this isomorphism here, 
but stop two steps before the end of the construction. 
This is the content of the following subsection.

 \subsection{Parametrization of $\C\setminus \Mone$}\label{s:parametrization}
\begin{figure}[h!]
\begin{center}
\includegraphics[height=1.52 in]{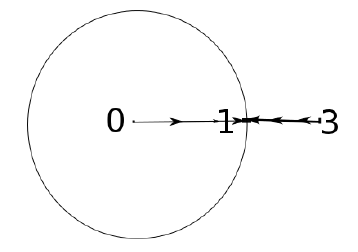}\end{center}
\caption{The attracting directions of the Blaschke product $\Bla$}
\label{the Blaschke-product-Bl}
\end{figure}

The  idea of Shishikura's proof is to parametrize $\C\setminus \Mone$ by the relative position of the critical values in suitable coordinates, namely in $\C\setminus\overline {\D}$ viewed 
as the parabolic  basin of the  standard parabolic  Blaschke product.  
For this purpose we introduce the  parabolic Blaschke product $\Bla \in [g_0]$ 
see \figref{the Blaschke-product-Bl} for an illustration. 
It acts as the external class of the maps $g_B$, similarly to $z^2$ for quadratic polynomials 
$$\Bla(z) = \frac{z^2+1/3}{1+z^2/3}.$$

Note that  $\Bla = \nu\circ g_0\circ \nu^{-1}$, 
where $\nu$ is the \Mobius~transformation $\displaystyle \nu(z) = (z+1)/(z-1)$. 
It follows immediately from the above that the Julia set $J(\Bla)$ is the unit circle,  
$\D$ and $\Chat\setminus\Dbar$ are the two components of the basin of the double parabolic fixed point $1$. 
The critical points are $0$ and $\infty$, with images $1/3$ and $3$ respectively. 
The map $\Bla$ admits  $\tau(z)= 1/\overline z$ as a symmetry interchanging the immediate basins.   
The arcs $[0,1]$ and $[\infty,1]$ form attracting axis for the attracting petals of $\Bla$.
 \begin{figure}[h]
\begin{center}
\includegraphics[height=6 cm]{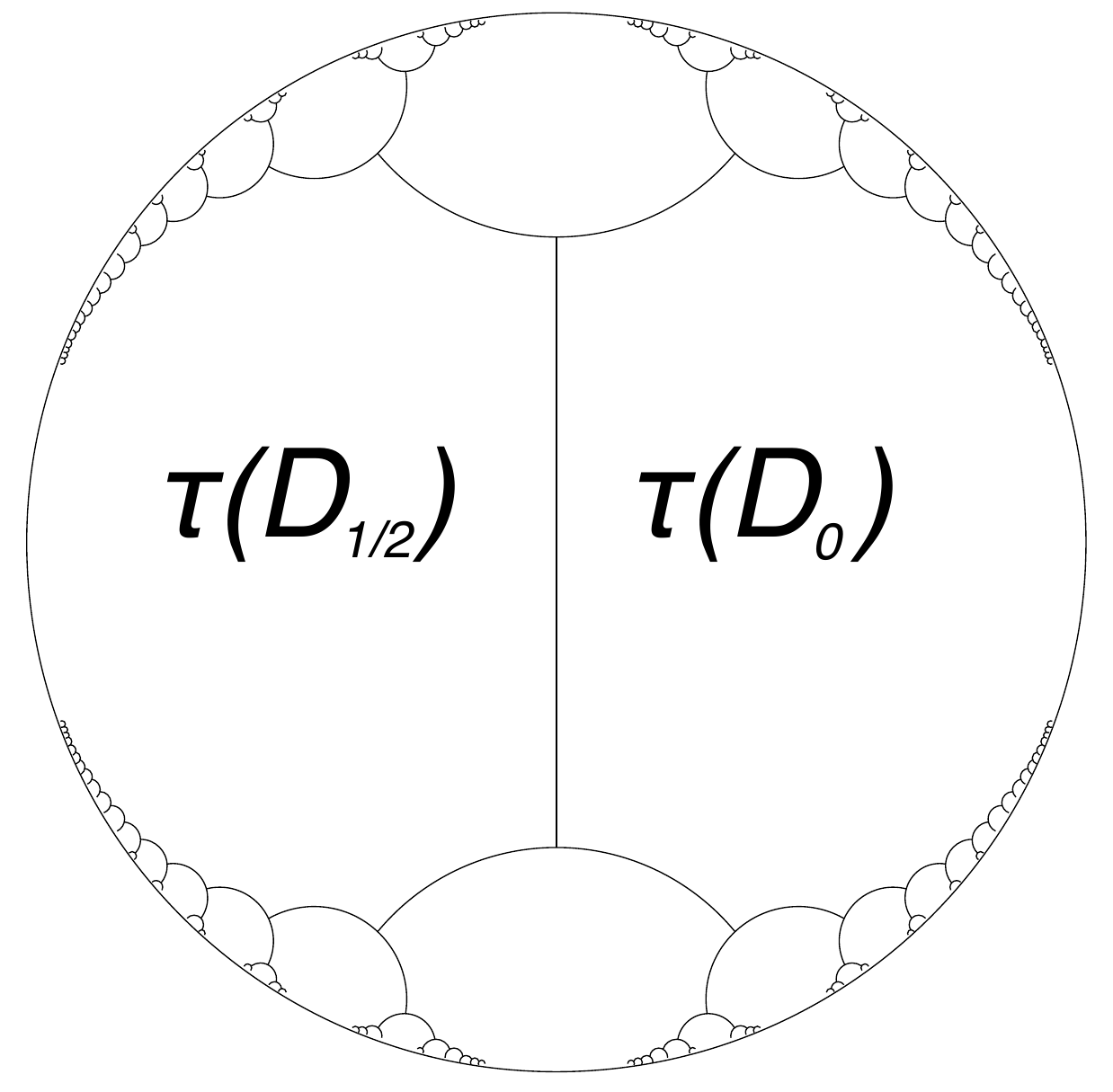}
\end{center}
\caption{The tiling of the internal parabolic bassin $\Bla$ by the connected components of $\phi^{-1}(\C\Sm]-\infty,0])$.\label{paratiles}}
\end{figure}

We denote by $\phi = \phi_0\circ \nu^{-1}:\Chat\setminus \Dbar\to\C$ 
the Fatou coordinate normalized as above, that is 
$\phi(\infty)=0$ and $\phi$ is univalent on $\Omega =\phi^{-1}( \H_+)$, 
where $\Bla(\Omega) \subset\Omega$. 
Let $D_0$ denote the connected component of $\phi^{-1}(\C\Sm]-\infty, 0])$.
Then $\phi^{-1}$ extends to a univalent map 
{\mapfromto {\phi^{-1}}{\C\Sm]-\infty, 0]} {D_0\supset\Omega}}, 
because $\Chat\Sm\Dbar$ contains only the critical point $\infty$ of $\Bla$. 
We write $D_{1/2}$ for the disk $-D_0$ and $D_1$ 
for the disk, which is the interior of the closure of $D_0\cup D_{1/2}$. 
Then $D_1= \Bla^{-1}(D_0)$.

We define topological disks $\Omega_n = \Bla^{-n}(\Omega)$ and 
$D_n=\Bla^{-n}(D_0)\supset\Om_n$ for each $n\geq 0$.
\begin{figure}[h!]
\begin{center}
\includegraphics[height=1.52 in]{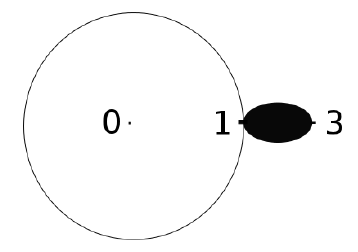}\end{center}
\caption{The image $\Bla(\Omega)$ of the petal $\Omega$ in black.}\label{Petalext}
\end{figure}

The biholomorphic parametrization of $\C\setminus \Mone$ easily follows from the following construction. Let $B\in\Hpb\Sm\Mone$ so that both critical points $\pm1$ of $g_B$ 
belong to $\Lambda_B$. 
For each $k\geq 0$ let $\Omega^B_k = g_B^{-k}(\Omega^B)$ and let
$n\geq 0$ be minimal with $v_B=g_B(-1)\in\Omega^B_n$. 
Then each $\Omega^B_k$ with $k\leq n$ is simply connected.
Define a biholomorphic conjugacy  $h_B$ by  
$$\begin{array}{cccc}h_B: &\Omega^B&\to& \Omega\\
&z&\mapsto &\phi^{-1}\circ\phi_B(z).
\end{array}$$
Then since $h_B$ sends the critical value $g_B(1)$ to the critical value $3=\Bla(\infty)$, 
and the domians $\Omega^B_k$ are simply connected for $k\leq n$ 
the map $h_B$ can be univalently lifted iteratively to define a conjugacy between $g_B$ 
and $\Bla$ on the domain $\Omega^B_n=g_B^{-n}(\Omega^B)$ 
containing the second critical value $v_B$, but not the second critical point $-1$.
Then 
 \[\xymatrix{\Omega^B_n\ar[r]^{g_B}\ar[d]_{h_B}&\Omega^B_{n-1}\ar[r]^{g_B}\ar[d]_{h_B}&\Omega^B_{n-2}\ar[d]_{h_B} \ar[r]^{g_B}&\ldots \ar[r]^{g_B}\ar[d]^{ \ldots}\ar[d]_{ \ldots}&\Omega^B_1\ar[d]^{h_B}  \ar[r]^{g_B} &\Omega^B\ar[d]^{h_B} \\
\Omega_{n}\ar[r]_{\Bla}&\Omega_{n-1}\ar[r]_{\Bla}&\Omega_{n-2} \ar[r]_{\Bla}&\ldots \ar[r]_{\Bla}&\Omega_{1}\ar[r]_{\Bla}&\Omega  }\]
Hence we get the following Lemma. 
\REFLEM{conj} 
For every $B\in\Hpb\Sm\Mone$ there exist $n = n_B \in \N$ 
and $h_B :  \Omega_n^B\to \Chat\Sm\Dbar$ a univalent conjugacy between $g_B$ and $\Bla$ such that $h_B(B+2)=3$ and $v_B = B-2 \in\Omega^B_n$.
\ENDLEM

\REM
Note that the map $(B,z) \mapsto h_B^{-1}(z)$ is complex analytic as a function of the 
pair of variables $(B, z)$. Because the Fatou-coordinates depend holomorphically on $B$ 
and the map $(B,z)\mapsto (B,h_B(z))$ is locally biholomorphic off the critical points of $h_B$.
\ENDREM
\REFDEF{Complcoord}
Let $\Omega' = \Bla(\Omega)$ and define a holomorphic map
$${\mapfromto \Upsilon {\Hplus\Sm\Mone} {\C\Sm(\Dbar\cup\overline{\Omega'})}} 
\qquad\textrm{by}\qquad
\Upsilon(B) := h_B(v_B).
$$
\ENDDEF
We shall see that this map is injective and analytically extends as homeomorphism 
$$
\mapfromto\Upsilon {\Hpb\Sm\Mone}{\C\Sm(\Dbar\cup\Omega')}. 
$$
\REM
Since $\si_1([g_B]) = 1 - B^2 = A$ the representation of the above map in the $A$-coordinate 
gives a holomorphic and in fact biholomorphic map (see also \remref{peroneparametrizations})
$${\mapfromto {\widehat\Upsilon} {\C\Sm(\Mone\cup[1,\infty[)} {\C\Sm(\Dbar\cup\overline{\Omega'})}}.
$$
\ENDREM
\REFPROP{milnor}
The map $\Upsilon$ is a proper holomorphic map of degree $1$, hence an isomorphism. 
\ENDPROP
\PROOF
We shall show that $B\to \partial (\Hplus\Sm\Mone)$ implies 
$\Upsilon(B)\to\partial(\Dbar\cup\overline{\Omega'})$. 
It follows that $\Upsilon$ is proper.

For $B\not= 0$ the linear map $z\mapsto z/B$ conjugates the map $g_B$ 
to $z\mapsto z + 1 + 1/(B^2z)$ with critical points at $\pm1/B$ and corresponding critical values $1\pm2/B$. It follows that both critical values for $g_B$ belong to $\Omega^B$, 
when $|B|$ is sufficiently large and that $\Upsilon(B) = h_B(v_B)$ converge to $h_B(B+2) = 3$, 
when $|B|\to\infty$. 
If $\Re(B)\to 0$, then $\Upsilon(B) = h_B(v_B)$ converge to $\partial \Omega'$, 
since $B\mapsto h_B(v_B)$ is continuous and belongs to $\partial\Omega'$, 
when $\Re(B) = 0$.
Finally suppose $\{B_k\}_k\in\Hplus\Sm\Mone$ is a sequence converging to $\partial\Mone$, 
but $h_{B_k}(v_{B_k})$ does not converge to $\partial\D$. 
Then passing to a subsequence if necessary we can suppose that $B_k\to B\in\partial\Mone$ and $h_{B_k}(v_{B_k})$ converge to $w\in\Chat\Sm(\Dbar\cup\Omega')$. 
Choose $N$ such that $w\in\Omega_N$ and thus $\Bla^N(w)\in\Omega$. 
Then $g_{B_k}^N(v_{B_k})\in\Omega^{B_k}$ for all $k$ large enough. But then by continuity 
also $g_B^N(v_B)\in\Omega^B$, contradicting that $B\in\partial\Mone$. 

Finally the degree is $1$ because it extends continuously and injectively to $\partial\Hplus$, 
which is mapped onto $\partial\Omega'$ and we showed above that 
if $B_n\to\partial\Mone$ then $\Upsilon(B_n)$ converge to $\partial\D$.
\ENDPROOF

\REFLEM{D1_ext_of_hBinv}
If $\Upsilon(B)\notin D_0$ then $h_B^{-1}$ extends as a biholomorphic conjugacy 
$$
\mapfromto {h_B^{-1}} {D_1} {h_B^{-1}(D_1)}. 
$$
\ENDLEM
\DEF 
In view of the Lemma above, when $\Upsilon(B)\notin D_0$
we define $D_1^B := h_B^{-1}(D_1)$, $D_{1/2}^B := h_B^{-1}(D_{1/2})$ and  $D_0^B := h_B^{-1}(D_0)$.
\ENDDEF
\PROOF
Note that $\bigcup_m \Omega_m\cap D_0 = D_0$,
$\bigcup_m \Omega_m\cap D_1 = D_1$ 
and that each set $\Omega_n\cap D_0$ is simply connected and does not contain 
$h_B(v_B)$. 
Hence we obtain an increasing sequence of extensions of $h_B^{-1}$ 
by iterated lifting:
 \[\xymatrix{\Omega_n\cap D_0\ar[r]^{\Bla}\ar[d]_{h_B^{-1}}&\Omega_{n-1}\cap D_0
 \ar[r]^{\Bla}\ar[d]_{h_B^{-1}}&\ \ \Omega_{n-2}\cap D_0\ \ar[d]_{h_B^{-1}} \ar[r]^{\ \ \ \ \Bla}&\ldots \ar[r]^{\Bla\ \ }\ar[d]^{ \ldots}\ar[d]_{ \ldots}
 &\Omega_1\cap D_0\ar[d]^{h_B^{-1}}  \ar[r]^{\ \ \Bla} &\Omega\ar[d]^{h_B^{-1}} & \\
h_B^{-1}(\Omega^B_n\cap D^B_0)\ar[r]_{g_B}&h_B^{-1}(\Omega^B_{n-1}\cap D^B_0)\ar[r]_{g_B}&h_B^{-1}(\Omega^B_{n-2}\cap D^B_0) \ar[r]_{\ \  \ \ \ \ \ g_B}&\ldots \ar[r]_{g_B   \ \  \ \ \ \ }&
h_B^{-1}(\Omega^B_{1}\cap D^B_0)\ar[r]_{\ \ \ \ \ g_B}&\Omega^B  
}\]
\ENDPROOF

%%%%%%%%%%%%%%%%%%%%%%%%%%%%%%%%%%%%%%%%%%%%%%%
\section{Parabolic Rays}\label{s:parabolicrays}
%%%%%%%%%%%%%%%%%%%%%%%%%%%%%%%%%%%%%%%%%%%%%%%
\subsection{Dynamical parabolic rays}
We first define parabolic rays in $\D$ and in $\C\setminus \overline \D$ for the model map $\Bla$.
 We then define parabolic  rays in the basin of $\infty$ for the maps $g_B$. 
In  the case where  the Julia set is connected, the conjugacy $h_B$ 
(between $g_B$ and $\Bla$) extends to the whole basin of infinity, 
so we just pull back the parabolic rays defined for the model map  $\Bla$.
In the non-connected case, we pull back when it is possible the beginning of the ray. 

%%%%%%%%%%%%%%%%%%%%%%%%%%%%%%%%%%%%%%%%%%%%
\subsubsection{Parabolic ray for the Blaschke product}\label{s:blaschke}

The notion of (external)  rays is well defined for quadratic polynomials, since on their basin of $\infty$  polynomials  are conjugated (in the connected case) 
 to $z\mapsto z^2$ on  $\C\setminus \overline {\D}$. The (external) rays are the  pull-back of  straight lines in  $\C\setminus \overline {\D}$.
 
The map $\Bla$ is a degree $2$ map on $\D$ and on $\widehat {\C}\setminus \D$, but it is not 
conjugate to $z\mapsto z^2$ on these domains. 
Nevertheless,  $\Bla$ is conjugate to $z^2$ on $\mathbb S^1$. 

\begin{lemma}\label{l:h}There exists a unique homeomorphism
{\mapfromto {h} \Sen \Sen} fixing $1$
and conjugating $z\mapsto z^2$ to $\Bla$, i.e.~$h(z^2)=\Bla\circ h$.  It commutes with $z\mapsto \overline z$.
\end{lemma}
\proof
Indeed, the map $\Bla$ is weakly expanding on $\Sen$,
as $|\Bla'(z)|\geq 1$ on $\Sen$ with equality iff $z^2=1$.
The rest of the proof  is a classical theorem for strongly expanding maps, for which the proof 
passes over to the weakly expanding case with out any essential changes.
(Define recursively 
{\mapfromto {h_n} \Sen \Sen} by $\Bla\circ h_n = h_{n-1}(z^2)$ and $h_n(1) = 1$, with $h_0=id$. 
The maps $h_n$ converge to  an order preserving bijection between the two sets 
of iterated preimages of $1$ and by the weakly expanding property both sets of iterated preimages 
are dense in $\Sen$ so that the limit of the $h_n$ exists on all of $\Sen$ and 
is the required topological conjugacy.)
\endproof

Similarly to the binary expansion of the angle, 
we will define rays for $\Bla$ using the itineraries. 
   
Let $\Sigma_2:={\{0,1\}}^{\N}$ denote the one-sided shift space
on $2$-symbols.  The angle $\theta$ is said to have  binary expansion $\veps=(\epsilon_1,\epsilon_2,\ldots,\epsilon_n,\ldots)$ if   
$$\bftheta = \sum_{n=1}^\infty \frac{\epsilon_n}{2^n}.
$$
Denote by {\mapfromto {\Pi_2} {\Sigma_2} \Sen} the 
projection map:
$
\Pi_2(\epsilon_1,\epsilon_2,\ldots,\epsilon_n,\ldots) =
\exp\left(2\pi i \theta \right)
$.
Obviously $\Pi_2$ conjugates the shift $\sigma_2$ to $z\mapsto z^2$ on $\Sen$
 with {\mapfromto {\sigma_2} {\Sigma_2} {\Sigma_2}} the
shift map:
$\sigma_2(\epsilon_1,\epsilon_2,\ldots,\epsilon_n,\ldots) =
(\epsilon_2,\epsilon_3,\ldots,\epsilon_{n-1},\ldots)$.
Moreover, 
we equip $\Sigma_2$ with the lexicographic ordering:
$\vepso=(\eps_1^1,\eps_2^1,\ldots,\eps_n^1,\ldots) <
(\eps_1^2,\eps_2^2,\ldots,\eps_n^2,\ldots)=\vepst$ iff
$\eps_k^1=\eps_k^2$ for $1\leq k < m$ and $\eps_m^1<\eps_m^2$
for some $m\in\N$.

Write the upper half-arc $I_0=[1,-1]\subset\Sen$  and lower $I_1=[-1,1]\subset\Sen$.
An itinerary of a point $z\in\Sen$ under the map $z\mapsto z^2$ is a sequence
$\veps=(\epsilon_1,\epsilon_2,\ldots,\epsilon_n,\ldots)$
with the property that for all $n\in\N$: $z^{2^n}\in I_{\epsilon_{n+1}}$.
The reader shall easily verify that for each $\veps\in\Sigma_2$
the point $\Pi_2(\veps)$ is the unique point of itinerary $\veps$ under $z^2$.
Moreover two sequences $\vepso<\vepst$ are the common itineraries of
a point $z$ if and only if $z^{2^n}=1$ for some minimal $n\geq 0$
and equivalently for this $n$
$\eps_k^1=\eps_k^2$ for $1\leq k < n$,
$0<\eps_n^2=\eps_n^1+1<2$
and $\eps_k^1=1$, $\eps_k^2=0$ for $n<k$.

Defining itineraries for $\Bla$ by the same algorithm as for $z^2$ above,
i.e.~$\Bla^n(z)\in I_{\epsilon_{n+1}}$,
we obtain exactly the same statements for $\Bla$.
For example $h\circ\Pi_2$ conjugates the shift $\sigma_2$ to $\Bla$,
any  itinerary for $\Bla$ determines a unique point of $\Sen$ and
a point has two itineraries if and only if $\Bla^n(z)=1$ for some $n$.

We shall now construct accesses to these points. Called parabolic rays, they sit in a tree. 
We explain the construction of this tree in $\D$ instead of $\widehat \C\setminus \overline \D$  to be  more   visual.

For each $j=0,1$ the open sector $S_j$ spanned by the arc
$I_j$, i.e. the interior of the convex hull of the union of $I_j$
and $0$, is mapped univalently onto $\D\Sm[1/3,1]$. The parts $[-1,0]\subset \R$ and $[0,1]\subset \R$ of the boundary 
are each mapped (homeomorphically) onto    $[1/3,1]$ which is 
forward invariant.
Let $z_\emptyset=0$ and $T_\emptyset := \Bla^{-1}([0,1/3]) = [0,z_0]\cup [0,z_1]$, 
where $z_0={\frac {i}{\sqrt 3}}$ and $z_1=-z_0$.
Since 
$
T_\emptyset\subset \D\Sm[1/3,1]$, define  $T_j$ to be  the connected component of
$\Bla^{-1}(T_\emptyset)$ containing $z_j$.
Define recursively after $n\in\N^*$ and for each
$(\eps_1,\eps_2,\ldots,\eps_n)\in\Sigma_2$ the
point $z_{\eps_1,\eps_2,\ldots,\eps_n}$ as the unique point of the
preimage $\Bla^{-1}(z_{\eps_2,\ldots,\eps_n})$ belonging to
$T_{\eps_1,\eps_2,\ldots,\eps_{n-1}}$. Define then
$T_{\eps_1,\eps_2,\ldots,\eps_n}$ to be the connected component of
the preimage $\Bla^{-1}(T_{\eps_2,\ldots,\eps_n})$ containing
$z_{\eps_1,\eps_2,\ldots,\eps_n}$.

 \begin{figure}[h]
\begin{center}
\includegraphics[height=5 cm]{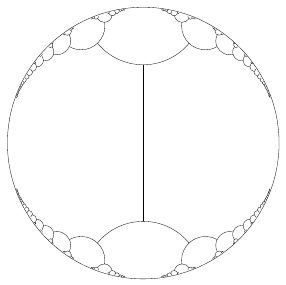}\quad \includegraphics[height=6 cm]{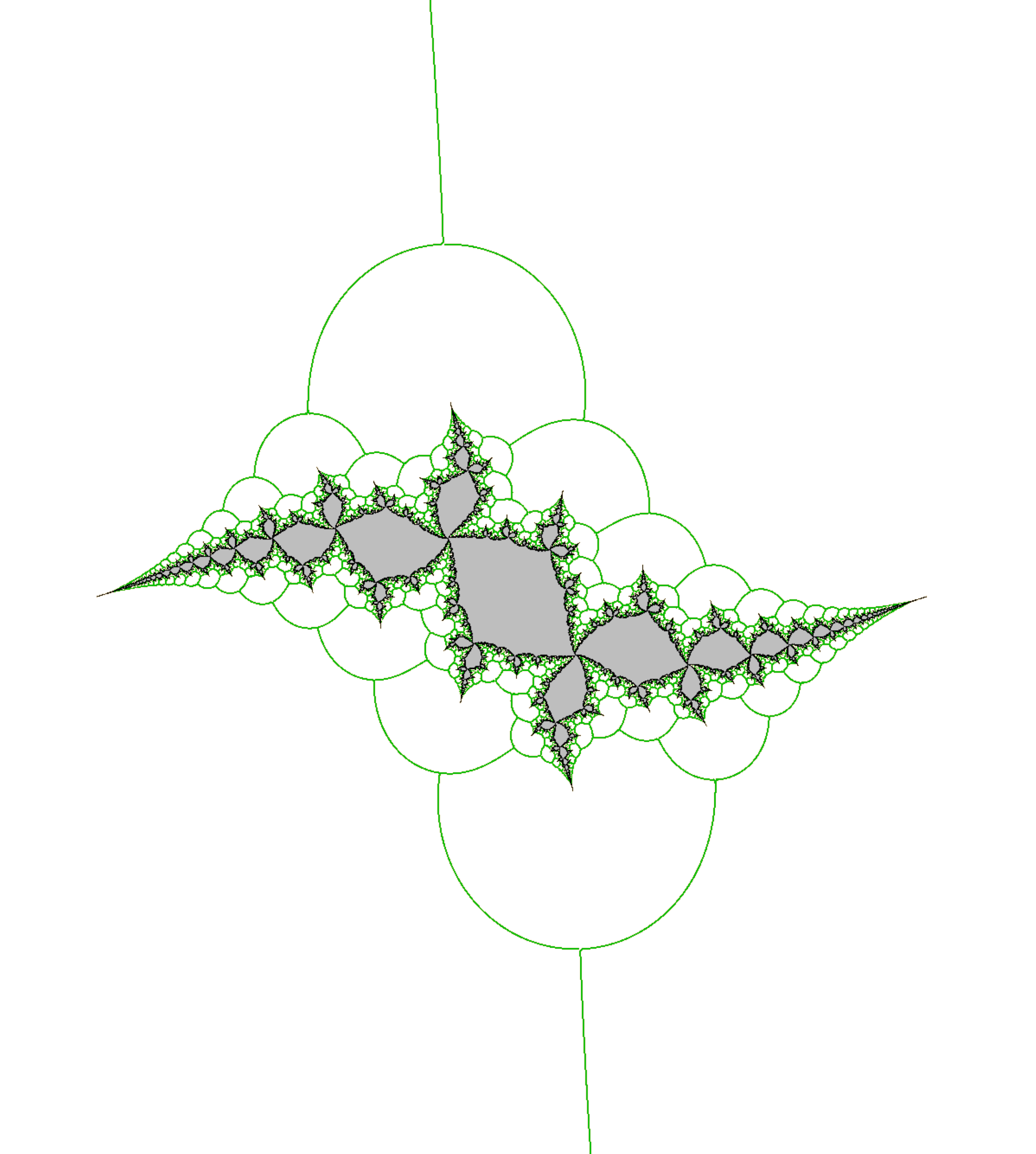}
\end{center}
\caption{The infinite dyadic tree $\TTT$
containing the internal parabolic  rays of $\Bla$ and the corresponding one for $g_B$.\label{pararaysinD}}
\end{figure}

 \begin{figure}[h] \label{BMonechess}
\begin{center}\includegraphics[height=2.8 in]{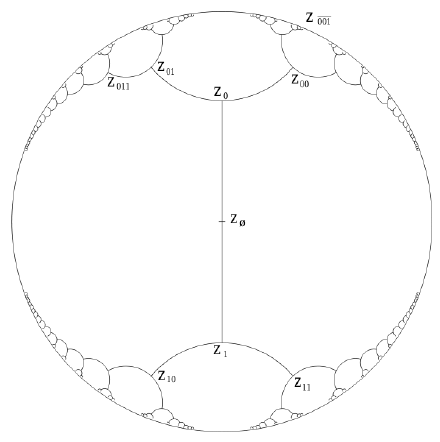}
\  \includegraphics[height=2.8 in]{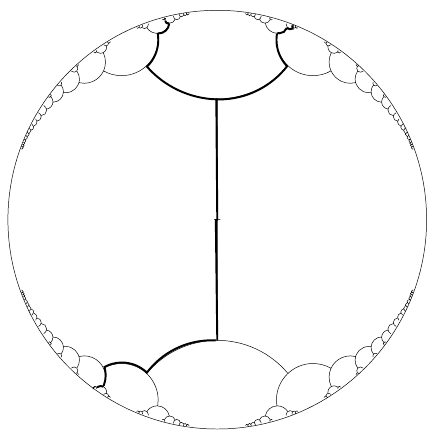}
\end{center}
\caption{Construction of a parabolic ray for $\Bla$ of rotation number $1/3$ }
\end{figure}

Define for each $n$ the  dyadic trees 
$$\TTT_n:=\bigcup_{k=0}^n\Bla^{-k}(T_\emptyset)\quad \hbox{   so that   }\quad   \TTT_n = \TTT_{n-1}\cup\bigcup_{(\eps_1,\eps_2,\ldots,\eps_n) \\
\in\{0,1\}^n} T_{\eps_1,\eps_2,\ldots,\eps_n}.$$
Then  define $$ \TTT:=\bigcup_{k=0}^\infty
\Bla^{-k}(T_\emptyset)$$ with boundary (in)  $\Sen$.

\REFDEF{InandOutRays}
For $\veps\in\Sigma_2$ a {\it parabolic internal ray } $\hat{R_\veps}$ is the
minimal connected subset of $\TTT$ containing the sequence of
points $z_{\eps_1,\eps_2,\ldots,\eps_n}$, $n\geq 0$ (enterpreting $n=0$ as $z_\emptyset$).

A parabolic external ray $R_\veps$ is the image of $ \hat{R_\veps}$ by $ z\mapsto 1/\overline{z}$. 
\ENDDEF
In order to stay close to the notations for quadratic polynomials, it will be convenient to identify $\T$ 
and the Julia set $\Sen$ for $\Bla$. This motivates the following definition.
\REFDEF{angledRays}
We shall say that $\theta\in \T$ is the external angle 
of the point $h(e^{i2\pi\theta})$. And write $R_\theta$ for the ray $R_\veps$, 
where $\veps$ is a binary expansion of $\theta$ mod $1$. 
In the special case where $\theta$ has two binary expansions $\veps_1$, $\veps_2$ 
we shall write $R_\theta := R_{\veps_1}\cup R_{\veps_2}$.
\ENDDEF

For $j\in\{0, 1\}$ the boundary $]-1,1[$ of $S_j$ in $\D$ is forward invariant and $z_j\in S_j$. 
It follows that the set $S_j\cup\{0\}$ contains any of the rays $R_\veps$ with $\eps_1 = j$. 
Moreover as $\Om$ is also forward invariant and disjoint from $T_\emptyset$ we even have  that $S_j\cup\{0\}\setminus\Omega$ contains any of the rays $R_\veps$ with $\eps_1 = j$. 
It follows that we may define parabolic rays in parameter space by 
$\RR^\Mone_\veps= \Upsilon^{-1}(R_\veps)$.
See also \defref{rayonsparam}.

 %%%%%%%%%%%%%%%%%%%%%%%%%%%%%%%%%%%%%%%%%%
\subsubsection{Parabolic rays for the rational map $g_B$}\label{s:paraborays}

Parabolic rays for $g_B$ are defined as pre-images of the external parabolic rays $R_\veps$ and $R_\theta$.

\REFDEF{rayonsdyn}
Let  $B\in \Mone$ and $\vepsilon$ be an itinerary. 
The {\rm parabolic dynamical ray} for $g_B$ of itinerary $\vepsilon $ is by definition $\RR^B_\vepsilon=h_B^{-1}(R_{\underline\epsilon})$. 
And the {\rm parabolic dynamical ray} for $g_B$ with angle $\theta$ is by definition $\RR^B_\theta=h_B^{-1}(R_{\theta})$. 
\ENDDEF
%\vspace{-3 cm}
\begin{figure}[h]
\begin{center}
\hbox{  \includegraphics[height= 8 cm]{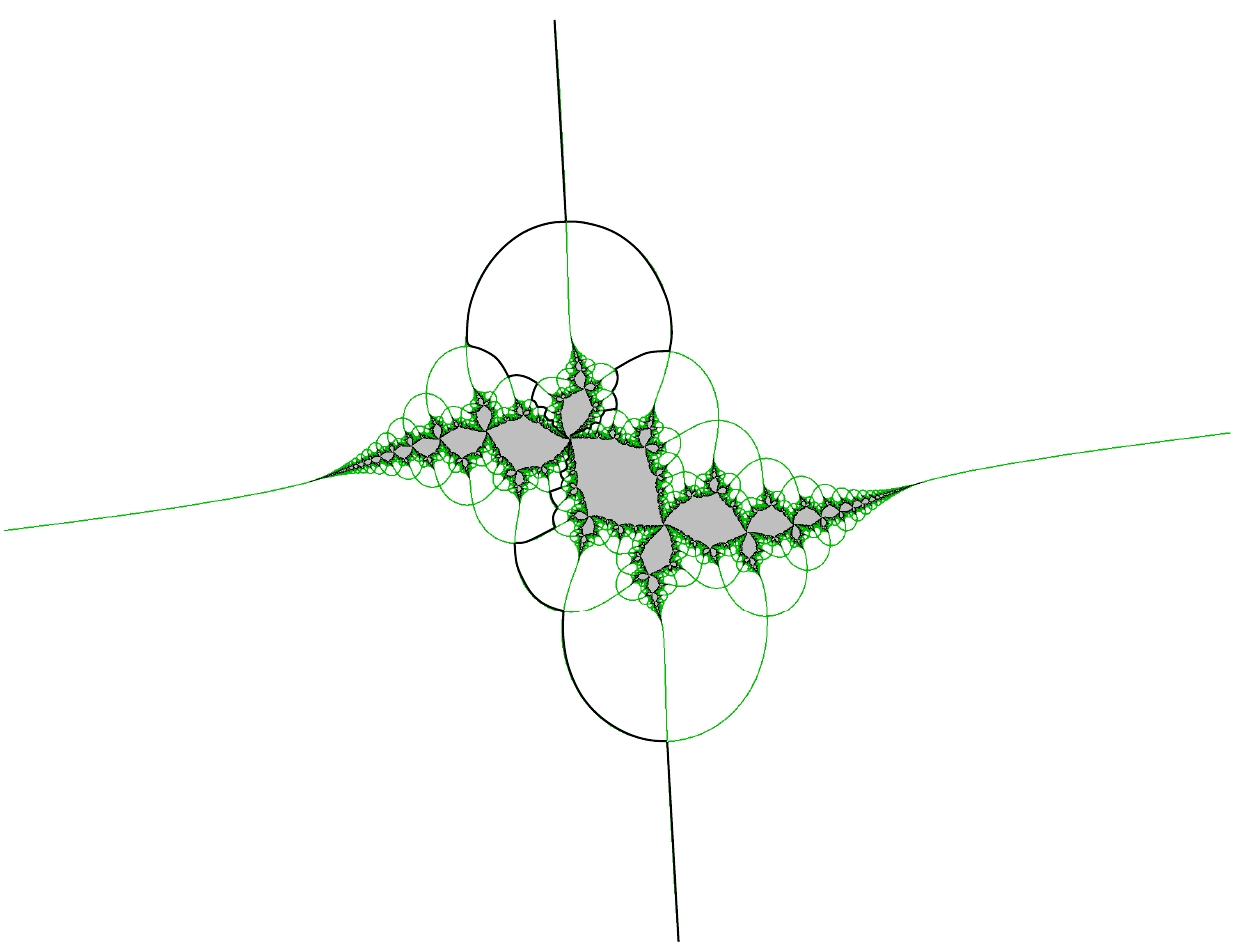}
\includegraphics[height= 8 cm]{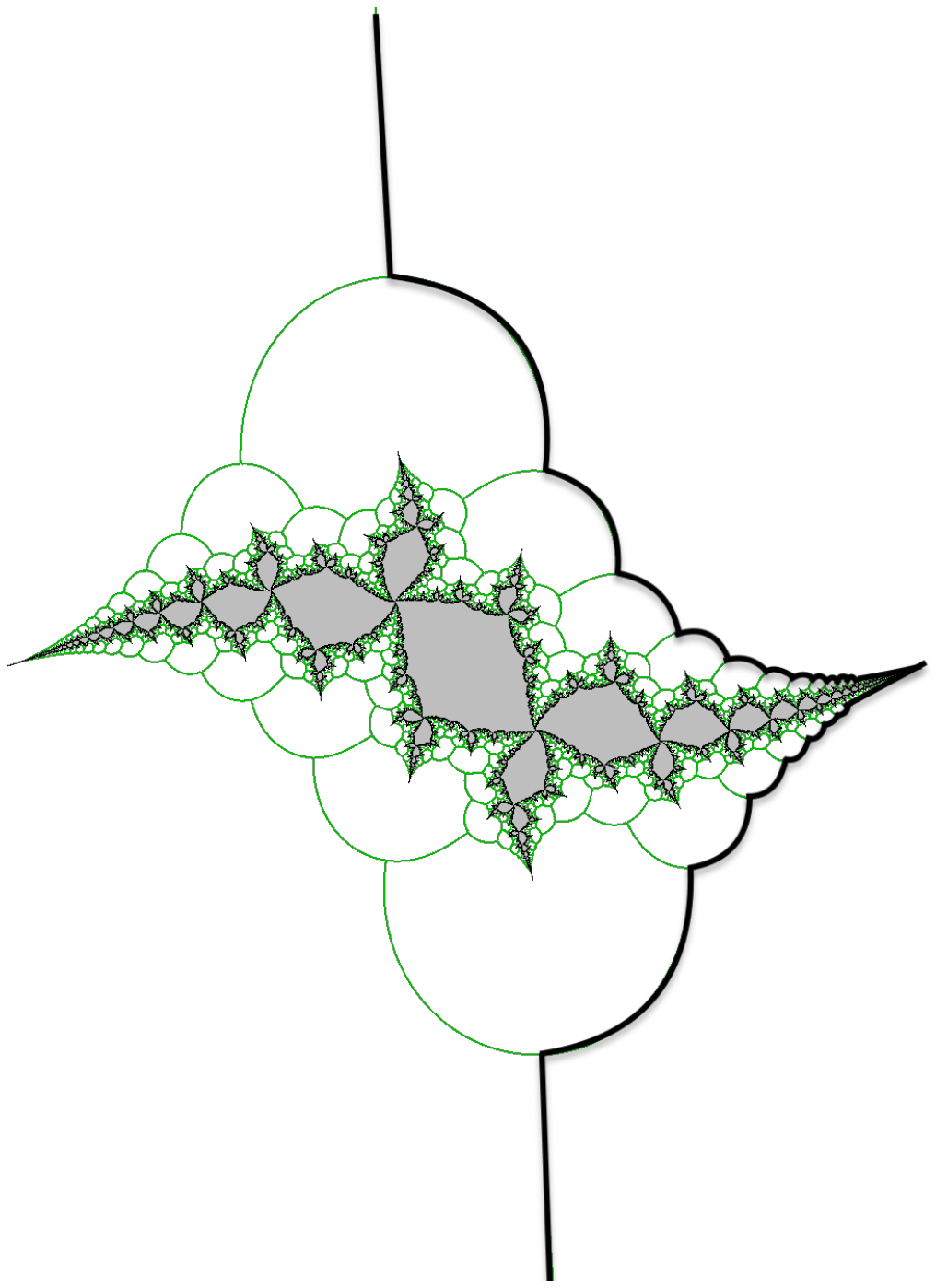}}
  \end{center}
%\vspace{-3 cm}
\caption{The $1/3$ cycle of parabolic rays $\RR^B_{1/7}, \RR^B_{2/7}$ and $\RR^B_{4/7}$ in $\C\setminus J_B$ on the left and 
the ray (pair) $\RR^B_0$ on the right (also in $\C\setminus J_B$).}
\end{figure}
 
 \REFDEF{rayonsparam} The {\rm parabolic parameter ray} of itinerary ${\underline\epsilon}$  is defined by $\RR^\Mone_{\underline\epsilon}=\Upsilon^{-1}(R_{\underline\epsilon})$. 
Similarly, the {\rm parabolic parameter ray} of angle $\theta$ is defined by 
$\RR^\Mone_\theta=\Upsilon^{-1}(R_\theta)$.
\ENDDEF

We say that a $q$ cycle of rays $\RR_0, \ldots \RR_{q-1}$ for $g_B$
landing on a common $k$ periodic point $z$ and numbered in the counter clockwise order
around $z$ defines the \emph{combinatorial rotation number} $p/q$, $(p,q)=1$
iff $g_B^k(\RR_j) = \RR_{(j+p)\mod q}$.
\REFTHM{simplelanding}
Let $B\in\Mone$. For any {\rm(}pre-{\rm)}periodic argument $\veps\in\Sigma_2$,
i.e.~$\sigma^k(\sigma^l(\veps))=\sigma^l(\veps)$,
the parabolic ray $\RR=\RR^B_\veps$ converges to a $g_B$ {\rm(}pre-{\rm)}periodic point
$z\in J(g_B)$ with $g_B^k(g_B^l(z))=g_B^l(z)$.
If  the argument is periodic {\rm(}i.e.~$l=0${\rm)}, let $k'$ denote the exact period of $z$ and let $q=k/k'$.
Then the ray $\RR$ defines the combinatorial rotation number $p/q$, $(p,q)=1$ for $z$.
The periodic point $z$ is repelling or parabolic with multiplier $\e^{i2\pi p/q}$.
Moreover any other external parabolic  ray landing at $z$ is also $k$-periodic and
defines the same rotation number.
\ENDTHM
This is a standard result which in its initial form is due to
Sullivan, Douady and Hubbard, for the polynomials. See the proof in  \cite[Th. A and Prop. 2.1]{Petersen1}, it goes through for parabolic rays.

And conversely
\REFTHM{douadylandingthm}
If $B\in\Mone$ and $z\in J(g_B)$ is any repelling or parabolic  periodic point. 
Then there is a periodic parabolic ray landing at $z$. 
It defines for $z$ its (unique) combinatorial rotation number. 
In particular for $B\in\Mone$ with $A=1-B^2\notin \Dbar$ the fixed point $\al_B$ has a 
combinatorial rotation number.
\ENDTHM
\PROOF
Since $g_B$ has degree $2$ the parabolic basin for $\be_B$ is completely invariant and 
thus the Theorem is a special case of \cite[Th. B]{Petersen1}.
\ENDPROOF
%%%%%%%%%%%%%%%%%%%%%%%%%%%%%%%%%%%%%%%
\vskip 1em \noindent{\bf The non-connected case:}

Assume now that $B\notin \Mone$ so that  the Julia set $J(g_B)$ is a Cantor set. 
The map  $h_B$ (\lemref{conj}) is well defined on $\wtOm^B=\Omega^B_{n_B}$, so that 
$\delta_\veps:=h_B^{-1}(R_\veps)$ is well defined 
(it is the pull-back of the part in $\tilde \Omega=h_B(\tilde\Omega^B)$). 
The part  $\delta_\veps$ is the beginning of the dynamical ray (as before). 
Using the relation $g_B\circ h_B=h_B\circ \Bla$ one can define the ray, 
until it bumps on an iterated  pre-image of the second critical value $v_B$, as follows. 
Define recursively $\delta^{n+1}_\veps$ as the connected component 
of $g_B^{-1}(\delta^n_{\sigma(\veps)})$ 
containing  $\delta^n_\veps$, with $\delta^0_\veps=\delta_\veps$.  
For $n\geq 0$ define $\de^n_\veps$ as the connected component of 
$g_B^{-n}(\de_{\si^n(\veps)})$ containing $\de_\veps$.

\REFLEM{prerays}
Let $\veps\in\Si_2$. If the critical value $v_B$ does not belong to  $\delta_{\si^j(\veps)}$ for any $0<j\leq n$, 
then the set $\delta^n_\veps$ is a  simple curve. 
Moreover $h_B$ has a univalent analytic extension to a neighbourhood of $\delta^n_\veps$.
\ENDLEM

\REFDEF{raynonconnected} If the critical value $v_B$ does not belong to  $\delta_{\si^n(\veps)}$ for any $n$, 
then we define {\it the dynamical ray of itinerary} $\veps$ by 
$$\RR^B_\veps:=\bigcup_n\delta^n_\veps.$$   

If the critical value $v_B$ belongs to $\delta_{\si^n(\veps)}$ for some $n$,
then we say that the ray $\RR^B_\veps$  bumps on some (pre)-critical point in $g_B^{-n}(-1)$ 
and that the ray is defined until this (pre)-critical point by the same procedure.
\ENDDEF
\REFREM{equirays}
As an immediate consequence of the definition of rays, the conjugacy $h_B$ has a unique analytic extension along rays. 
In particular we have 
$$v_B\in \RR^B_\veps \iff B\in\RMone{\veps}.$$
\ENDREM

The landing property given by Theorem~\ref{simplelanding} in the connected case, translates in the non connected case as follows\,: 

\REFTHM{landing} 
Let $B\in\Hpb$. For any {\rm(}pre-{\rm)}periodic argument $\veps\in\Sigma_2$, either 
the parabolic ray $\RR^{B}_\veps$  bumps on the critical point $-1$ or it converges to a $g_B$ {\rm(}pre-{\rm)} periodic point
$z\in J(g_B)$ with $g_B^k(g_B^l(z))=g_B^l(z)$.
If periodic {\rm(}i.e.~$l=0${\rm)}, let $k'$ denote the exact period of $z$ and let $q=k/k'$.
Then the ray $\RR^{B}_\veps$ defines a combinatorial rotation number $p/q$, $(p,q)=1$ for $z$.
The periodic point $z$ is repelling or parabolic with multiplier $\e^{i2\pi p/q}$.
Moreover any other external parabolic  ray landing at $z$ is also $k$-periodic and
defines the same rotation number.
\ENDTHM
The following  stability statement will be crucial in the sequel\,:
\REFLEM{stability} 
Let $\Bc\in\Hpb$ and assume that the critical point $-1$ 
is not on the forward orbit of $\overline {\RR^\Bc_\veps}$ 
and that $\RR^\Bc_\veps$ is landing on either a pre-image of $\infty$ 
or on a point, which is pre-periodic to a repelling periodic point. 
Then, there exists a neighborhood $U$ of $\Bc$ such that 
for any $B\in U$, the ray $\RR_\veps^B$ lands at a pre-periodic point 
and there exists a holomorphic motion $\psi : U\times \overline \RR^\Bc_\veps\to \C$ 
such that $\psi_B(\overline \RR^\Bc_\veps)=\overline \RR_\veps^B$. 
\ENDLEM
\PROOF In the case where the landing point is pre-repelling, the proof is similar to the one of Douady-Hubbard in the case of quadratic polynomials :  it is based on   the implicit function Theorem.
Note that  $\RR^{B}_{\overline{0}}$ always  lands at $\infty$ when it is defined. 
Hence, if $v_B$ is not on the closure of $\RR^{B}_{\overline{0}}$, 
this ray varies holomorphically (in this family) and the pre-images $\RR^{B}_{\overline{0}}$, 
$\RR^{B}_{1\overline{0}}$ cannot break on the critical point $-1$. 
Indeed, on any disk in $\C\setminus \RR^\Mone_{\overline{0}}$, 
we have a holomorphic motion of the arc $h_\Bc^{-1}([0,1])$ connecting the 
critical point $1$ to its critical value $g_\Bc(1)$. 
Pullback by iteratively along $\RR^{B}_{1\overline{0}}$ 
by the dynamics we never encounter the second critical value $v_B$ 
and so lifting the holomorphic motion gives a holomorphic motion 
of all the ray parameterized by this disk. 
By the $\lambda$-Lemma it extends to the closure of $\RR^{B}_{\overline{0}}$.
Note that  the only parameter $B\in\Mone$ for which $v_B$ is on the closure of 
$\RR^{B}_{\overline{0}}$ is $B=0$. 
The similar statement hold for $\RR^\Mone_{\overline{1}}$.
\ENDPROOF
\begin{figure}[h]
\vspace{-5.6cm}
\centerline{\includegraphics[height=21cm]{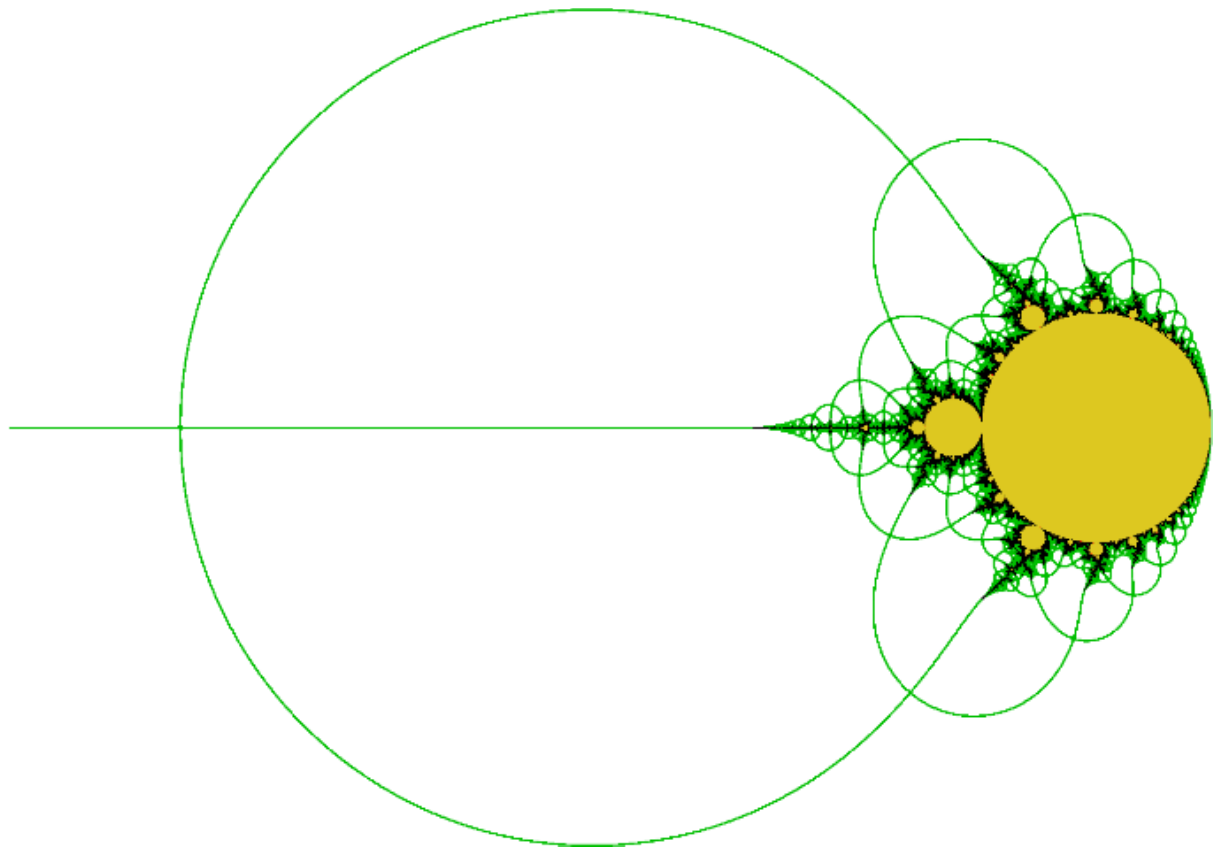}}
\vspace{-5.6cm}
\caption{Parabolic chess board outside $\Mone$, viewed in the $A$-parameter plane.\label{Monechess}}
\end{figure}
\REFCOR{stabilitydomain} In any disk contained in the complement of 
$\Chat\setminus \cup_i\overline  \RR^{\Mone}_{2^i\theta}$, 
the ray $\RR^B_{\theta}$ admits a holomorphic motion and so does its closure 
(by the $\lambda$-Lemma).
\ENDCOR

%%%%%%%%%%%%%%%%%%%%%%%%%%%%%%%%%%%%%%
\subsection{Limbs  of $\Mone$}\label{s:limbs}

Similarly to Douady and Hubbard  description of $\Mbrot$, the parabolic Mandelbrot set $\Mone$ 
can be described in terms of limbs sprouting out of the central, 
period $1$ (relative) hyperbolic component $\Hz$. 
\REFDEF{parabo_limbs}
For $0< p/q<1$ an irreducible rational we define the $p/q$ limb $\LL^\Mone_{p/q}$ as 
$$
\LL^\Mone_{p/q} = \{ B\in\Mone \mid \al_B \textrm{ has rotation number } p/q \}.
$$
\ENDDEF
By uniqueness of the rotation limbs of different rotation numbers are disjoint and 
moreover $B_{p/q}$ with $1-B_{p/q}^2 = A_{p/q}:=\e^{i2\pi p/q} \in\LL^\Mone_{p/q}$ is called the root of the limb.

\REFTHM{Limbsdecomposition} 
\REFEQN{M1_limbdecomposition}
\Mone= \Hbz\cup\bigsqcup_{\frac pq \neq \frac 01}\LL^\Mone_{p/q}
\ENDEQN
and each limb $\LL^\Mone_{p/q}$ is compact and connected with $\Hbz\cap\LL^\Mone_{p/q} = \{B_{p/q}\}$. 
\ENDTHM
Though this is similar to the Mandelbrot case we give in this section a complete proof 
following the approach by Milnor in~\cite{Milnorasterisque}. 
\PROOF
The decomposition \eqref{M1_limbdecomposition} of $\Mone$ 
is an immediate consequence of \thmref{douadylandingthm} 
and $\Hbz\cap\LL^\Mone_{p/q} = \{B_{p/q}\}$ follows from the combination 
of \thmref{douadylandingthm} and \thmref{simplelanding}. 
The compactness and precise localisation of the limb follows from \corref{limbs_compact} of \thmref{wake} below.
\ENDPROOF

The following discussion is most conveniently taken in the $A$-parametrization, where $\Hz = \D$. 
For this subsection we shall henceforth use the $A$-parameterization.
For  $A\in \C\setminus \{1\}$, 
the multiplier of the finite fixed point $\al(A)$ is $A\in \C$, so it is attracting when $A\in\D=D(0,1)$, 
neutral when $A\in\partial \D$ and repelling when $A\in\C\setminus\overline \D$.
For each $p/q\neq 1$, with  $(p,q)=1$, the parameter
$A_{p/q}=\e^{i2\pi p/q}\in\partial \D\setminus\{1\}$  belongs to $\Mone$,  so that   there is a parabolic external ray converging to $\alpha(A)$ by \thmref{douadylandingthm} and  \thmref{simplelanding}.  
This ray has rotation number $p/q$. 
Let us denote by $\RR^B_{\theta_-(p/q)}$ and $\RR^B_{\theta_+(p/q)}$, recall $A = 1-B^2$, 
the rays in the cycle that are adjacent  to  the critical value ({\it i.e.} to the Fatou component containing it). 
 In~\thmref{wake}   we  prove that the corresponding parameter rays $\RR^\Mone_{\theta_\pm(p/q)}$ lands at $A_{p/q}$ and that they cut off  a wake $\W^\Mone(p/q)$. We call $(p/q-)$ (derooted) limb of $\Mone$, 
 the set $\LL^*_{p/q}:=\Mone \cap \W^\Mone(p/q)$.
 
For $\theta = p/2^l, l\geq 0$ a dyadic angle we shall say that $\RR^*_\theta$, $*\in\{\Mone, B\}$ lands if the two rays 
$\RR^*_{\veps_i}$ land on the same point, where $\veps_0, \veps_1$ are the two dyadic expansions of $\theta$.
We obtain  quite  precise properties of the landing  in the parameter plane in the following\,:

\begin{figure}[h]
\vspace{-4.6cm}
\centerline{\includegraphics[height=21cm]{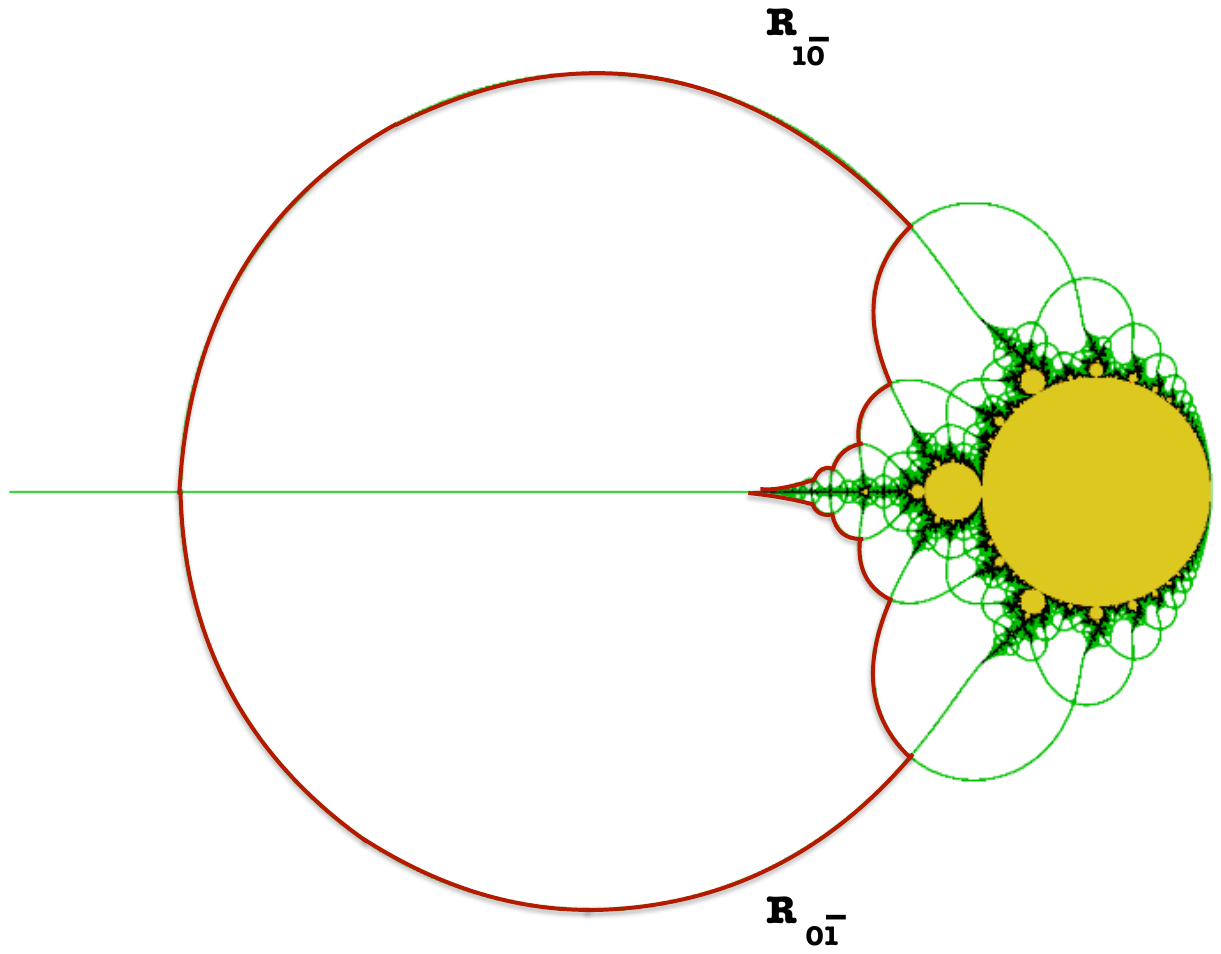}}
\vspace{-4.6cm}
\caption{The dyadic rays $\RR^\Mone_{0\overline{1}}$, $\RR^\Mone_{1\overline{0}}$ corresponding to the dyadic angle $1/2$ They bound the disk $D_{1/2}^\Mone := \Upsilon^{-1}(D_{1/2})$, 
viewed in the $A$-parameter plane.\label{RMonechess}}
\end{figure}

\REFTHM{DouadyHubbardparamland} For every
pre-periodic {\rm(}i.e. rational{\rm)} angle $\theta$, the parameter ray 
 $\RMone{\theta}$  lands.
More precisely, suppose $2^{k+l}\theta \equiv 2^l\theta \mod 1$
with  {\rm(}period {\rm)} $k>0$ and $l\geq 0$ minimal. 
\ENUM \item If $\theta$ is periodic {\rm(}$l=0${\rm)} then
$\RMone{\theta}$ lands on a parameter $A=1-B^2$ for which the
corresponding dynamical ray $\RR_\theta^B$ lands at a parabolic
periodic point $z(B)$, with exact period $k'|\, k$ and with
multiplier $\lambda=(g_B^{k'})'(z(B))$ a primitive $k/k'$-th root
of unity. 
\item If
$l>0$ then $\RMone{\theta}$ lands at a parameter $A$ for which the
corresponding dynamical ray $\RR_\theta^B$  lands on $v_B$ and
for which $g_B^{k+l}(v_B)=g_B^l(v_B)$ is a  periodic point of
exact period $k$. This periodic point  is repelling if $k>1$ and is the parabolic fixed point $\infty$ for $k=1$.
Moreover for any dynamical ray $\RR_{\theta'}^B$ landing on $v_B$,
the corresponding parameter ray $\RMone{\theta'}$ lands at $A\in\Mone$.
\ENDENUM
\ENDTHM

\REFREM{Dyadicrays}
In particular in point 1. $\theta=0$ both rays $\RR_{\overline 0}^B, \RR_{\overline 1}^B$ land at the parabolic fixed point $\infty$ of multiplier $1$, and there is no other ray landing at $\infty$.
In point 2., $\theta = p/2^l$, $l>0$ the landing point of $\RMone{\theta}$ is not on $\partial \Hz$ (but at the so called ``$\theta$-dyadic tip'' of $\Mone$), so that these two rays landing at the same point do not define a limb but bounds a disk. 
\ENDREM

%%%%%%%%%%%%%%%%%%%%%%%%%%%%%%%%%%%%%%%%

\subsubsection{Landing of parameter rays} 

\REFLEM{periodiclanding}Let $\theta$ be a $k$-periodic angle,  then
$\RMone{\theta}$ lands at a parameter $A\in\Mone$. Moreover,   
 the
corresponding dynamical ray $\RR_\theta^B$ lands at a parabolic
periodic point $z(B)$ where $A=1-B^2$. If $k=1$,  $z(B)=\infty$. If $k>1$,  $z(B)$ has exact period 
$k'|\, k$ and multiplier $\lambda=(g_B^{k'})'(z(B))$ a primitive $k/k'$-th root
of unity. 
\ENDLEM
\PROOF The argument is classical. Let $A$ be any accumulation point
of $\RMone{\theta}$. Since  $A$ is in  $\Mone$ the Julia set is connected.
 For $B$ such that  $A=1-B^2$, the ray $R=\RR^B_\theta$ lands at a $k'$-periodic point  $z(B)$ of $J(g_B)$ with $k'|\, k$. It is either repelling or parabolic. If it is repelling or if it lands at the parabolic point $\infty$,  \thmref{stability} gives a  holomorphic motion of $\overline R$ in a neighborhood of $A$ since   $z(B)$ cannot be critical (it is periodic).
But this contradicts the fact that if $A'$ is close to $A$ on $\RMone{\theta}$,
the critical value is on the ray $\RR^{B'}_\theta$  (\lemref{equirays}) so that   the two preimages (one  is in the cycle) bump on the critical point (since $\theta$ is periodic).
Hence, $z(B)$ is a parabolic point of period $k'$ dividing $k\ge1$, with
multiplier $\lambda=(g_B^{k'})'(z(B))$ a primitive $k/k'$-th root
of unity.

The  set of parameters $B\in \C$ such that $g_B$ has a parabolic cycle of period $k'|\, k$   with $k>1$ is  included in  
$\{B\in\C\mid\exists z\in\C, \  g_B^{k'}(z)=z,\;  (g_B^{k'})'(z)=e^{2i\pi jk'/k },\  j\in\{0,\cdots,k/k' \}\}$. This set  is finite since it is defined by the equations in $(z,B)$ of two (relatively prime) polynomials.

Therefore, the accumulation set of  $\RMone{\theta}$ is finite, so it reduces to one point.
\ENDPROOF

\REFLEM{preperiodiclanding} Let  $\theta$ be a strictly preperiodic angle  :  
$2^{k+l}\theta \equiv 2^l\theta \mod 1$ with $k>0$ and $l> 0$ minimal. The parameter  ray  $\RMone{\theta}$  lands.  Moreover, if $k=1$, the
corresponding dynamical ray $\RR_\theta^B$  lands on the critical value $v_B$ and
 $g_B^{l}(v_B)=\infty$. 
\ENDLEM
\PROOF
As before let $A$ be an accumulation point of the ray  $\RMone{\theta}$. The dynamical ray $\RR^B_\theta$ lands at a strictly preperiodic point $z(B)$ of $J(g_B)$. The  point 
$g_B^l(z(B))$ is periodic,  either repelling or parabolic.  Assume first that $k>1$. As in previous lemma, the set of parameters such that $z(B)$ is parabolic is finite. Now if $z(B)$ is repelling, the critical point is in the orbit of $z(B)$, by the stability Lemma.
This situation also corresponds to a finite number of $B\in \C$ since  $B$ has to satisfy a 
polynomial equation (the critical point is pre-periodic).
Since the accumulation set  of $\RMone{\theta}$ is connected and finite, it reduces to one point. Therefore, the parameter ray  lands. 

We consider now the case $k=1$.    The angle is  dyadic and $\RR^B_\theta$ lands at the critical value $v_B$, so  $g_B^{l}(v_B)=\infty$. This equation also gives a finite number of parameters so that the parameter ray lands. \ENDPROOF

To achieve the proof of   \thmref{DouadyHubbardparamland} we need to define Wakes as in the Mandelbrot case.

%%%%%%%%%%%%%%%%%%%%%%%%%%%%%%
\subsubsection{Wakes }
We consider now for $p/q\notin\{0, 1\}$, the parabolic parameter $A_{p/q} = \e^{i2\pi p/q}$ and the  $q$-cycle of external parabolic rays landing to the $\alpha$ fixed point, with  angles 
$0<\theta_0<\theta_1<\ldots<\theta_{q-1}<1$ of combinatorial
rotation number $p/q$ {\rm(}{\it i.e.} $2\theta_i \equiv \theta_{(i+p)\mod q}
\mod 1${\rm)} (defined at the beginning of the \subsecref{s:limbs}). Denote by  $\mathcal I=(\theta_-,\theta_+)$ the smallest 
interval in $\displaystyle \Sen\setminus \bigcup_{i\ge 0}\theta_i$.

\REFLEM{existencerotation} Let $A$ be a parameter outside of $\Mone$ on some external ray of angle $t$.  The dynamical rays 
$\RR^B_{\theta_0},\RR^B_{\theta_1},\ldots,\RR^B_{\theta_{q-1}}$ land at the repelling fixed point if an only if $t$ belongs to $\mathcal I$.
\ENDLEM

\PROOF The proof is  the same as Lemma~2.9 of~\cite{Milnorasterisque} (which deals 
with all kind of cycles of quadratic polynomials). We recall it briefly. Note that the two rays $\overline{\RR^B_{t/2}}\cup \overline{\RR^B_{t/2+1/2}}$ crash on  the critical point and so partition the plane $\C$ in two sides.  It follows  that  two dynamical rays land at the same point if and only if they have the same itinerary with respect to this partition of $\C$.

On the other hand, the arcs in the complement of the cycle $\theta_0,\cdots,\theta_q$   are mapped injectivly to another complementary arc, except for one complementary arc that double covers $\mathcal I$.
Therefore, if $t\notin \mathcal I$, its preimages $t/2$ and $(t+1)/2$  belongs to different complementary arcs, or to the cycle, so that the rays of the cycle cannot land at the same point, or are even not defined until the end.

Now, for $t\in \mathcal I$,   $t/2$ and $(t+1)/2$  belongs to a complementary interval of length greater than $1/2$. So all the rays of the cycle  are  in the same component of the partition, they land at the same point. 
\ENDPROOF

\begin{figure}[h]
\vspace{-5.6cm}
 \begin{center}\includegraphics[height=21cm]{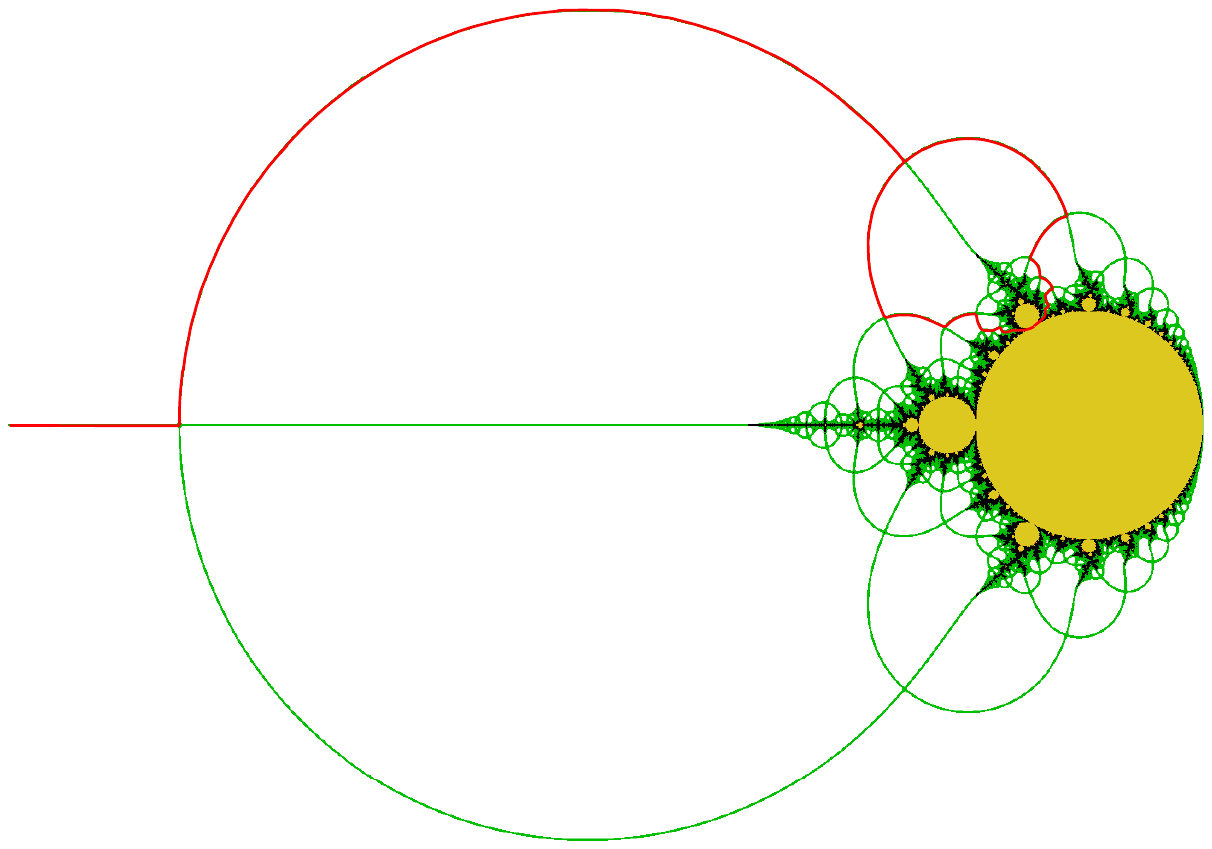}\end{center}
\vspace{-5.6cm}
\caption{Wake $1/3$}
 \end{figure}

\REFTHM{wake}
The parameter rays $\RMone{\theta_-},
\RMone{\theta_+}$ both land at  $A_{p/q}$.
 Moreover, the  curve $\RMone{\theta_-}\cup
\RMone{\theta_+}\cup A_{p/q}$ cuts the sphere into two connected components. 
Denote by 
$\W^{\Mone}(p/q)$ the one not containing $\D$. The dynamical rays $\RR^B_{\theta_i}, 0\le i\le q-1$, land at
a common repelling  fixed point if and only if $A=1-B^2\in \W^{\Mone}(p/q)$.
\ENDTHM
\PROOF The proof is similar to the one in Theorem~3.1 of~\cite{Milnorasterisque}.
Let $W$ be the set of parameters in $ \C$ such that the dynamical rays $\RR^B_{\theta_i}$ all  land at the same point which is a repelling  fixed point.  Note that such  rays cannot land at the fixed point $\infty$. By \lemref{existencerotation} $W$ is non empty; it  is an open set by the stablity property of such rays.
From the charaterization given by  \lemref{existencerotation}, a parameter ray 
 $\RR^\Mone_t$ belongs to $W$ if and only if $t\in \mathcal I$. The boundary of $W$ consists of parameters for which there is no stability of the dynamical  rays.  That is, either the critical point $c_1$ is on the cycle $\RR^B_{\theta_i}$
(it cannot be  at the landing point which is periodic) or the landing point of the rays is parabolic (\lemref{stability})\,: 
 $\displaystyle \partial W\subset (\bigcup_{0\le i\le q-1} \RMone{\theta_i})\cup \mathcal F$,  where $\mathcal F$ corresponds to the  finite set  of parameters  for which there is a parabolic point of period $q$.
 We deduce from this description of $\partial W$ and  from \lemref{existencerotation} that 
 $W$ is connected. Then the parameter rays  $\RMone{\theta_-}$ and   $\RMone{\theta_{+}}$  have to  land at a common  point $A$ of $\mathcal F$, since $W$ does not contain parameter rays of angle outside $\mathcal I$. 
For this parameter $A$,  the rays  $\RR^B_{\theta_i}$ land at a   parabolic $q'$-cycle   different from $\infty$ by \lemref{periodiclanding}, with $q'$ dividing $q$.
Assume that $A\neq  A_{p/q}$, 
then the map $g_B$ has a repelling fixed point 
(since it can have only one parabolic cycle in $\C$). 
There is a cycle of external rays, different from  $\RR^B_{\theta_i}$, landing at this fixed point.
 This cycle is stable (\lemref{stability}) in a small neighborhood since the critical point is in the basin of the parabolic cycle. This contradicts the fact that  $A$ is in the boundary of $W$ where these rays do not land at the fixed point.
 Therefore $A=A_{p/q}$ and  the statement follows.
\ENDPROOF

 \REFCOR{limbs_compact}
 The limbs $\LL^\Mone_{p/q}$ are compact, more pricesely 
 $$
 \LL^\Mone_{p/q} = \ov{\Mone\cap\W^\Mone(p/q)} = \LL^*_{p/q}\cup\{A_{p/q}\}.
 $$
 \ENDCOR

\REFDEF{ParDynWake}
For $A= 1 - B^2 \in\W^\Mone(p/q)$ we define the {\it dynamical wake} $\W_B(p/q)$ 
as the connected component of $\overline\C\setminus \overline{\RR^B_{\theta_-}\cup \RR^B_{\theta_{+}}}$ containing  the rays $\RR^B_{t}$ for $t\in(\theta_-,\theta_{+})$. 
\ENDDEF

Note that when $A\notin\W^\Mone(p/q)$, $p\not=0$ then 
$\overline{\RR^B_{\theta_-}\cup \RR^B_{\theta_{+}}}$ does separate $\C$ into two sub disks.

\REFCOR{vlocation}The parameter $A$ in $\W^{\Mone}(p/q)$, if and only if the second critical value 
$v_B$ is in the dynamical wake $\W_B(p/q)$.
\ENDCOR 
\PROOF It follows from the construction.  As explained in the proof of \lemref{existencerotation}, the only non injective interval of angles   double covers $\mathcal I$. One deduces easily that the sector/wake corresponding $\mathcal I$ contains the critical value $v_B$.
\ENDPROOF

\REFCOR{par_shringking_limbs}
The diameter of the limbs $\LL^{\Mone}_{p/q}$ tends to zero, when 
$q$ tends to $\infty$. 
\end{corollary}
\proof
Assume to get a contradiction that there is a sequence of Limbs  $\LL^\Mone_{p_n/q_n}$, 
with $q_n\to\infty$ whose  diameter does not go to zero. Then, one can find points $x_n,y_n\in \LL^\Mone_{p_n/q_n}$ 
converging to $x\neq y$  respectively. 
By  \thmref{Limbsdecomposition} these two points cannot both belong to $\overline \D$ (they would be separated by some wake). So one at least, say $y$,  correspond to a map that has a repelling fixed point  
and therefore is in some Limb $\LL(p/q)$. But this implies that the sequence 
$y_n$ enters in $\W^{\Mone}(p/q)$ for $n$ large. The contradiction comes from the fact that 
the wakes  $\W^{\Mone}(p_n/q_n)$ and $\W^{\Mone}(p/q)$ are disjoint.
Alternatively apply Yoccoz inequality in the form \cite[Theorem C]{Petersen1} to all $g_B$, 
$B\in\LL^\Mone_{p/q}$. To obtain that $\log$ of the limb is contained in the closed Euclidean disk 
of radius $r(q)$ and center $r(q)+i2\pi p/q$, where
$$
r(q) = \frac{\log 4}{q}
$$
Here the argument $4$ comes from the inequality $|\Bla'(z)| \leq 4$ on the unit circle.
\endproof

\vskip 1em \noindent 
{\bf Proof of \thmref{DouadyHubbardparamland}:}

Note that this Corollary \ref{vlocation} together with \lemref{periodiclanding} achieve the proof of 
part 1 of \thmref{DouadyHubbardparamland} in the case $k=1$. Thus it suffices to consider the case $k>1$.

\REFLEM{Misiurewicz} For a pre-periodic angle $\theta$, {\it i.e.}
$2^{k+l}\theta \equiv 2^l\theta \mod 1$ with  {\rm(}period {\rm)} $k>1$ 
and $l> 0$ minimal, if $A$ denotes the landing point of $\RR^\Mone_\theta$, 
 the corresponding dynamical ray $\RR_\theta^B$
  lands at $v_B$ and $g_B^{k+l}(v_B)=g_B^l(v_B)$ is a  repelling  periodic point of
exact period $k$. 
\ENDLEM
\PROOF Assume to get a contradiction that the external rays of the cycle of  angles $2^{i+l}\theta$ land at a parabolic periodic point. 
The parameter $A$ belongs to some wake $\W^{\Mone}(p/q)$, so that there is a cycle of external rays landing at the repelling fixed point. Let us denote by $\Gamma$ the union of these external rays  together with 
the fixed point.  The iterated pre-images $\Gamma_n=g_B^{-n}(\Gamma)$ give a partition of $\C$ that separates for $n$ large enough,  the external rays of the cycle of  angles $2^{i+l}\theta$ with $i\ge 0$  from the ray of angle $\theta$. Therefore $\Gamma_n$ separates the critical value $v_B$ from the external ray of angle $\theta$ (since it is in a Fatou component adjacent to one of the rays in the cycle  $2^{i+l}\theta$). 
Now the graphs  $\Gamma_n$ are stable, so for parameters $A'$ in a neighborhood of $A$,  but on the  ray  $\RR^\Mone_\theta$, the graph still separates the critical value $v_{B'}$  from the  ray $\RR^{B'}_\theta$ where $A'=1-B'^2$ (the elements stays in some different sectors). This contradicts the fact that $v_{B'}$ has to be on $\RR^{B'}_\theta$.
Therefore $\RR_\theta^B$
  lands at $v_B$. 
\ENDPROOF

Recall from page \pageref{Renormalization_Cantor_set_of_arguments} 
that the pair of periodic arguments $\theta, \theta'$, $0<\theta < \theta' < 1$ 
of a pair of parameter rays $\RR^\Mbrot_\theta$ and $\RR^\Mbrot_{\theta'}$ 
co-landing on a parabolic parameter defines a Cantor set $C(\theta,\theta')$. 
And that this lead to the definition of dyadic wakes and limbs of the corresponding copy 
of the Mandelbrot set. In view of \lemref{periodiclanding} and \lemref{Misiurewicz} 
we generalize this to $\Mone$ as follows.
Let $M^\Mone(\theta,\theta')$ be the copy of $\Mbrot$ in $\Mone$ with root rays $\RR^\Mone_\theta$ 
and $\RR^\Mone_{\theta'}$. 
\REFDEF{para_p_q_dyadic_wakes}
Define the dyadic wake  $\W^\Mone(\theta,\theta',r,s)$ of $M^\Mone(\theta,\theta')$ 
as the set bounded by the pair of co-landing parameter rays 
$\RR^\Mone_{\theta_0}, \RR^\Mone_{\theta_1}$, where $\theta_0$, $\theta_1$ are the 
the arguments bounding the $r/2^s$ gap in the Cantor set $C(\theta,\theta')$. 
And define the dyadic limb $\LL^\Mone(\theta,\theta', r, s)$ as the intersection 
$$
\LL^\Mone(\theta,\theta', r, s) := \W^\Mone(\theta,\theta',r,s)\cap\Mone.
 $$
 Moreover as for polynomials in the special case where $M^\Mone(\theta,\theta') = \Mbrot^\Mone_{p/q}$ 
 we shall write $\W^\Mone(p/q,r,s)$ for $\W^\Mone(\theta,\theta',r,s)$ and 
 $\LL^\Mone(p/q, r, s)$ for $\LL^\Mone(\theta,\theta', r, s)$
dyadic wakes and limbs associated with $\Mbrot^\Mone_{p/q}$. 
\ENDDEF
The root $B$ of the dyadic wake $\W^\Mone(\theta,\theta',r,s)$ is the $r/s$ dyadic tip of 
$\Mbrot^\Mone(\theta,\theta')$, 
that is the parameter such that the $q$-renormaliztion of $g_B$ is hybridly equivalent to 
the $r/s$ tip, i.e. the landing point of the parameter ray $\RR^\Mbrot_{r/s}$ of the Mandelbrot set.

%%%%%%%%%%%%%%%%%%%%%%%%%%%%%%%%%%%%%%%
\section{Parabolic Puzzles and Parabolic Para-puzzles}\label{s:ParabolicPuzzles}
We shall state and prove in \secref{Yoccoztheoremformone} a theorem for the parabolic Mandelbrot set $\Mone$ ananalogous to the Yoccoz parameter puzzle theorem for the Mandelbrot set (see~\cite{PascaleinTanLeisbook}).  The idea underlying the proof is also in this case to transfer the result obtained in the dynamical plane to the parameter plane using the trick of Shishikura to control the dilatation of the holomorphic motion in puzzle pieces.

Yoccoz theorem for the parabolic map $g_B$ was proved in~\cite{PR2}. We recall  briefly the proof here since we need the detailed construction of the parabolic puzzle.  Before let us  recall  the classical Yoccoz puzzle. Then, the construction of the parabolic puzzle will appear more natural even  in the parameter plane. 
%%%%%%%%%%%%%%%%%%%%%%%%%%%%%%%%%%%%%%%
\subsection{Yoccoz puzzle for Quadratic polynomials} 
For $c\in\Mbrot\setminus \overline {\Card}$,  $c$ belongs to some derooted limb $L^\star_{p/q}$. 
For the rest of this section we fix the reduced rational $p/q$,
but we shall only occasionally make reference to $p/q$. 
This motivates the following. 
Let $0<\theta_0 < \theta_1 < \ldots \theta_{q-1} <1$ denote the arguments of the unique $q$-cycle of 
rotation number $p/q$ for $Q_0$. 

Recall that the wake parameter wake $\WMbrotpq$ is the subset of parameter space $\C$ 
bounded by the parameter rays of arguments $\theta_{p-1}, \theta_p$. 
Recall further that $c\in\WMbrotpq$ if and only if the cycle of dynamical rays 
$\RR^c_{\theta_0},\ldots, \RR^c_{\theta_{q-1}}$ co-land on $\al_c$ and then also $c\in\W^c_{p/q}$, 
the dynamical wake bounded by the dynamical rays of arguments $\theta_{p-1}, \theta_p$. 
Moreover the strictly pre-periodic pre-image rays of arguments 
$\theta_0+\frac12 < \ldots \theta_{q-}+\frac12\subset ]\theta_{q-1},\theta_0[$ 
co-lands on $\al'_c$. 
It follows immediately that all $2q$ rays together with their landing points $\al_c, \al'_c$ 
move holomorphically with the parameter $c\in\WMbrotpq$. 

We shall fix an arbitrary choice of potential $l_0=1$. 
For $n\in\N$ we define the dynamical sets $V_n^c := \{ z \in\C\;|\; G_c(z) < l_0/2^n\}$ 
bounded by the $l_0/2^n$ level set $\EE_n^c = \{ z \in\C\;|\; G_c(z) = l_0/2^n\}$. 
And we define the restricted parameter wakes $\W^\Mbrot_n(p/q) := \{c\in\WMbrotpq| c\in V_n^c\}$.

For $c\in\W^\Mbrot_0(p/q)$ we define the Yoccoz puzzle as follows.
Let $\GG\YY_0^c$ denote graph
$$
\GY_0^c= \EE^c_0\cup\{\al_c, \al'_c\}\cup\bigcup_{i=0}^{q-1}
((\RR^c_{\theta_i}\cup\RR^c_{\theta_i+1/2})\cap V^c_0).
$$
That is the union of the equipotential $\EE^c_0$ together with $\al_c$, $\al'_c$
and the segments, inside $\EE_0$,
of the external rays landing on these two points.
(Note that the original construction involved only the cycle of rays landing on $\al_c$
and not the preimages landing on $\al'_c$, that we add here for convenience).

Following up un the remark above we note that the graph $\GY_0^c$ 
move holomorphically with $c\in\W^\Mbrot_0(p/q)$. 

We define recursively the level $n$-Yoccoz graph $\GY^c_{n+1} := Q_c^{-1}(\GY^c_n)$. 

The {\it level-$0$ puzzle pieces} are the bounded connected components of
$\C\Sm\GG\YY_0^c$. Denote by $\YY_0^c$ the level-$0$ puzzle : the collection of these $2q-1$ puzzle pieces.
Define the level-$n\in\N$  puzzle $\YY_n^c$ as the collection
of connected components of $Q_c^{-n}(Y)$,
where $Y$ ranges over all of the level-$0$ puzzle pieces or equivalently 
as the set of bounded connected components of $\C\Sm\GY^c_n$. 
The ($p/q$-Yoccoz) Puzzle for $Q_c$ is
the union $\YY^c=\cup_{n\geq 0}\YY_n^c$ of the puzzles at all levels. 
We shall also use the finite unions of puzzles $\YY^c(N) := \cup_{0\leq n\leq N} \YY_n^c$, $N\in\N$.

Denote by  $\GY^c(n)$ the union $\bigcup_{j= 0}^n \GY^c_j$ and let 
$\GY^c$ be the union of these graphs of all levels.
       
\begin{figure}[h]
  \begin{center}\includegraphics[height=3.1cm]{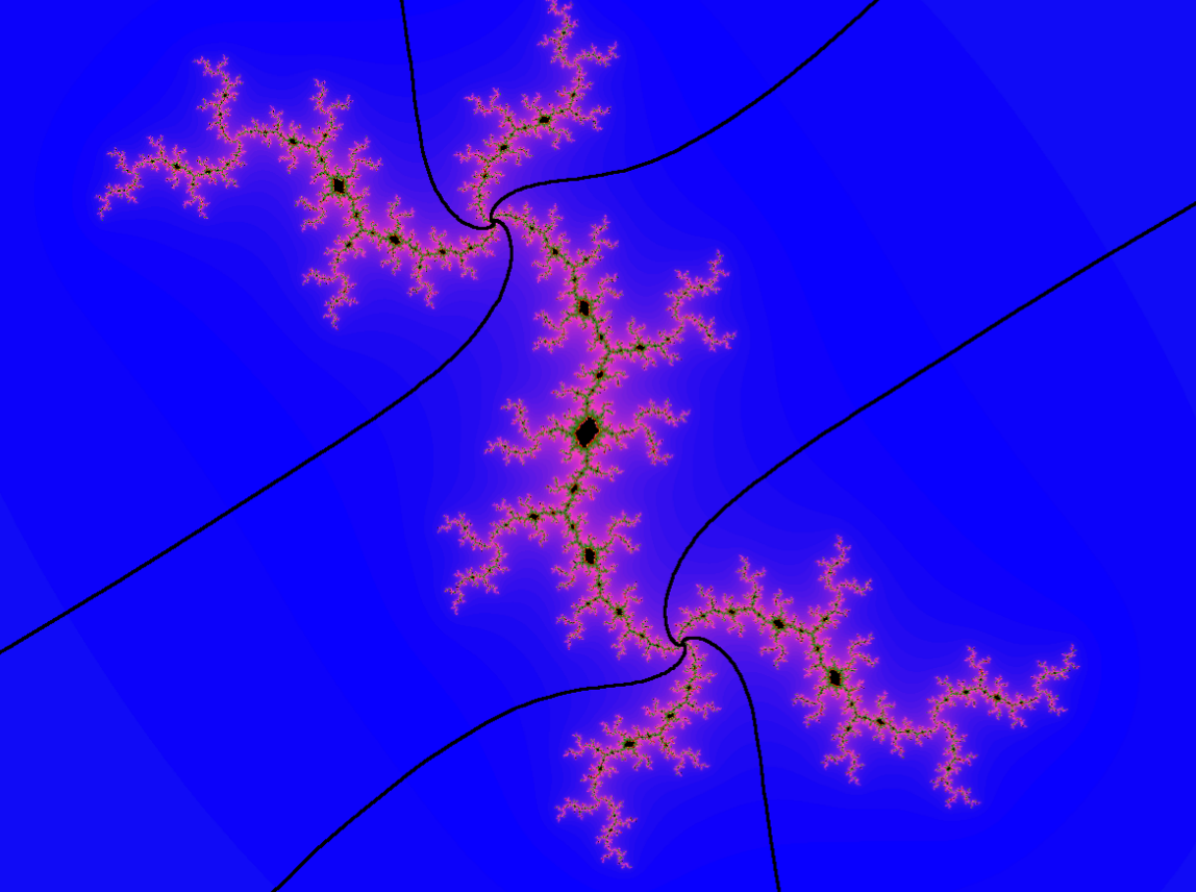}\ \includegraphics[height=3.1cm]{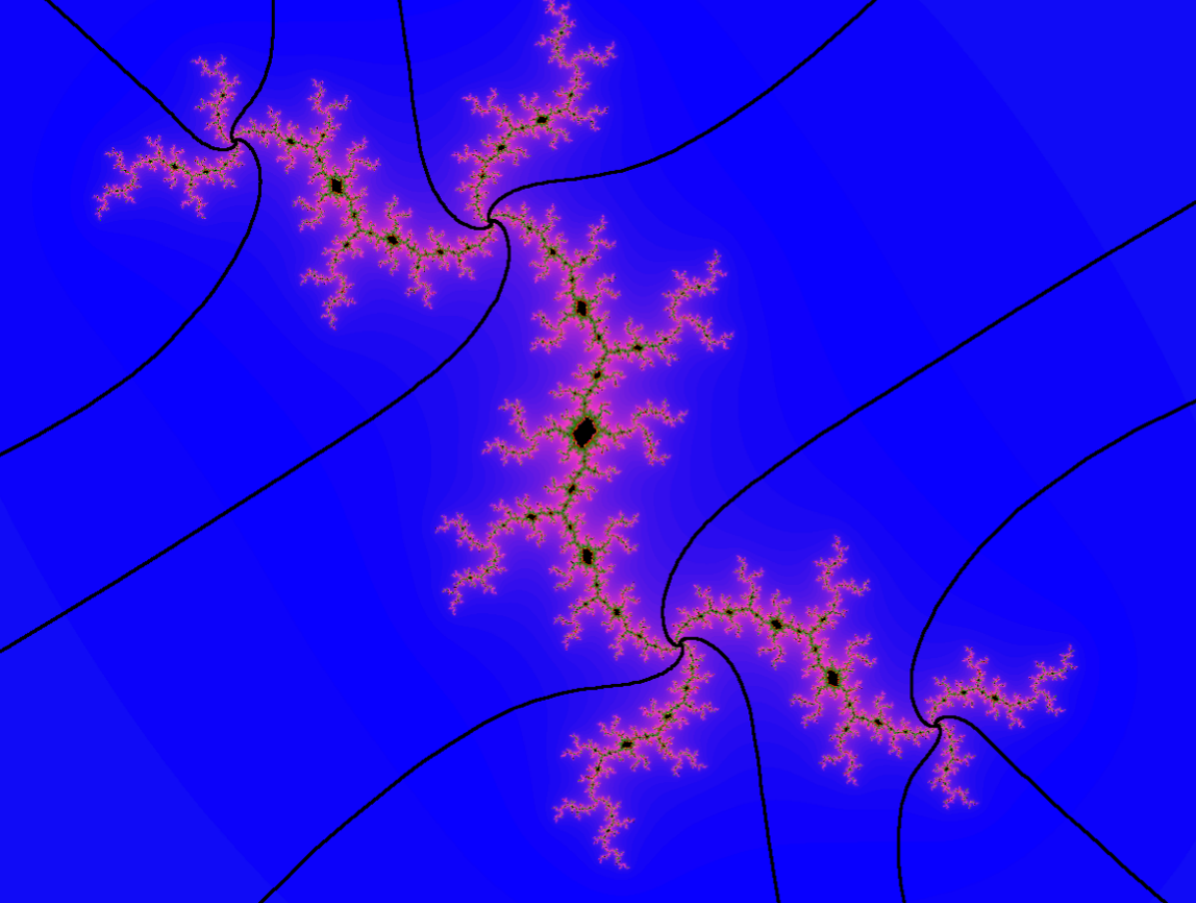} \ \includegraphics[height=3.1cm]{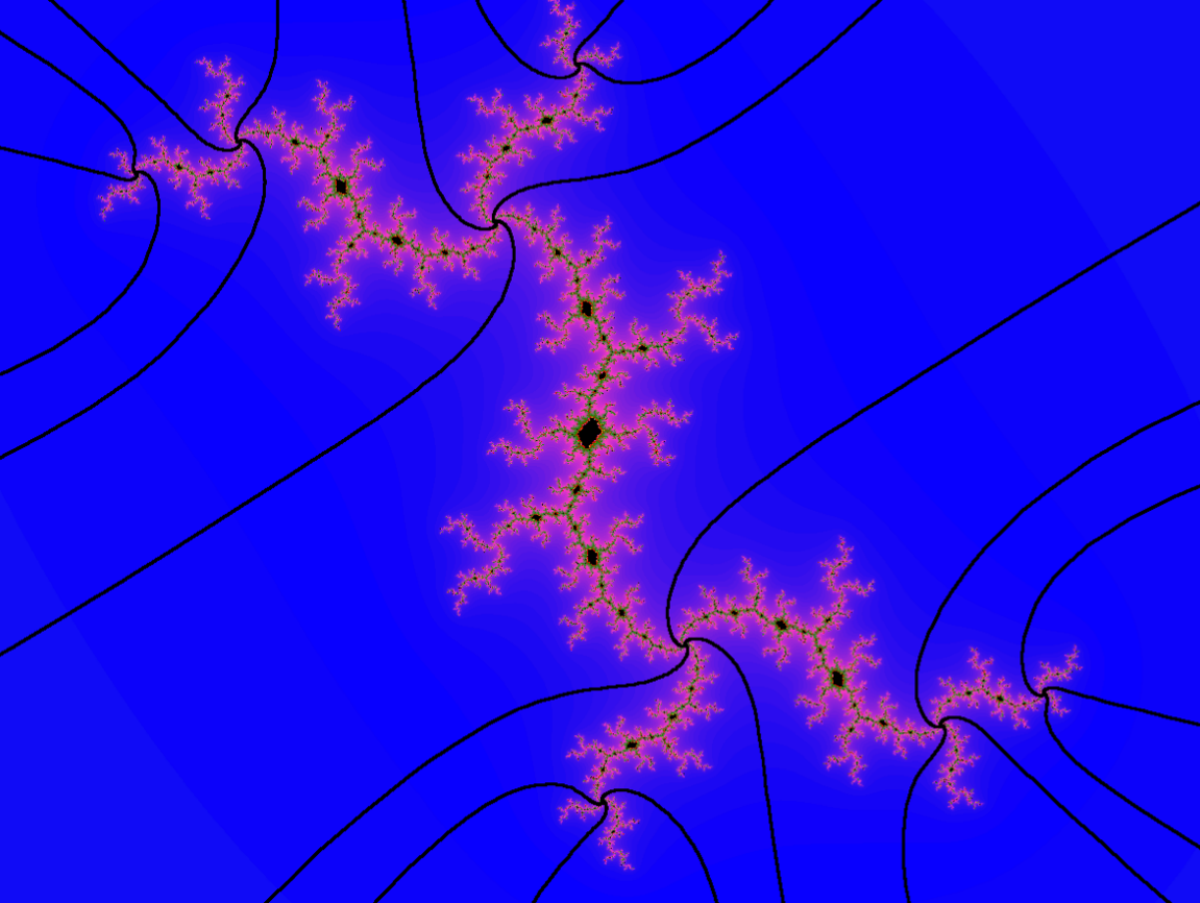}\end{center}
\caption{ Yoccoz dynamical puzzles (without equipotential) in the wake $1/3$}
 \end{figure}

Any two puzzle pieces $Y\in\YY_n^c$ and $Y'\in\YY_m^c$, $m\leq n$
are either interiorly disjoint or nested with $Y\subseteq Y'$ (because the
potential is multiplied by two under the dynamics and
the set of rays in the construction of $\YY^c_0$ is forward invariant).

A \emph{nest}, i.e.~a sequence
$\NN={\{Y_n\}}_n$, $Y_n\in\YY^c_n$ with $Y_{n+1}\subseteq Y_n$,
is called \emph{convergent} iff $\End(\NN) := \bigcap_{n\in\N} \ov{Y}_n = \{z\}$
a singleton and is called \emph{divergent} otherwise. 
When wanting to emphasise $z$ we say the nest $\NN$ is convergent to $z$. 
A nest $\NN$ is called \emph{critical} iff $0\in\End(\NN)$ and
called a \emph{critical value nest} iff $c\in\End(\NN)$.

%%%%%%%%%%%%%%%%%%%%%%%%%%%%%%%%%%%%%%%%%%%%
  {\bf Universal Yoccoz Puzzle} 

Associated to a rotation number $p/q$, one defines the universal Yoccoz Puzzle  on the complement of the disk. It is a model of all $p/q$ Yoccoz Puzzle 
using the B\"ottcher conjugacy.
Let $\displaystyle \ZZ_{0}=\bigcup_{j=0}^{q-1} \e^{i2\pi\theta_j}\cup
\e^{i2\pi(\theta_j+\frac 1 2)}$ be  the unique $q$-cycle for $Q_0$
of combinatorial rotation number $p/q$ in $\Sen$ and its preimage.

Let $l_0$ denote the choice of equipotential above and define $E_0=\{z\mid |z|=e^{l_0}\}$.
Let $\GU_0$ denote the union of the equipotential $E_0$, the unit circle, together with the segments, 
of radial lines through the points of $\ZZ_0$ between $E_0$ and the unit circle.
 
The Universal  Yoccoz graph is then 
 $$
\GU_0 = E_0\cup \Sen\cup\left(\bigcup_{i=0}^{q-1}
\RR_{\theta_i}
\cup\bigcup_{i=0}^{q-1}
\RR_{\theta_i+1/2})\right)\cap  \{z\in \C\mid 1<|z|<e^{l_0}\}
$$
and define the universal ($p/q$-Yoccoz) puzzle $\UY_0$ as the set consisting of the $2q$ 
bounded connected components of the complement of $\GU_0$ in $\C\setminus\D$.

Define $\GU_n$ recursively as follows\,:
$$ \GU_n=Q_0^{-1}(\GU_{n-1})$$ 
and the universal ($p/q$-Yoccoz) puzzle $\UY_n$ as the set
consisting of the  bounded connected components of the complement of $\GU_n$ in $\C\setminus\D$.
 
Finally define $\UY(n) = \cup_{j\leq n} \UY_j$ and $\UY= \cup_{n\in\N} \UY_n$.
We call $\UY$ \emph{the universal $p/q$-Yoccoz puzzle}.
Remark that if $\phi_c(c)\notin \UY(n)$, in particular if $c\in L^\star_{p/q}$, 
then 
$$
\GY^c(n) = \ov{\psi_c(\UY(n))}.
$$
%%%%%%%%%%%%%%%%%%%%%
\subsection{Parabolic dynamical puzzle}
Similarly to the polynomial case above, we have for each irreducible rational $p/q$ defined in \cite{PR2} 
a universal parabolic $p/q$-puzzle using the parabolic rays described in section~\ref{s:parabolicrays} above. 
We shall for completeness briefly review the construction here
The universal parabolic $p/q$-puzzle is the puzzle for the Blaschke product $\Bla$, 
which is the model map for the external class of the maps $g\in Per_1(1)$. 
Recall that the notation $R_\veps$ refers to the \emph{external parabolic ray} of $\Bla$ 
with argument $\veps\in\Sigma^2$, as defined in section~\ref{s:parabolicrays}.
Note that with the parabolics, the difficulty is that we do not have equipotentials. Therefore 
we will define the shortest path from one ray to the next (in the cycle) and use this path 
as an equipotential. In \cite{PR2} we compared the two universal puzzles and 
showed that there is a natural dynamics preserving bijection between the two 
puzzles. In fact there is a modified Universal Yoccoz puzzles inducing Yoccoz puzzles similarly as above,  
puzzles which yield the same puzzle results as standard Yoccoz puzzles associated with the Universal Yoccoz puzzle $\UY$ and such that the modified Universal Yoccoz puzzle is homeomorphic to the Universal parabolic Yoccoz puzzle $\PP$ to be introduced below.
 
\subsubsection{Shortcuts}\label{shortcuts}
Recall that the  restriction  {\mapfromto {\phi^+} {D_0} {\C\setminus ]-\infty,0]}} of $\phi$
is a conformal isomorphism which extends continuously to the boundary. 
For $0< n_0, n_1$ let $\ckga(n_0, n_1)$ be the arc which
is mapped by $\phi^+$ to the Archimedean spiral/circle of center $0$ connecting $-n_0$ 
and $-n_1$ counter clockwise through $\C\Sm\Rminus$. 
And let $\ga(n_0,n_1) = -\ckga(n_0, n_1)\subset D_{1/2}$. 
Since any branch of $\Bla^{-n}$ for any $n\geq 1$ is univalent on $D_{1/2}$, 
we may use such branches to define short-cuts in any of the pre-images $D_{r/2^n}$ of $D_{1/2}$
under $\Bla^{n-1}$, where $n\geq 1$ and $r$ is the odd number such that $h(\exp( i2\pi r/2^n))$ 
belongs to the boundary of $D_{r/2^n}$. Short-cuts were introduced in \cite{PR2} in order to 
produces parabolic Yoccoz graphs and puzzles, which are topologically similar to 
standard polynomial Yoccoz graphs and puzzles. 

The basic observation is that if $\veps^0\in\Si_2$ has $n_0 > 1$ leading $0$'s followed by a $1$
and $\veps^1\in\Si_2$ has $n_1 > 1$ leading $1$'s followed by a $0$. 
Then the two rays $R_{\veps^0}$ and $ R_{\veps^1}$ follow the boundary of the disk $D_0$ 
precisely down to times $-n_0$ and $-n_1$ respectively. 
When forming e.g. puzzles where the two rays $R_{\veps^0}, R_{\veps^1}$ are adjacent and 
so are destined to bound a puzzle piece we shall replace the subarc of $R_{\veps^0}\cup R_{\veps^1}$ 
between $R_{\veps^0}(-n_0)$ and $R_{\veps^1}(-n_1)$ with the short-cut $\ckga(n_0, n_1)$. 
Similarly if $\veps^0$ has a leading $1$, followed by $n_0-1$  digits $0$ with  $n_0-1\geq 1$ and then a $1$ and 
$\veps^1$ has a leading $0$ followed by $n_1-1\geq 1$ leading $1$'s and then a $0$. 
Then the two rays $R_{\veps^0}$ and $ R_{\veps^1}$ follow the boundary of the disk $D_{1/2}$ 
precisely down to times $-n_0$ and $-n_1$ respectively. 
In this case when  the two rays $R_{\veps^0}, R_{\veps^1}$ are adjacent in a graph or puzzle 
we shall replace the subarc of $R_{\veps^0}\cup R_{\veps^1}$ 
between $R_{\veps^0}(-n_0)$ and $R_{\veps^1}(-n_1)$ with the short-cut $\ga(n_0, n_1)$. 
And finally if $\veps^0, \veps^1$ coincide up to digit $n-1$, but differ on the $n$-th digit, say 
 $\si_2(\veps^0)$ has a leading $1$, followed by $n_0-1\geq 1$ $0$'s and then a $1$ and 
$\si_2(\veps^1)$ has a leading $0$ followed by $n_1-1\geq 1$ leading $1$'s and then a $0$. 
Then the two rays $R_{\veps^0}$ and $ R_{\veps^1}$ coincide down to time $n-1$ and 
follow the boundary of the disk $D_{r/2^n}$ 
precisely down to times $-n_0-n$ and $-n_1-n$ respectively, where $r$ 
has the binary representation given by the first $n$ digits of $\veps^0$. 
Similarly to the above we can short-cut $R_{\veps^0}\cup R_{\veps^1}$ through $D_{r/2^n}$. 

 \begin{figure}[h]\label{fshortcut}
\begin{center}
\includegraphics[width=0.5\textwidth]
{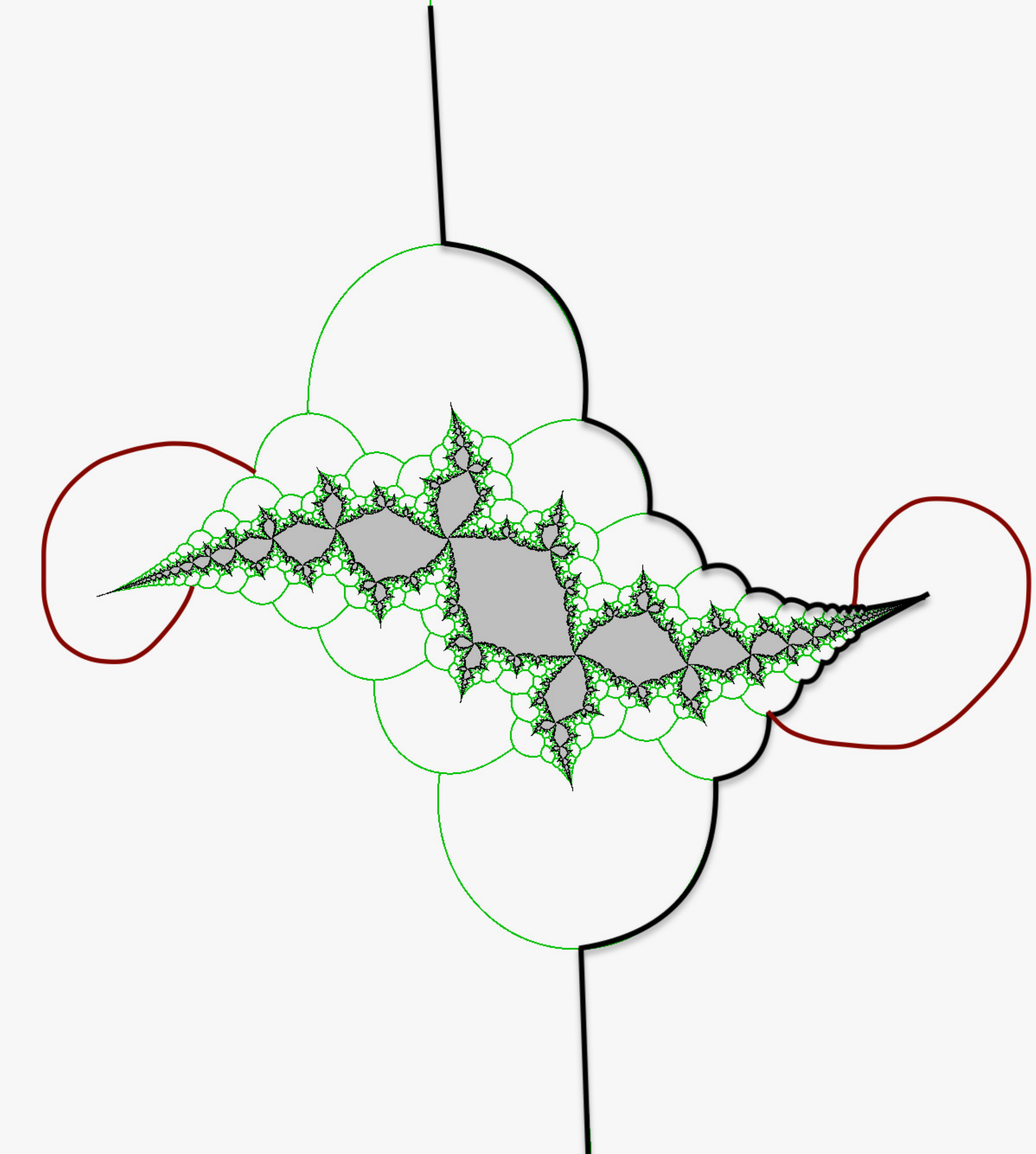}
 \end{center}
\caption{Example of shortcut on the rays $\RR^B_{\overline 0}$ and $\RR^B_{\overline 1}$ : on the right shortcut  $\ckga(n, 3)$ and on the left $\ga(n, 3)$.}
\end{figure}

In any of the three cases we denote by $\widehat\gamma(\veps^0,\veps^1)$ 
the arc obtained from $R_{\veps^0}\cup R_{\veps^1}$ by short-cutting through the appropriate 
$D_{r/2^n}$.

\subsubsection{The Universal Parabolic $p/q$ Yoccoz Puzzle.}
As above let $\displaystyle \ZZ_{0}=\bigcup_{j=0}^{q-1} \e^{i2\pi\theta_j}\cup
\e^{i2\pi(\theta_j+\frac 1 2)}$ be  the unique $q$-cycle for $Q_0$
of combinatorial rotation number $p/q$ in $\Sen$ and its preimage. Then the set  $h(\ZZ_0)$ corresponds to  the unique $p/q$ orbit of $\Bla$
together with its preimage under $\Bla$ (recall that $h$ is the
 conjugacy between $\Bla$ and $z^2$ satisfying $h(z^2)=\Bla\circ h$ defined in Lemma~\ref{l:h}). 
Let
$\overline{0}< \veps_0 < \veps_1 < \ldots < \veps_{2q-1}< \overline{1}$
denote the unique itineraries of these points and
let $\GP_0$ denote the graph
$$
\GP_0 = \Sen\cup\bigcup_{i=0}^{2q-1}
\widehat\gamma(\veps_i,\veps_{(i+1)\mod 2q})
$$
and define the universal parabolic ($p/q$-Yoccoz) puzzle $\PP_0$ as the set
consisting of the $2q$ bounded connected components
of the complement of $\GP_0$ in $\C\setminus\overline{\D}$. 
Denote by $P_{1,0}$ and $P_{-1,0}$ the puzzle pieces with 
$1$ and $-1$ respectively on the boundary. 
Then by construction all the level $0$ puzzle pieces except $P_{1,0}$ 
are pre-images of $P_{-1,0}$ under some iterate of $\Bla^k$, $0\leq k \leq q$. 
This is different from the level $0$ universal Yoccoz-puzzle, where all puzzle pieces 
are bounded by the same equipotential. 

Define $\PP_n$ recursively as follows\,:
$$ \PP_n=\{\Bla^{-1}(P)\mid P\in\PP_{n-1}, 1\notin \partial P\}\cup \{P_{1,n}, P_{-1,n}\},$$
where  $P_{1,n}$, resp. $ P_{-1,n}$,
is  the component bounded by
$$
\widehat\gamma( (\underset{n\text{
times}}{\underbrace{0,\ldots,0}}\,,\veps_0), (\underset{n\text{
times}}{\underbrace{1,\ldots,1}}\,,\veps_{2q-1}))
\quad\text{resp. by }\quad \widehat\gamma(
(0,\negthickspace\underset{(n-1)\text{
times}}{\underbrace{1,\ldots,1}} \negthickspace,\veps_{2q-1}),
(1,\negthickspace\underset{(n-1)\text{
times}}{\underbrace{0,\ldots,0}} \negthickspace,\veps_{0}))
$$
together with the corresponding arc on the unit circle. 

We shall write $\ckga_n$ for the short-cut $\partial P_{1,n}\cap D_0$ and 
$\ga_n$ for the short-cut $\partial P_{-1,n}\cap D_{1/2}$.

By construction the only non dynamical parts of the universal parabolic $p/q$ graph and puzzle are the short-cuts 
$\ckga_n$ and $\ga_n$, $n\geq 0$, i.e.
$$
\Bla(\GP_{n+1}\Sm(\ckga_{n+1}\cup\ga_{n+1})) = \GP_n\Sm\ckga_n.
$$ 

Finally define $\PP = \cup_{n\in\N} \PP_n$.
We call $\PP$ \emph{the (quadratic) universal parabolic $p/q$-Yoccoz puzzle}. 
Denote by $\GP_n=\bigcup_{P\in\PP_n}\partial P$.
Denote by  $\GP(n)$ the union $\bigcup_{0\le k \le n}\GP_k$; it coincides with  the union
of the boundaries of puzzle pieces of all levels up to and including $n$.
Let  $\GP$ be the union of  these graphs of all levels. 
 \begin{figure}[h]\label{fshortcutgraf}
\begin{center}
 \includegraphics[height=2 in, angle=0]{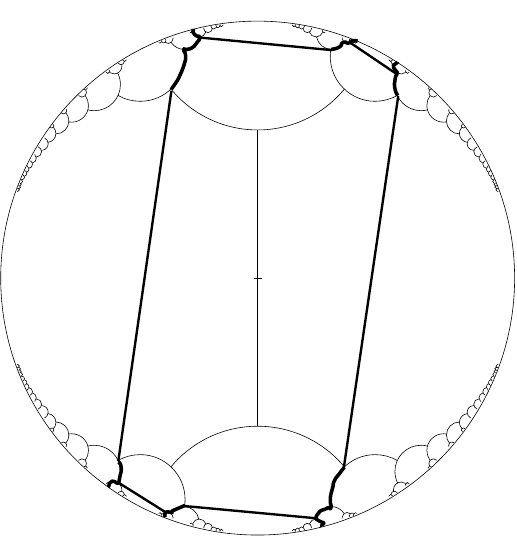}   \includegraphics[height=2 in, angle=0]{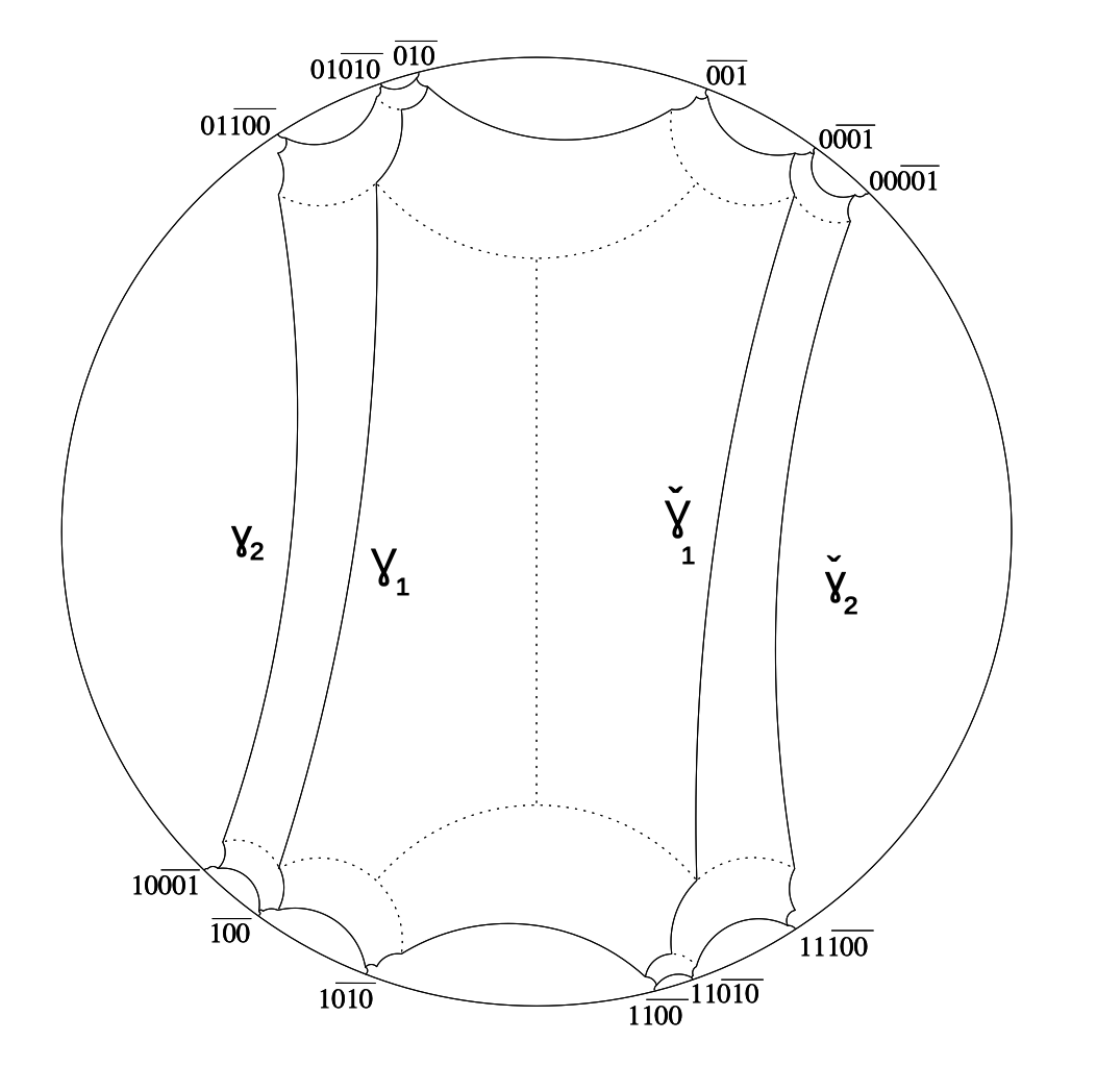}
\end{center}
\caption{Shortcuts In the model, graphs of level 0 and 1.}
\end{figure}

\begin{figure}[h]\label{fpuzzlesmodel}
\begin{center}  \includegraphics[height=2 in]{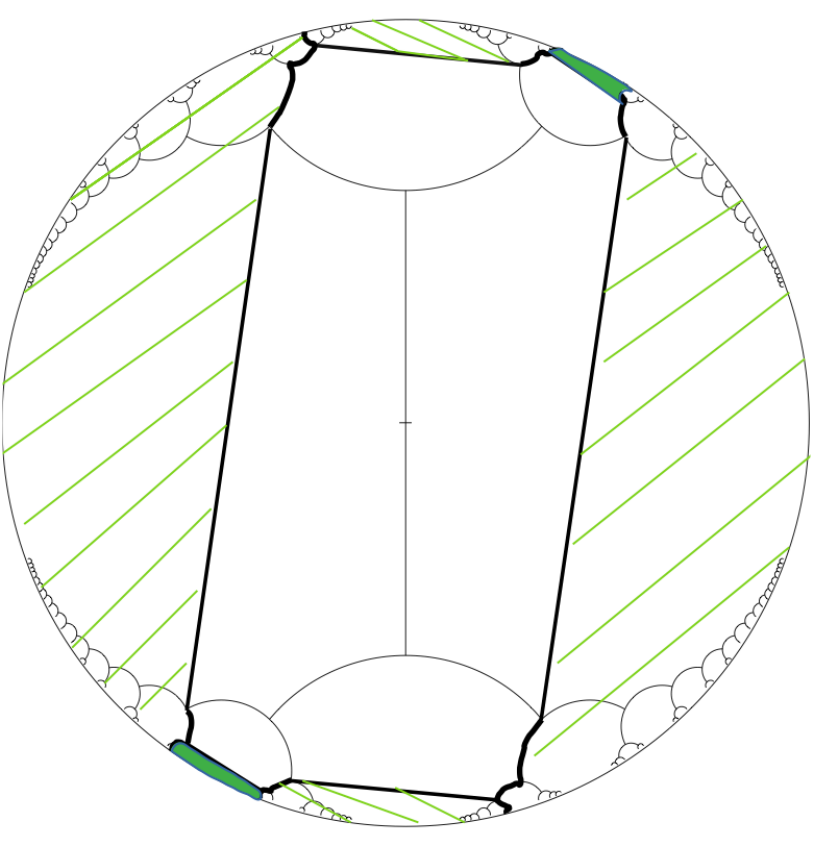}\ \includegraphics[height=2 in]{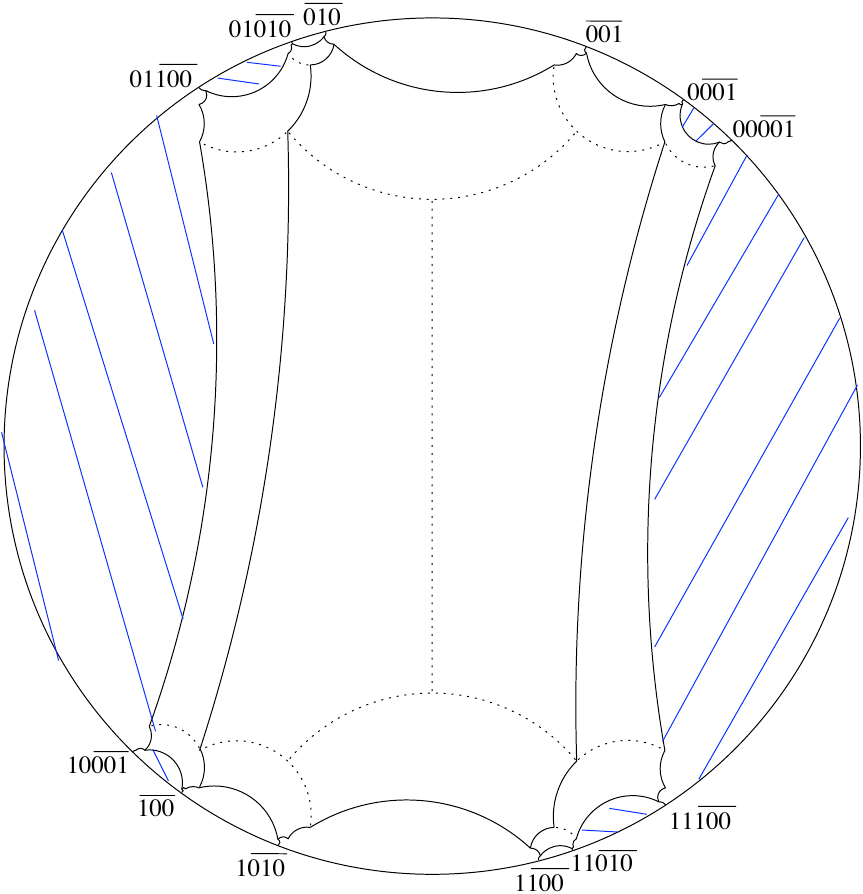}
\end{center}
\caption{Puzzles pieces of level 0 and  of level 1, the ones with color define non degenerate annuli.}
\end{figure}

For every $p/q$, there is a correspondence between the Universal Yoccoz puzzle 
and Universal Parabolic puzzle. 
For any universal Yoccoz puzzle piece of depth n bounded by external rays of argument 
$\{t_1,t_2\}$, the corresponding universal parabolic puzzle piece of depth $n$ 
is bounded by the parabolic rays of argument $\{h(t_1),h(t_2)\}$.

\subsubsection{Parabolic $p/q$-Puzzle }\label{s:parabopuzzle}
Let $p/q$ be an irreducible rational and let $B\in\W^{\Mone}(p/q)$ 
(for the definition of the $p/q$ wake $\W^{\Mone}(p/q)$ see Lemma~\ref{wake}).  
The parabolic $p/q$ puzzle for $g_B$ is derived from the universal parabolic $p/q$ puzzle, 
in a maner similar to how the Yoccoz-puzzle for $Q_c$, $c\in\W^\Mbrot(p/q)$ is derived 
from the universal Yoccoz puzzle. 

For each $n\geq 0$ let $V_n^\PP$ be the interior of the union of closures of level $n$ 
universal parabolic puzzle pieces. We define reduced wakes
$$
\W^\Mone_n(p/q) := \LL^\Mone_{p/q}\cup \{B\in\WMonepq | h_B(v_B)\in V_n^\PP\}, 
$$
which are similar to the reduced wakes $\W^\Mbrot_n(p/q)$ though the phrasing of the definition is different.

Recall from \corref{vlocation} that for $B\in\W^{\Mone}(p/q)$, 
the unique $p/q$-cycle of parabolic rays for $g_B$ 
with rotation number $p/q$ co-lands on $\al_B$ the unique finite fixed point for $g_B$, 
which is repelling. 
And moreover the second critical value $v_B$ for $g_B$ belongs to the dynamical wake 
$\W^B(p/q)$.
If $B\in\W^{\Mone}(p/q)\Sm\Mone$ then we may extend $h_B^{-1}$ 
analytically to a univalent map on a neighbourhood in $\Chat\Sm\Dbar$ 
of $(\GP^0\Sm\Dbar)\cup D_1 \cup U_{p/q}$, where $U_{p/q}$ is the unbounded connected 
component of $\Chat\Sm\GP^0$.
And if $B\in\Mone\cap\W^{\Mone}(p/q)$ then $h_B^{-1}$ even 
extends to a biholomorphic map {\mapfromto {h_B^{-1}} {\Chat\Sm\Dbar} {\La_B}}. 
Thus for every $B\in\WMonepq$ we may define short cuts $\ckga_n^B = h_B^{-1}(\ckga_n)\subset D_0^B$ 
and $\ga_n^B = h_B^{-1}(\ga_n)\subset D_{1/2}^B$.
So we may define parabolic Yoccoz graphs:

\REFDEF{Dynamic_Graphs}
For $B\in\W^\Mone_0(p/q)$ we define the dynamical graph $\GP_0^B$ as 
$$
\GP_0^B = h_B^{-1}(\GP_0)\cup\{\al_B, \al'_B\}.
$$
We define the parabolic ($p/q$-Yoccoz) puzzle $\PP_0^B$ for $g_B$ as the set
consisting of the $2q-1$ connected components
of the complement of $\GP_0^B$ intersecting the Julia set of $g_B$.

Define $\GP_n^B$ recursively by 
$$
\GP_{n+1}^B := g_B^{-1}(\GP_n^B\Sm \ckga_n^B)\cup \{\ckga^B_{n+1}\cup\ga^B_{n+1}).
$$
We denote by  $\GP^B(n)$ the union $\bigcup_{0\le k \le n}\GP_k^B$; 
and we write $\GP^B$ for the union of these graphs of all levels.

And define the parabolic ($p/q$-Yoccoz) puzzle $\PP_n^B$ for $g_B$ as the set
of complementary components $P$ of $\GP_n^B$ 
intersecting the Julia set for $g_B$, $\PP^B(N) := \cup_{0\leq n\leq N} \PP_n^B$ and 
$\PP^B := \cup_{0\leq n} \PP_n^B$. 

We take over the vocabulary from Yoccoz puzzles and write $P_n^B(z)$ 
for the level $n$ puzzle piece containing $z$, if one such exists, 
$\NN = \{P_n\}_n$ for a nest of puzzle pieces, 
$\End(\NN) = \cap_n \ov{P}_n$ for the end of $\NN$. 
And moreover that $\NN$ is convergent to $z$ iff $\End(\NN) = \{z\}$ and 
divergent if $\End(\NN)$ is not a singleton.
\ENDDEF

 \begin{figure}[h]\label{fpuzzlesmodelGB}
\begin{center}  \includegraphics[angle=90, height=4 in]{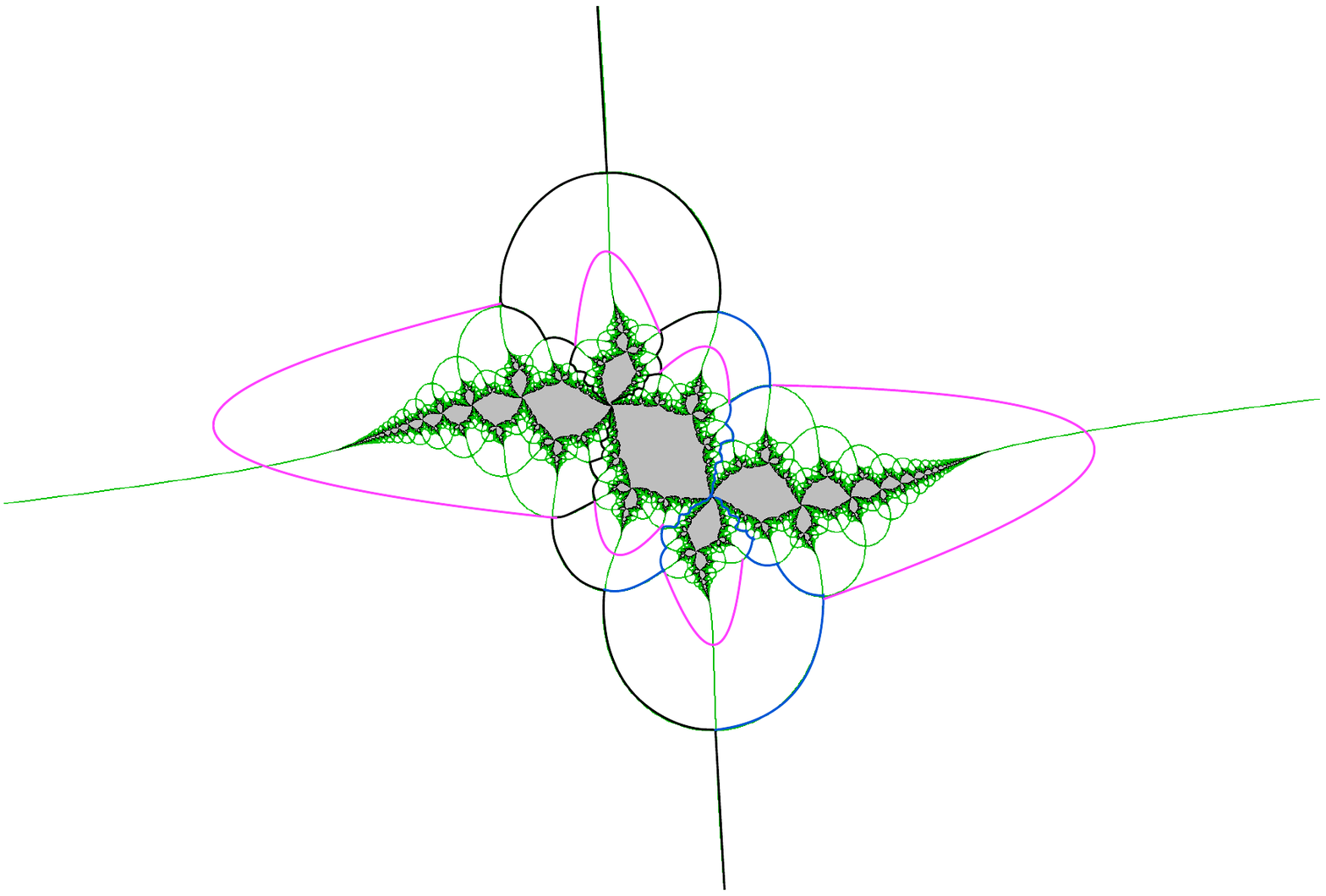}  \end{center}
\caption{Puzzle pieces of level 0.}
\end{figure}

\REFREM{bord}
Note that by construction the second critical point and value, $-1, v_B$ as well as $\be_B$ 
are pairwise separated by $\GP^B_0$ for every $B\in\W^\Mone_0(p/q)$. 
Moreover the dynamical wake $\W_B(p/q)$ is disjoint from $D_0^B$, for any $B\in\W^\Mone(p/q)$ 
so that the second critical point $-1$ does not belong to $D_1^B$ and hence 
for every $n$ the puzzle pieces $P_n(\beta_B)$ and $P_n(\beta'_B)$ are defined and the restrictions 
$$
\mapfromto{g_B}{\partial P_{n+1}(\be_B)\Sm\ckga_{n+1}^B}{\partial P_{n}(\be_B)\Sm\ckga_{n}^B} 
\quad\textrm{and}\quad
\mapfromto{g_B}{\partial P_{n+1}(\be'_B)\Sm\ga_{n+1}^B}{\partial P_{n}(\be_B)\Sm\ckga_{n}^B} 
$$
are diffeomorphisms.
\ENDREM

For $B\in\W^\Mone_0(p/q)$ define the graph
$$
\GP^B_\beta := \GP^B_0\cup\{\be_B, \be'_B\} \cup\bigcup_{n\geq 0} (\partial P_n(\be_B)\cup\partial P_n(\be'_B)).
$$
\REFPROP{base_para_holomorphic_motion}
Let $B_*\in\W^\Mone_0(p/q)$ be arbitrary. Then there is a holomorphic motion 
$$
\mapfromto{\psi^\Bstar_\beta}{\W^\Mone_0(p/q)\times\GP^\Bstar_\beta}\Chat
$$
with base point $\Bstar$ such that $\psi^\Bstar_\be(B, \GP^\Bstar_\beta) = \GP^B_\beta$ and 
$g_B\circ\psi^\Bstar_\be(B, z) = \psi^\Bstar_\be(B, g_\Bstar(z))$ 
for every $B\in\WMonepq$ and $z\in\GP^\Bstar_\beta$.
\ENDPROP
\PROOF
Since the family $g_B$, $B\not=0$ has a persistent parabolic fixed point of fixed parabolic multiplicity, 
the normalized Fatou coordinates $\phi_B$ for $g_B$ depends holomorphically on both $B$ and $z$. 
Hence also the coordinates $h_B$ depends holomorphically on the two variables $(B,z)$. 
For the same price the short-cuts $\ckga_n^B, \ga_n^B\subset \ov{D}_1^B$ move holomorphically with $B$. 
Since $\GP^B_0$ is the closure of $h_B^{-1}(\GP_0\Sm\Dbar)$ the graph $\GP^B_0$ moves 
holomorphically with $B\in\W^\Mone_0(p/q)$. By construction $\be_B\equiv \infty$ and $\be'_B\equiv 0$ 
and so move holomorphically. 
Finally by \remref{bord} above there is no critical point for $g_B$ on any of the boundary arcs 
$\partial P_{n+1}(\be_B)\Sm\ckga_{n+1}^B$, $\partial P_{n+1}(\be'_B)\Sm\ga_{n+1}^B$ for any $n\geq 0$ 
and any $B\in\W^\Mone_0(p/q)$. Hence also these move holomorphically with $B$. 
From this the proof follows.
\ENDPROOF

\section{ Tower of laminations, combinatorial invariants }\label{s:Tower}
\subsection{Abstract Towers}\label{abstract_towers}
The following presentation is an excerpt from \cite{PR1}.
For more details the reader is referred to this paper.

Let $E_0$ denote the unique $p/q$ cycle for $Q_0$, then denote by $\ZZ_n=Q_0^{-(n+1)}(E_0)$ for $n\geq -1$ and
$\ZZ=\cup_{n\geq 0} \ZZ_n$.  Remark that   $\ZZ_0=E_0\cup (-E_0)$ and  that $\ZZ_n=Q_0^{-n}(\ZZ_0)$ for $n\geq 0$.
 
For $E\subset\Sen$ we let $H(E)$ denote $E$ union its hyperbolic convex hull in $\D$, $H(E)$ is therefore a closed set in $\overline \D$.

 \begin{figure}[h]\label{flam}
\begin{center}  \includegraphics[height=2 in]{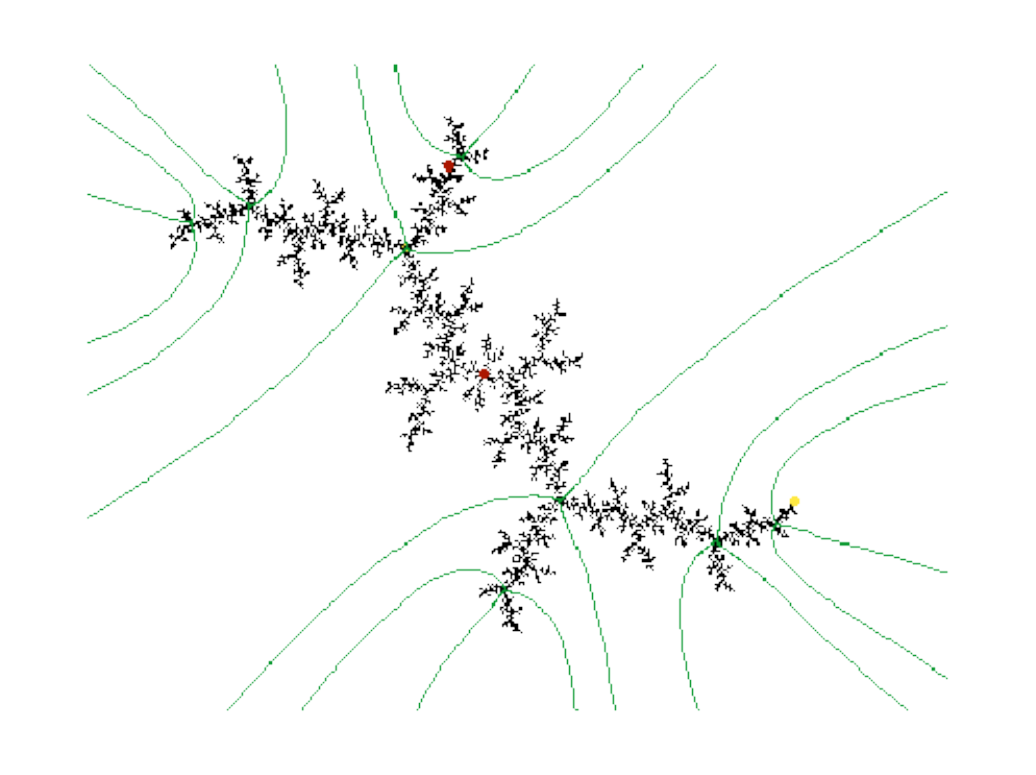}\ \includegraphics[height=2.53 in]{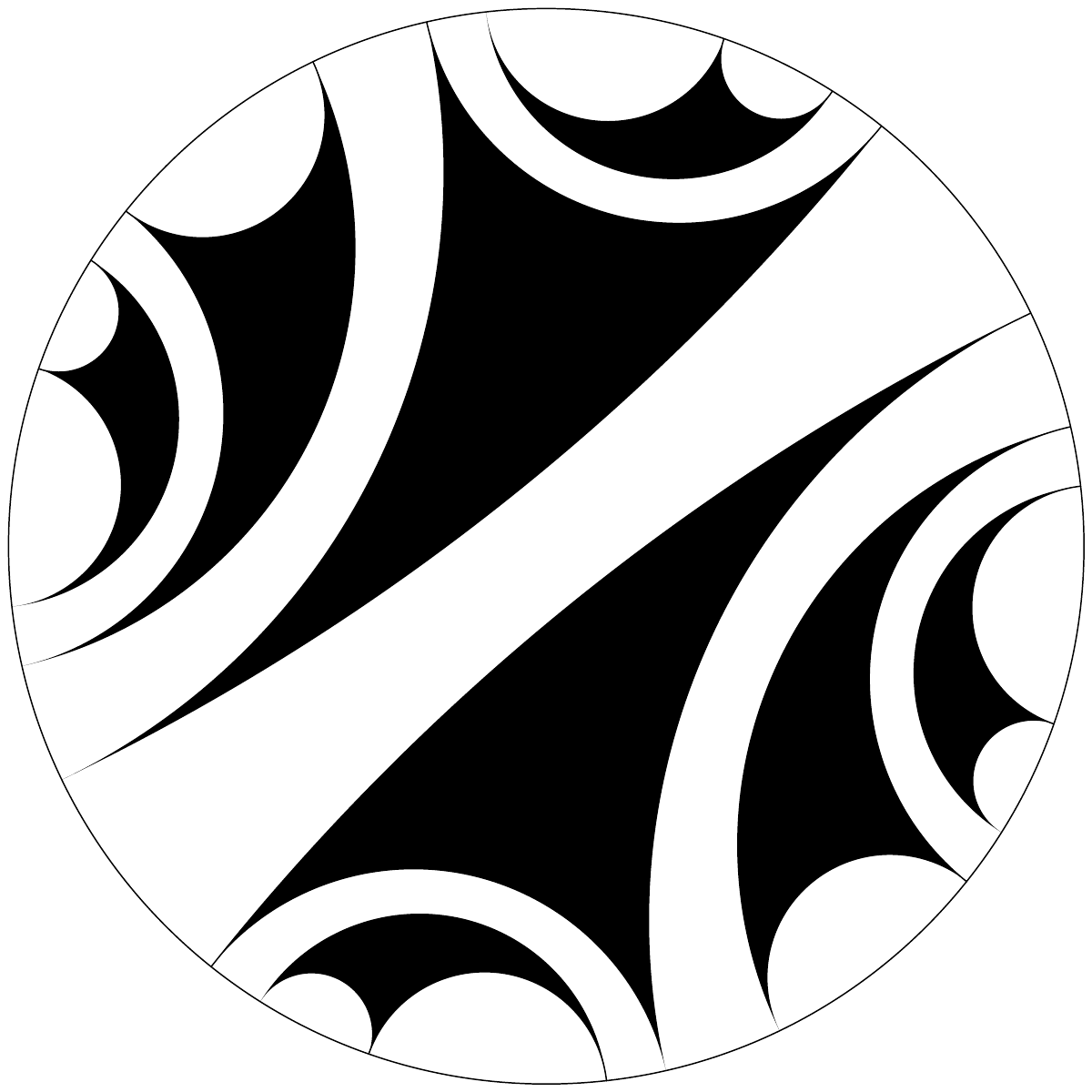}
 \end{center}
\caption{Lamination associate to a Yoccoz graph, the critical point and value in red.}
\end{figure}
 
\DEF
A  tuple of equivalence relations ${(\sim_n)}_{0\leq n \le  N}$, with $N\in\N\cup\{\infty\}$,
is called a tower if it satisfies the following admissibility conditions (see also \cite{Kiwi}):
\ENUMi
\item For each $n$:  $\sim_n$ is an equivalence relation on $\ZZ_n$.
\item $\sim_0$ has the two classes $E_0$ and $-E_0$.
\item  For any class $E$ of $\sim_n$ with $0\leq n\leq N$
the set $Q_0(E)$ is a class of $\sim_{(n-1)}$\,;
\item\label{increasingunion}
$\displaystyle{\sim_N = \bigcup_{n=0}^N \sim_n}$ so that $\displaystyle{\sim_n|_{\ZZ_m} = \; \sim_m}$
for any $m,n$ with  $0\leq m < n\leq N$\,;
\item
For  any two distinct classes
$E$ and $E'$ of $\sim_n$, with $0\leq n\leq N$,
$\displaystyle{H(E)\cap H(E')=\emptyset}.$
\ENDENUM
\ENDDEF
By property \Itemref{increasingunion} $\sim_N$ imposes $\sim_n$ for $n\leq N$.
We shall thus abbreviate and write simply $\sim_N$ for the tower ${(\sim_n)}_{0\leq n \le  N}$.

The \emph{level} of a class $E$ is the minimal $n\geq 0$ for which $E\subset\ZZ_n$.

The finite towers $\sim_N$ are the nodes of a tree with root $\sim_0$
and with a branch connecting each child $\sim_N$ back to its parent $\sim_{N-1}$.
We denote this tree by $\TTT$.
The infinite towers $\sim_\infty$ on the other hand are the infinite branches of this tree
starting at $\sim_0$. We denote the set or space of all infinite branches $\TTT_\infty$.

For a tower $\sim_N$, we denote by the \emph{graph of $\sim_N$} the set
$$
\GG_{\sim_N} = \bigcup_{E \text{ a class of }\sim_N} H(E) \subset \Dbar
$$

A \emph{gap} $G$ of a finite tower $\sim_n$ is any connected component of
$\Dbar\Sm\GG_{\sim_n}$. We denote by \emph{essential boundary} of a gap $G$
the set $\delta G=G\cap\Sen$.
The image of the gap $G_n$ of $\sim_n$ is defined as the gap $G_{n-1}$ of
$\sim_{n-1}$ with $\delta G_{n-1} = Q_0(\delta G_n)$.

A class $E$ or a gap $G$ is said to be \emph{critical} iff $0\in H(E)$, resp. $0\in G$.
Clearly any finite tower has either a (unique) critical class or gap.
We shall denote the critical class/gap of $\sim_n$ by $E^*_n/G^*_n$ (or just
$E^*/G^*$). The image of the critical class or gap will be called the
\emph{critical value class or gap of $\sim_n$} and denoted $E'_n/G'_n$.
Note that the critical value class or gap is a class or gap of $\sim_{n-1}$
and (provided the level of the critical class is $n$)
is a subset of the critical value gap of $\sim_{n-1}$.

For a finite tower $\sim_N$ with critical gap $G^*_N$ 
define the \emph{critical period} $k\geq 1$ of $\sim_N$ as 
the minimal $k\geq 1$ for which $Q_0^k(G^*_N)$ is again a critical gap 
(of $\sim_{N-k}$). 
Note that in fact $k\geq q$ always. 
Also in order to ensure that a critical gap always has a critical period, we may formally define 
$\ZZ_n= E_0$ and $\sim_n$ as the equivalence relation with only one class $E_0$ for any $n$ with 
$-q<n<0$.

Let $\sim_n$ be a finite tower.
If $\sim_n$ has a critical class $E$ it has a unique child
and we say that $\sim_n$ is a \emph{terminal} tower.

If $\sim_n$ has a critical gap with critical value gap $G'_n$
and if $E\subset G'_n$ or $G\subseteq G'_n$
is any class or gap of $\sim_n$ within $G'_n$.
Then $\sim_n$ has a unique extension $\sim_{n+1}$ with critical value class $E$ respective
critical value gap $G$. For this reason we say $\sim_n$ is a \emph{fertile} tower, 
when it has a critical gap.

An infinite tower $\sim_\infty$ is said to be
\emph{renormalizable with combinatorics} $\sim_N$ and
renormalization period $k$ if for every $n\ge N$, $\sim_n$ has critical period $k$
and $N$ is the minimal height with this period.

Suppose $\siminfT$ is an infinite terminal tower with {\it critical value class} $\cvE{n}$,
that is $\sim^T_{n}=\siminf^T|_{\ZZ_{n}}$ has a critical class $ \cE{n}$ with image $\cvE{n}$ and
$\sim_{n-1}=\sim_{\infty}^T|_{\ZZ_{n-1}}$ has a critical gap $\cG{n-1}$
with image the {\it critical value gap} $\cvG{n-1}$ containing $H(\cvE{n})$.
Then $\cvG{n-1}$ contains exactly $q$ gaps $G_n^1, \ldots , G_n^q$ of $\sim_{n-1}$,
which are adjacent to $\cvE{n}$, i.e.~with $H(\cvE{n})\cap\partial G_n^j\not=\emptyset$,
because $H(\cvE{n})$ is a $q$-gon.
In the light of the above discussion let $\sim_n^j$ denote the unique extensions of $\sim_{n-1}$
with critical value gaps $G_n^j$ for $j=1, \ldots , q$.
Define recursively for $m>n$ unique extensions $\sim_m^j$ of $\sim_{m-1}^j$ with critical
value gap $G_m^j\subset G_{m-1}^j$ adjacent to $\cvE{n}$.
Finally denote by $\siminf^j = \cup_{m\geq n} \sim_m^j$
the corresponding infinite towers for $j=1, \ldots , q$.

We shall say that $\siminfT$ is \emph{adjacent} to any of the $q$ towers
$\siminf^1, \ldots ,\siminf^q$ and vice versa. 

\subsection{The natural tower puzzle relation}
Recall that if $c\in \Mbrot\setminus\Card$ then $c$ belongs to a limb $L_{p/q}$, 
so that the $p/q$ cycle of rays co-land at the $\alpha$ fixed point. 
Fix $p/q$ with $(p,q)=1$ and let $\ZZ_n, \ZZ$ be given by $p/q$ as in \subsecref{abstract_towers}.
\REFDEF{Q_c-Puzzle-Tower-relation}
Let $c$ belong to the limb $L_{p/q}$.
Define $t, t'\in\ZZ_n$, $n\in\N\cup\{\infty\}$ to be equivalent, $t\sim_n^c t'$ 
if and only if the rays $\RR_t^c$ and $\RR_{t'}^c$ co-land. 
And define $\sim^c$ as the corresponding tower of equivalence relations.
\ENDDEF

Note that the arguments $t\in\ZZ_n$ are precisely the arguments of external rays in the level $n$ Yoccoz graph and puzzle. 
And moreover $t\sim_n^c t'$ if and only if $\RR_t^c$ and $\RR_{t'}^c$ co-land on a point of $Q_c^{-n}(\al_c)$.
We can thus view $\sim_n^c$ as an abstract version of the Yoccoz graph/puzzle, 
where the gaps $G$ corresponds to level $n$ puzzle pieces and the hulls $H(E)$ of classes $E$ 
corresponds to the unions of segments of co-landing rays. In view of this for any class $E$ of $\sim^c$ 
say of level $n$, we shall refer to the closure of the complete set of co-landing rays with arguments in $E$ as $\RR_E$. 
Then the Yoccoz graph $\GY_n$ is also the first graph containing the lower ends of the rays in $\RR_E$. 

Similarly any $g\in \Mone\setminus \overline \D$ belongs to a Limb $\LL^{\Mone}_{p/q}$.
\REFDEF{g_B-Puzzle-Tower-relation}
Let $g$ belong to the Limb $\LL^{\Mone}_{p/q}$. 
Define  $t, t'\in\ZZ_n$, $n\in\N\cup\{\infty\}$ to be equivalent, $t\sim_n^B t'$ 
if and only if the rays $R_{h(t)}^B$ and $R_{h(t')}^B$ co-land.  
And define $\sim^B$ as the corresponding tower of equivalence relations.
\ENDDEF

\REFREM{sim_ext_to_compl}
In both the polynomial and the parabolic case we shall extend the definition of $\sim_n^c$ respectively $\sim_n^B$ 
to maps $Q_c$ with $c\in\W^\Mbrot_{n-1}(p/q)$ with $c\notin\GY^c_{n-1}$ respectively 
maps $g_B$ with $B\in\W^\Mone_{n-1}(p/q)$ with $v_B\notin\GP^B_{n-1}$.
\ENDREM

Recall that {\mapfromto{h}{\Sen}{\Sen}} is the conjugacy between $Q_0$ and $\Bla$ (see Lemma~\ref{l:h}).

In other words $t\sim_n^B t'$ are equivalent if and only if the
parabolic rays of the level $n$ parabolic Yoccoz graph $\GP^n_g$ corresponding to $h(u)$ and $h(v)$ co-land at
the samepoint of $g^{-n}(\al(g))$. 
As in the polynomial case the equivalence relations $\sim_n^B$ can be viewed as abstract parabolic Yoccoz graphs/puzzles 
with the gaps $G$ corresponds to level $n$ puzzle pieces and the hulls $H(E)$ of classes $E$ 
corresponds to the unions of segments co-landing rays. 

We proved in \cite[Lemma 5.7]{PR2} that among all $p/q$ towers adjacent to some terminal tower $\siminfT$ only the renormalizable tower $\siminfstar(p/q)$ with renormalization period $q$ is realised as $\siminfc$ for some $c\in\Mbrot$ and 
similarly as $\sim_\infty^B$ for some $g\in\Mone$.
Moreover we proved in \cite[Theorem 5.6]{PR2} that any other infinite tower is realized 
as $\siminfc$ for some $c\in\Mbrot$. We summarize this as the following Theorem
\REFTHM{quadratic_realisation}
Let $B\in\W^\Mone(p/q)$, then there exists $c\in\W^\Mbrot(p/q)$ with 
$$\sim_\infty^c = \sim_\infty^B.$$
\ENDTHM

\REFDEF{abstract_g}
Define an abstract map $\ckg$ from $\PP_B$ to itself given by 
$\ckg(P_n(\beta')) = \ckg(P_n(\beta)) = P_{n-1}(\beta)$ for every $n\ge 1$ and $\ckg(P_n) = g(P_n)$ 
for every other level $n\ge 1$ puzzle piece $P_n$.
\ENDDEF
The following Proposition says that if $g_B$ with $B\in\Mone\Sm\Dbar$ and $Q_c$ with $c\in\Mbrot\Sm\Card$ 
define the same infinite tower $\siminf$, then their puzzles are similar: 

\REFPROP{samebehaviour}
Let $p/q$ be an irreducible rational, let $B \in \LL_{p/q}^\Mone$ and $c\in L_{p/q}$ be parameters
such that $\sim_N^c=\sim^B_N$, $N\in\N\cup\{\infty\}$.
Then there is a dynamically defined bijection {\mapfromto {\Xi_B} {\YY^c(N)}{\PP^B(N)}}
between the Yoccoz puzzle $\YY^c(N)$ for $Q_c$ and the parabolic Yoccoz puzzle  $\PP^B(N)$ for $g_B$ 
such that
\ENUM
\item\label{puzzle_dynamics} 
For any puzzle piece $Y_n\in\YY^c$ of level $n, 1\le n\le N$ the puzzle piece $P_n = \Xi_B(Y_n)$ 
also has level $n$ and
$\ckg\circ\Xi_B(Y_n) = \Xi_B\circ Q_c(Y_n)$.
\item 
In particular critical puzzle pieces correspond to critical puzzle pieces.
\item Any annulus of the parabolic puzzle $\PP_B$ is non degenerate
if and only if the corresponding annulus in the Yoccoz puzzle $\YY_c$ is non degenerate.
\ENDENUM
\ENDPROP
\PROOF
Combine \thmref{quadratic_realisation} with \cite[Lemma 5.1]{PR2} for the proof in the case $N=\infty$. 
The easier proof in the finite case follows the same line.
\ENDPROOF

Note that if $Y$ in the Propostition above contains either $\beta(c)$ or $\beta'(c) = -\beta(c)$, 
then $P=\Xi_B(Y)$ contains the parabolic fixed point $\infty$ 
respectively its preimage $0$.

\section{Transfering Yoccoz results to maps in $\Mone$}\label{s:Transfer_to_Mone}
We first recall the basic steps in the proof of Yoccoz theorem of local connectivity for quadratic polynomials. 
Then we show that a similar proof can be made for maps $g_B$ in $\Mone$. 
In the following chapters we use this setup to transfer results to the parameter spaces.
Fix for the rest of this section an arbitrary irreducible rational $p/q$. 

\subsection{Basic Yocooz puzzle theory and  estimates}
In this subsection we set-up the machinery for the proof of Yoccoz theorem on local connectivity of the Mandelbrot set at any non renormalizable parameter in any de-rooted limb $L^\star_{p/q}$ 
and for the same price local connectivity of the Julia set of such polynomials.  (see for instance in~\cite{Milnorlc}). 
Recall that $l_0$ was the equipotential level in the definition of Yoccoz puzzles, 
that $V_n^c := \{ z \in\C\;|\; G_c(z) < l_0/2^n\}$ is the dynamical set bounded by the $l_0/2^n$ level set. 
And moreover  $\W^\Mbrot _n:= \{ c \in\WMbrotpq\;|\; c\in V^c_n\}$, $n\geq 0$.

Note that for $c\in\W^\Mbrot_0$ the set $V_n^c$ is the interior of the union of closures of all level $n$ puzzle pieces.

Let $c\in\W^\Mbrot_0(p/q)$.
For $Y_n$ a puzzle piece of some level $n$ and $z\in Y_n$ we write $Y_n(z) :=Y_n$. 
We shall furthermore use the abbreviations $Y_n^0 := Y_n(0)\in\YY_n$
and $Y_n^c:=Y_n(c) = Q_c(Y_{n+1})\in\YY_n$ whenever there is such a puzzle piece, 
that is whenever $Q_c^n(c)$ belongs to a level $0$ puzzle piece. 
If $K_c$ is connected, this only fails whenever $Q_c^{n+1}(c) = \al_c$.
Since the dynamics is quadratic any non-critical puzzle piece $Y$ has a unique dynamical twin $\ckY$ of the same level with 
$Q_c(\ckY) = Q_c(Y)$, in fact $\ckY = -Y$. 
And the critical puzzle pieces are siamese twins in the sence that $Y_n^0 = \ckY_n^0$. 

\REFLEM{q_renormalization}
Let $c\in\W^\Mbrot_{q-1}(p/q)$. 
Then the map $Q_c^q$ has a quadratic-like restriction 
{\mapfromto {f_c  := Q_c^q} U {U'}} with $\ov{Y_0^0}\cap V_q^c \subset U\subset V_q^c$. 

Moreover the filled-in Julia set $K'_c$ of $f_c$ is contained in $\{\al, \al'\}\cup(Y_0^0\cap V_q^c)$ and 
$K'_c$ is connected if and only if $f_c^n(0) = Q_c^{nq}(0)\in \{\al, \al'\}\cup(Y_0^0\cap V_q^c)$ for all $n$.
\ENDLEM
\PROOF
Apply a small thickening of $\ov{Y_0^0}\cap V_q^c \subset U$ at the ends, se e.g. \cite[Corollary~1.7]{Milnorlc}. 
The proof given there in the case $Q_c^{nq}(0)\in\ov{Y_0^0}$ for all $n$ works for all $c\in\W^\Mbrot_{q-1}(p/q)$.
\ENDPROOF

In the following fix $c\in L_{p/q}$.Then precisely one of the following two cases occur
\ENUMD
\item\label{satelite_case}
For all $n\in\N$ : $Q_c^{nq}(0)\in\ov{Y_0^0}$.
\item\label{post_satelite_case}
There exists $m\geq 1$ minimal such that $Q_c^{mq}(0) \notin \ov{Y_0^0}$. 
\ENDENUM

In the first case \ItemDref{satelite_case} it follows from \lemref{q_renormalization} above that $Q_c$ is $q$-renormalizable.
That is, there exists a quadratic like restriction  {\mapfromto {Q_c^q} U {U'}} with 
connected filled-in Julia set $K'_c \subset Y_0^0\cup\{\al, \al'\}$. 

We shall henceforth focus on the second case \ItemDref{post_satelite_case} 

In order to describe better the second case we set up some additional notation. 
For $0 <k <q$ let $Y_0^k$ denote the level $0$ puzzle piece contained in  $Q_c^k(Y_0^0)$. 
Then $Q_c^k(0)\in Y_0^k$ for any $c\in\W^\Mbrot_{q-2}(p/q)$.
Each $Y_0^k$ is adjacent to $\al$, $Y_0^1 = Y_0^c$ and $Y_0^{q-1} = Y_0(\be')$.
The corresponding twins $\ckY_0^k$ are adjacent to $\al'$ and $\ckY_0^{q-1} = Y_0(\be)$. 
Denote by $X_0$ the interior of $\cup_{k=1}^{q-1}\ov{Y_0^k}$ and by $\ckX_0$ its twin. 
Then the common univalent image $Q_c(X_0)=Q_c(\ckX_0)$ covers all of the level $0$ puzzle 
except $\ov{Y}_0^c$. 

Note that the condition $Q_c^{mq}(0) \notin \ov{Y_0^0}$ in \ItemDref{post_satelite_case} is equivalent to
$Q_c^{mq}(0) \in \ckX_0$ and to $Q_c^{mq}(c)\notin \ov{Y_0^c}$.

When studying parameter space we shall also be interested in the following extension of condition \ItemDref{post_satelite_case} to all of $\W^\Mbrot_0(p/q)$. 
\ENUMDp
\setcounter{enumi}{1}
\item\label{wake_post_satelite_case}
There exist $m\geq 1$ minimal with $Q_c^{mq}(0) \in \ckX_0$.
\ENDENUM
Parameters $c$ satisfying \ItemDpref{wake_post_satelite_case} 
belongs to a dyadic sub-wake of the satellite-copy $\Mbrot_{p/q}$:
\REFPROP{wake_post_satelite_case_in_dyadic_wake}
A parameter $c\in\W^\Mbrot_0(p/q)$ satisfies \ItemDpref{wake_post_satelite_case} if and only if 
$$
c\in \W^\Mbrot_{mq-1}(p/q,r,m) := \W^\Mbrot_{mq-1}(p/q)\cap\W^\Mbrot(p/q,r,m)
$$
where $m\geq1$ is from \ItemDpref{wake_post_satelite_case} 
and $r$ is odd with $0< r < 2^m$.
\ENDPROP
\PROOF
Let $c\in\W^\Mbrot_0(p/q)$ satsify $Q_c^{mq}(0) = Q_c^{mq-1}(c)\in \ckX_0$ for some minimal $m\geq 1$. 
Then $c\in V_{mq-1}^c$ and thus also $c\in \W_{mq-1}^\Mbrot$. 
Moreover $Q_c^{mq}(c)\notin \ov{Y}_0^c$, hence $Q_c^{(m-1)q}(c)$ belongs to the $1/2$ dyadic wake 
$\W_c(p/q,1,1)$ and thus $c\in\W_c(p/q,r,m)$ for some odd $r$ with $0< r < 2^m$ 
by induction and minimality of $m$. 
Hence also $c\in\W^\Mbrot(p/q,r,m)$.
\ENDPROOF

Recall that for $m\geq 1$ and $r$ odd $0<r<2^m$ the de-rooted dyadic decoration is given $L^*(p/q,r,m)$ by 
$$
L^*(p/q,r,m) = L_{p/q}\cap\W^\Mbrot(p/q,r,m).
$$
This gives the following decomposition of the limb $L_{p/q}$, first observed by Douady and Hubbard.
\REFCOR{MLimb_stratification}
The limb $L_{p/q}$ has a natural stratification as
$$
L_{p/q} = \Mbrot_{p/q}\cup\bigcup_{\frac{r}{2^m}} L^*(p/q,r,m).
$$
\ENDCOR
\PROOF
This follows immediately from the dichotomy, \ItemDref{satelite_case}, \ItemDref{post_satelite_case} above.
\ENDPROOF

For any $c\in\W^\Mbrot_0(p/q)$ the common image of puzzle pieces $Q_c(Y_0(\be')) = Q_c(Y_0(\be))$ 
univalently covers $V_0\Sm\ov{X}_0$, and does not intersect the set $\ov{X}_0$. 
It follows that the level $0$ twin puzzle pieces $Y_0(\be')$ and $Y_0(\be) = \ckY_0(\be')$ 
each contain $q$ level $1$ puzzle pieces, 
which are mapped homeomorphically onto $Y_0^0=\ckY_0^0$ and $\ckY_0^k$, $0<k<q$. 
By similar reasoning every other level $0$-puzzle piece contains a unique level $1$ puzzle piece
$Y_1^k\subset Y_0^k$ respectively $\ckY_1^k\subset\ckY_1^k$ which is mapped properly onto 
$Y_0^{k+1}$. 

\REFPROP{Mbeta_c_nest_moves_holomorphically}
For any $c\in\W^\Mbrot_0(p/q)$ the boundaries of all puzzle pieces in the $\be$-nest move holomorphically 
with $c\in\W^\Mbrot_0(p/q)$.
\ENDPROP
\PROOF
The level $0$ Yoccoz graph $\GY_0^c\subset\bd{Y}_0(\be)$ moves holomorphically 
with $c$ over $\W^\Mbrot(p/q)$, since the $c\in\W_c(p/q)\cap V_0^c$, see also \secref{polynomial_wakes}.
Moreover the restriction {\mapfromto {Q_c}{Y_1(\be)}{Y_0(\be)\supset\supset Y_1(\be)}} is biholomorphic 
with a univalent extension to a neighbourhood of $Y_0(\be)$. 
Hence by induction $Y_{n+1}(\be)\subset\subset Y_n(\be) = Q_c(Y_{n+1}(\be))$.  
From this the proposition follows by induction.
\ENDPROOF

By construction there are $q$ puzzle pieces of level $n$ adjacent to $\al_c$ for every $n$. 
And thus $q$ sequences of nested puzzle pieces adjacent to $\al$, $\NN^{\al,k}:= \{Y_n^{\al,k}\}_{n\geq 0}$,
defined by $Y_0^{\al,k}=Y_0^k$. 
Moreover $Q_c$ maps $Y_{n+1}^{\al,k}$ properly onto $Y_n^{\al,(k+1)\!\!\!\mod\! q}$ for every $n$ and $k$ 
and the degree is $1$ unless $k=0$ and $Y_{n+1}^{\al,k}=Y_{n+1}^0$, in which case the degree is $2$. 
It follows immediately that either all $q$ nests are convergent to $\al$ or none is convergent.

In order to create a fundamental system of nested neighbourhoods of $\al_c$ we denote by $Y_n^\al$ 
the interior of $\cup_{k=0}^{q-1}\ov{Y_n^{\al,k}}$, 
so that $Y_n^\al$ is an open neighbourhood of $\al_c$ for all $n$. 
However no $Y_n^\al$ is a puzzle piece. 
Let $r/2^m$, $r$ odd and $0< r < 2^m$ be a dyadic rational and let 
$c\in\ov{\W^\Mbrot(p/q,r,m)}$. 
Then $Q_c$ maps $Y_n^\al$ biholomorphically onto $Y_{n-1}^\al$ and 
$Y_n^\al\subset\subset Y_{n-q}^{\al}$ for all $n\ge (m+1)q$. 
Moreover the boundaries $\partial Y_n^\al$ move holomorphically over $\W^\Mbrot(p/q,r,m)$ 
and continuously over the closure.

\REFLEM{Cantor_little_Julia}
Suppose $c\in\W^\Mbrot_0(p/q)$ satisfies \ItemDpref{wake_post_satelite_case} for some $m\geq 1$.
Let {\mapfromto {f_c = Q_c^q} U {U'}} be a quadratic like map as in \lemref{q_renormalization} 
with $\ov{Y_0^0}\cap V_{mq} \subset U\subset V_{mq}$.
Then 
\ENUM
\item\label{Kprime_Cantor_containment}
the filled Julia set $K'_c\subset \ov{Y_{mq}^{\al,0}}\cup\ov{\ckY_{mq}^{\al,0}}\subset \ov{U}$, 
\item\label{Standard_Cantor_Construct}
the restriction {\mapfromto {f_c} {\ov{Y}_{mq}^{\al,0}}{\ov{Y}^0_{m(q-1)}\subset\ov{\ckX_0}}} 
is a holomorphic diffeomorphism,
\item\label{vanishing_sizes}
$\diam(\ov{Y}_{(n+m)q}) \to 0$ as $n\to\infty$ uniformly over all connected components $Y_{(n+m)q}$ of 
$f_c^{-n}(Y_{mq}^{\al,0}\cup\ckY_{mq}^{\al,0})$. 
\item\label{almost_compactly_contained}
If $Y_{(n+1+m)q}\subset Y_{(n+m)q}$ are nested puzzle pieces with 
$f_c^n(Y_{(n+m)q}), f_c^{n+1}(Y_{(n+1+m)q})\in\{Y_{mq}^{\al,0}, \ckY_{mq}^{\al,0}\}$,  
then $\partial Y_{(n+1+m)q} \cap \partial Y_{(n+m)q}\cap K_c' \subset f_c^{-(n+1)}(\al_c)$.
\item\label{good_neighbourhoods}
In particular for any $z\in K_c'$ either $z$ is not prefixed to $\al_c$ under $f_c$ and 
the nest $\{Y_n(z)\}_n$ is convergent to $z$ or $Q_c^{lq}(z) = \al_c$ for some minimal $l$ and 
there are $q$ nests $\{Y_n^{z,k}\}_n$, $0\le k < q$ convergent to $z$, where 
$Q_c^{lq}(Y_{n+lq}^{z,k}) = Y_n^{\al,k}$ for each $n$ and $k$. 
\ENDENUM
\ENDLEM
\PROOF
The set $U\Sm(\ov{Y_{mq}^{\al,0}}\cup\ov{\ckY_{mq}^{\al,0}})$ consists of points $z$ 
with $f_c^n(z)\in X_0\cup\ckX_0$ for some minimal $n, 0\leq n\leq m$. Thus all such points escapes and so 
$K' \subset (\ov{Y_{mq}^{\al,0}}\cup\ov{Y_{mq}^{\al,0}})$. 
The post-critical orbit $\OO_f$, the forward orbit of $0$ under $f_c$ is finite and disjoint from 
$(\ov{Y_{mq}^{\al,0}}\cup\ov{\ckY_{mq}^{\al,0}})$. 
Hence the latter has finite hypberbolic diameter in $U'\Sm\OO_f$. 
Thus the hyperbolic diameter of $\ov{Y_{(n+m)q}}$ for $Y_{(n+m)q}$ 
any connected component of $f_c^{-n}(Y_{mq}^{\al,0}\cup\ckY_{mq}^{\al,0})$ converges geometrically to $0$, 
as $n\to\infty$. 

By construction $\partial Y^{\al,0}_{mq}\cap\partial Y^0_{(m-1)q} \cap K_c' = \{\al_c\}$ and 
$\partial \ckY^{\al,0}_{mq}\cap\partial Y^0_{(m-1)q}  \cap K_c' = \{\al'_c\}$ so that \Itemref{almost_compactly_contained} 
follows by induction. 

Finally if $z\in K_c'$, then $z\in\End(\NN)$ for a unique nest $\NN = \{Y_l(z)\}_{l\geq 0}$ 
with $f_c^{(n-m)}(Y_{nq})\in\{Y_{mq}^{\al,0},\ckY_{mq}^{\al,0}\}$ for any $n\geq m$. 
And by the above this nest is convergent to $z$. 
Let $\NN^{\al,k}$ denote the $q$ nests adjacent to $\al$. 
Then the nest $\NN^{\al,0}$ is convergent to $\al$ by the first part of the proof and hence all are. 
Thus if $z\in K'_c$ is prefixed to $\al$ by $f_c$, then also all $q$ nests adjacent to $z$ are convergent.
 \ENDPROOF
 
 By \propref{wake_post_satelite_case_in_dyadic_wake} the hypothesis \ItemDpref{wake_post_satelite_case} 
of \lemref{Cantor_little_Julia} is equivalent to $c\in\W^\Mbrot_{mq-1}(p/q,r,m)$ for some odd $r$ with $0<r<2^m$. 
 Fix such $r$ and define for $c\in\W^\Mbrot_{mq-1}(p/q,r,m)$ the set 
 $$
 \Ga'_c := K_c'\cup\bigcup_{n\geq 0} f_c^{-n}(\bd{Y_{mq}^{\al,0}}\cup\bd{\ckY_{mq}^{\al,0}}) 
 = \ov{\bigcup_{n\geq 0} f_c^{-n}(\bd{Y_{mq}^{\al,0}}\cup\bd{\ckY_{mq}^{\al,0}})}.
 $$
 
 \REFPROP{Kprime_nest_move_holomorphically_over_dyadic_wakes}
 Let $m\geq 1$, let $r$ be odd with $0<r<2^m$ and fix $c_*\in\W^\Mbrot_{mq-1}(p/q,r,m)$. 
 Then there exists a holomorphic motion
  $$
 \mapfromto{\psi^{c_*}_{r,m}}{\W^\Mbrot_{mq-1}(p/q,r,m)\times\Ga'_{c_*}}\C
 $$
 with base point $c_*$ such that $\psi^{c_*}_{r,m}(c,\Ga'_{c_*}) = \Ga'_c$ and 
 $f_c\circ\psi^{c_*}_{r,m}(c,z) = \psi^{c_*}_{r,m}(c, f_{c_*}(z))$ for every $c\in\W^\Mbrot_{mq-1}(p/q,r,m)$ 
 and every $z\in\Ga'_{c_*}$. 
 \ENDPROP
 \PROOF
 For $c\in\W^\Mbrot_{mq-1}(p/q,r,m)$ the twin boundaries $\bd{Y_{mq}^{\al,0}}$ and $\bd{\ckY_{mq}^{\al,0}}$ 
 move holomorphically with $c$, because $f_c^m$ maps each boundary 
 onto the boundary of $Y_0^0$ by degree $2^{m-1}$ without passing the critical point, 
 so that $f_c^m$ is a local diffeomorphism around each boundary point. 
 Indeed, if the common forward orbit of the two boundaries were to pass the critical point $0$, 
 then the critical point would end up on the $q$ periodic rays on the boundary of $Y_0^0$, 
 so that $Q_c^{mq}(0)\in\ov{Y}_0^0$, which contradicts that $c\in\W^\Mbrot_{mq-1}(p/q,r,m)$. 
The boundary of $Y_0^0$ moves holomorphically over the larger set $\W^\Mbrot_0(p/q)$, 
which compactly contains $\W^\Mbrot_{mq-1}(p/q,r,m)$. 
And $f_c(z)$ is a holomorphic function of $(c,z)$. 
Thus  $\bd{Y_{mq}^{\al,0}}$ and $\bd{\ckY_{mq}^{\al,0}}$ 
move holomorphically with $c\in\W^\Mbrot_{mq-1}(p/q,r,m)$. 
Secondly for $c\in\W^\Mbrot_{mq-1}(p/q,r,m)$ the map $f_c$ sends 
each puzzle piece $Y_{mq}^{\al,0}$ and $\ckY_{mq}^{\al,0}$ univalently onto the larger (i.e. containing) puzzle piece $Y_{(m-1)q}^0$ and extends as a  diffeomorphism of neighbourhoods of the closures. 
Hence also all pre-images of $\bd{Y_{mq}^{\al,0}}$ and $\bd{\ckY_{mq}^{\al,0}}$ under iterates of $f_c$ 
move holomorphically with $c\in\W^\Mbrot_{mq-1}(p/q,r,m)$. 
Finally $K_c'$ is contained in the closure of the union of all pre-image puzzle boundaries in $\Ga'_c$. 
Thus the proposition follows by the $\la$-lemma for holomorphic motions.
 \ENDPROOF
 
Note that $z\in K_c$ with $Q_c^l(z)=\al_c$ will be adjacent to $2q$ nests if $0$ belongs to the orbit of $z$.

For $0<k<q$ let $X_1^k\subset Y_0(\be')$ 
denote the unique such level $1$ puzzle piece with (univalent) image $\ckY_0^k$ and let 
$X_1$ denote the interior of $\cup_k \ov{X}_1^k$, so that $Q_c$ maps $\ov{X}_1$ diffeomorphically onto ${\ov{\ckX}}_0$.

\REFPROP{basic_trichotomy} 
Let $c\in L_{p/q}$, then for any $z\in K_c$ the orbit falls in precisely one of the following three categories:
\ENUMi
\item\label{z_prebeta}
There exists $l\geq 0$ such that $Q_c^l(z) = \be$.
\item\label{z_capture}
There exists $l\geq 0$ such that $Q_c^l(z)\in K_c'$.
\item\label{z_level_zero_recurrent}
There exists a strictly increasing sequence $\{l_n\}_{n\geq 0}$ 
with $Q_c^{l_n}(z) \in X_1$ for all $n$.
\ENDENUM
\ENDPROP
\PROOF
Notice at first that for any point $z\in K_c$ which does not satisfy \Itemiref{z_prebeta} there exists $l\geq 0$, 
such that $Q_c^l(z)\in \ov{Y}_0(\be')$. 
If $Q_c^l(z)\notin X_1$, then $Q_c^{l+1}(z)\in \ov{Y_0^0}$, and thus $Q_c^{l+q}(z)\in \ov{Y}_0(\be')$. 
Hence either $Q_c^{l+1+nq}(z)\in \ov{Y_0^0}$ for all $n\geq 0$ so that the orbit of $z$ satisfies \Itemiref{z_capture}, 
by \lemref{q_renormalization}, or there exists some $n$ such that $Q_c^{l+nq}(z)\in X_1$, set $l_0=l+nq$. 
Then apply the same argument recursively to first $Q_c^{l_0}$, noting that $Q_c(Q_c^{l_0}(z))\notin Y_0(\be')\supset X_1$, 
to obtain the desired strictly increasing sequence with $Q_c^{l_n}(z)\in X_1$. 
\ENDPROOF

\REFPROP{simpel_convergnece_of_z_Nests}
Let $c\in L_{p/q}$ satisfy \ItemDref{post_satelite_case}
In the first two cases \Itemiref{z_prebeta} and \Itemiref{z_capture} of \propref{basic_trichotomy} 
any nest $\{Y_n\}_n$ such that $z\in\ov{Y}_n$ for every $n$ is convergent to $z$. 
Moreover if $z=c$ then there exists $N\geq l$ such that the restriction {\mapfromto {Q_c^l} {Y_N} {Q_c^l(Y_N)}} is univalent.
\ENDPROP
Recall that there is a unique $Y_n = Y_n(z)$ with $z\in\ov{Y_n}$ except if $Q_c^k(z) = \al_c$ for some $k$, 
in which case there are precisely $q$ such nests if the orbit of $z$ avoids the critical point $0$ 
and $2q$ such nests if not.
\PROOF
If there exists $l\geq 0$ such that $Q_c^l(z) = \be$, take $l$ minimal with this property. 
Then the restrictions {\mapfromto {Q_c^l} {Y_{n+l}}{Y_n(\beta)}} are proper maps of non-increasing degrees $d_n$ for every $n\geq 0$. 
Since the nest $\{Y_n(\beta)\}_n$ is convergent to $\beta$ the nest $\{Y_n\}_n$ is convergent to $z$. 
If $0 = Q_c^r(z)$ for some $r\leq l$ then $d_n=2$ for every $n\geq 0$. 
Otherwise, since the nest $\{Y_n\}_n$ is convergent to $z$, 
there exists $N\geq l$ so that $0\notin Q_c^k(Y_N)$ for any $k<l$ and so 
the restriction {\mapfromto {f^l} {Y_{N+l}} {Y_N(\beta)}} is univalent.
In particular if $z=c$, then this restriction is univalent.

For the remaining cases let $m\geq 1$ be given by $c$ satisfying \ItemDref{post_satelite_case}. 

If there exists $l\geq 0$ such that $Q_c^l(z) := w \in K_c'$, with $l$ minimal and $z$ is not prefixed to $\al$. 
Then the restrictions {\mapfromto {Q_c^l} {Y_{n+l}}{Y_n(w)}} are proper maps of non-increasing degrees $d_n$ for every $n\geq 0$. 
Since the nest $\{Y_n(w)\}_n$ is onvergent to $w$ by \lemref{Cantor_little_Julia}, the nest $\{Y_n\}_n$ is convergent to $z$. 
If $0 = Q_c^r(z)$ for some $r\leq l$ then $d_n=2$ for every $n\geq mq$. 
Otherwise, since the nest $\{Y_n\}_n$ is convergent to $z$, there exists $N\geq mq$ so that $0\notin Q_c^r(Y_N)$ for any $r<l$ 
and thus $d_n=1$ for $n\geq N$. 
In particular if $z=c$ then the restriction {\mapfromto {Q_c^l} {Y_{N+l}}{Y_N(w)}}is univalent.

Finally if there exists $l\geq 0$ such that $Q_c^l(z) = \al \in K_c'$, with $l$ minimal. 
Then there exists $k, 0\leq k <q$ such that  the restrictions {\mapfromto {Q_c^l} {Y_{n+l}}{Y^{\al,k}_n}} are proper maps of non-increasing degrees. 
Since the nest $\{Y_n^{\al, k}\}_n$ is convergent to $\al$ by \lemref{Cantor_little_Julia}, the nest $\{Y_n\}_n$ is convergent to $z$. And hence also any other nest adjacent to $z$ is convergent.
\ENDPROOF

The above Proposition immediately gives the following Corollary for parameterspace.
\REFCOR{sub_post_satelite_case}
Let $c\in L_{p/q}$ and suppose there exists $m\geq 1$ minimal such that $Q_c^{mq}(0) \notin \ov{Y_0^0}$. 
Then precisely one of the following three cases occur
\ENUMi
\item\label{prebeta}
There exists $l\geq mq$ such that $Q_c^l(0) = \beta$.
\item\label{capture}
There exists $l> mq$ such that $Q_c^l(0) \in K'$.
\item\label{critical_level_zero_recurrent}
There exists a strictly increasing sequence $\{l_n\}_{n\geq 0}$ with $l_0 = mq-1$ 
with $Q_c^{l_n}(0) \in X_1$ for all $n$.
\ENDENUM
Moreover in both cases \Itemiref{prebeta} and \Itemiref{capture} any nest $\{Y_n\}_n$ such that $c\in\ov{Y}_n$ 
for every $n$ is convergent to $c$.
\ENDCOR 
In the last statement there is a unique such $Y_n = Y_n(c)$ except if $Q_c^k(c) = \al_c$ for some $k$, 
in which case there are precisely $q$ such nests, which all are convergent to $\al_c$ by the discussion above.

Note that $Y_0^{q-1}\Sm\ov{X}_1$ is a non degenerate annulus contained in any of the annuli 
$Y_0^{q-1}\Sm\ov{X}_1^k$, $0<k<q$.

\REFTHM{p:yocc} 
Let $c\in L_{p/q}$ satisfy the hypotheses of \corref{sub_post_satelite_case} and the property 
\Itemiref{critical_level_zero_recurrent} therein. 
Then there exists a non degenerate annulus $A_{n_0} = Y_{n_0}\Sm\overline{Y}_{n_0+1}$ 
between nested puzzle pieces of the Yoccoz puzzle for $Q_c$ with $Q_c^{n_0}(Y_{n_0}) = Y_0^{q-1}$, 
$Q_c^{n_0}(Y_{n_0+1}) = X_1^k$ for some $0<k<q$. 
And there exists a nested sequence of annuli $A_{n_i}^c=Y_{n_i}^c\setminus
 \overline{ Y_{n_i+1}^c}$, $i>0$ with $n_0 < n_1\nearrow\infty$ surrounding the critical value $c$ 
 such that\,:
 \begin{itemize}
  \item the map $Q_c^{n_i-n_0} : A_{n_i}^c \to A_{n_0}$ is a covering map of degree $2^{d_i}$, $d_i\ge 0$ for $i\ge 1$\,;
  \item in particular also all the annuli $A^c_{n_i}$, $i>0$ are non degenerate\,;
 \item
 \begin{itemize}
\item  either the sum $\displaystyle \sum _{i\ge 1}\Mod(A_{n_i}^c) = \Mod(A_{n_0})\sum _{i\ge 1} \frac{1}{2^{d_i}}$ is infinite 
and\\
the intersection $\displaystyle\bigcap_{n\ge 0} \overline{Y_n^c}$ reduces to a point, 
\item or there exists $k>0$ such that for all $n$ large enough the map\\ 
$Q_c^k:Y_{n+k}^c\to Y_n^c$ is quadratic-like with connected filled-in Julia set. 
\end{itemize} 
 \end{itemize}
\ENDTHM

\noindent Note that the annulus $A_{n_0}$ does not necessarily surround the critical value. 
 
\thmref{p:yocc} follows from a classical tableaux argument, see e.g..~\cite{Milnorlc}. 
The degree $2^{d_0}$ of the restriction {\mapfromto {Q_c^{n_0}} {Y_{n_0}} {Y_0^{q-1}}} 
may be larger than the degree of the restriction {\mapfromto {Q_c^{n_0}} {Y_{n_0+1}} {X_1^k}}. 
This happens precisely when $Y_0^{q-1}\Sm\ov{X}_1^k$ contains one or more critical values for 
the restriction of $Q_c^{n_0}$ to $A_{n_0}$
However as the long composition of $Q_c$ with itself has degree either $1$ or $2$ in each step, 
it easily follows that 
$$
\Mod(A_{n_0}) \ge \Mod(Y_0^{q-1}\Sm\ov{X}_1^k)/2^{d_0} \ge \Mod(Y_0^{q-1}\Sm\ov{X}_1)/2^{d_0}.
$$

Note that \thmref{p:yocc} can also be proved  using the following results\,: 
\cite[Lemma~1.22]{R} for the non-recurrent case and \cite[Lemma~1.25 case~1)]{R} for the recurrent case.     

\subsection{Yoccoz-type estimates for the parabolic maps $g_B$}
In this section we port the results above for the quadratic polynomials $Q_c$ 
to the parabolic quadratic rational maps $g = g_B = z + 1/z + B$, $\Re(B) >0$ 
with $\be_B = \infty$ a parabolic fixed point of multiplier $1$ and 
a unique finite fixed point $\al_B= -1/B$ of multiplier $A(B) = 1 -B^2$. 
We let $\tau(z) = 1/z$ denote the covering involution for $g$.

As for quadratic polynomials we denote by $\be' = \be'_B = 0$ the finite preimage of $\be$. 
Similarly we denote by $\al' = \al'_B = -B$ the non fixed pre- image of $\al_B$.
The critical point $1$ for $g_B$ is first attracted in the sense that the extended 
attracting Fatou coordinate $\phi_B$ for $g_B$ maps $1$ to $0$ 
and the domain $\Om_B$ with $1\in\bd{\Om}_B$ univalently onto $\Hplus$. 
The other or second critical point $-1$ for $g_B$ and its critical value $v_B =g_B(-1) = -2+B$ 
play the same role for $g_B$ as the critical point $0$ and its critical value $c$ plays for $Q_c$. 
In particular the second critical point and value belong to the filled-in Julia set $K_B$, 
if and only if $K_B$ is connected and is otherwise in $\La_g\Sm g_B(\Om_B)$. 

In the rest of this subsection we shall fix an irreducible rational $p/q$ and consider $B\in\W^\Mone(p/q)$. 
We setup notation for special puzzle pieces for $g_B$ corresponding to the notation for special puzzle pieces for $Q_c$.

Recall that $V_n^\PP$ is the interior of the union of closures of level $n$ universal parabolic puzzle pieces 
for $n\geq 0$. And that 
$$\W_n^\Mone(p/q) := \LL^\Mone_{p/q}\cup \{ B \in \W^\Mone(p/q)\;|\; h_B(v_B) \in V_n^\PP\}$$ 
for $n\geq 0$.

Let $B\in\W^\Mone_0(p/q)$. 
We shall use the abbreviations $P_n^0 := P_n(-1)\in\PP_n$ for the critical puzzle piece of depth $n$
and $P_n^B:=P_n(v_B) = g_B(P_{n+1}^0)\in\PP_n$ for the critical value puzzle piece of depth $n$, 
whenever there is such a puzzle piece. 
We shall use the symbol $\ckP= \tau(P)$ for the dynamical twin of the puzzle piece $P$, i.e.~$g_B(\ckP) = g_B(P)$. 

\REFLEM{par_q_renormalization}
Let $B\in\W^\Mone_{q-1}(p/q)$. 
Then the map $g_B^q$ has a quadratic-like restriction 
{\mapfromto {f=f_B = g_B^q} U {U'}} with $\ov{P_0^0}\cap V_q^B \subset U\subset V_q^B$. 

Moreover the filled-in Julia set $K_B'$ of $f$ is contained in $\{\al, \al'\}\cup(P_0^0\cap V_q^B)$ and 
$K_B'$ is connected if and only if $f^n(-1) = g_B^{nq}(-1)\in \{\al, \al'\}\cup(P_0^0\cap V_q^B)$ for all $n$.
\ENDLEM
\PROOF
See the proof of the similar \lemref{q_renormalization} above for the corresponding polynomials $Q_c$.
\ENDPROOF

In the following fix $B\in \LL^\Mone_{p/q}$. Then just as for quadratic polynomials precisely one of the following two cases occur
\ENUMDB
\item\label{par_satelite_case}
For all $n\in\N$ : $g_B^{nq}(-1)\in\ov{P_0^0}$.
\item\label{par_post_satelite_case}
There exists $m\geq 1$ minimal such that $g_B^{mq}(-1) \notin \ov{P_0^0}$.
\ENDENUM

In the first case \ItemDBref{satelite_case} it follows from \lemref{par_q_renormalization} above that $g_B$ is $q$-renormalizable.
That is, there exists a quadratic like restriction {\mapfromto {f_B=g_B^q} U {U'}} with 
connected filled-in Julia set $K_B' \subset P_0^0\cup\{\al, \al'\}$. 

As for polynomials we shall henceforth focus on the second case \ItemDBref{post_satelite_case} 

We continue to set up notation analogous to the polynomial case. 
For $0\leq k <q$ let $P_0^k$ denote the level $0$ puzzle piece contained in $g_B^k(P_0^0)$. 
Then $g_B^k(-1) \in P_0^k$ for any $B\in\W_{q-2}^\Mone(p/q)$.
Each $P_0^k$ is adjacent to $\al$ and $P_0^{q-1} = P_0(\be')$.
The corresponding twins $\ckP_0^k$ are adjacent to $\al'$ and $\ckP_0^{q-1} = P_0(\be)$. 
Denote by $S_0$ the interior of $\cup_{k=1}^{q-1}\ov{P_0^k}$ and by $\ckS_0$ its twin, 
that is $S_0$ plays the role of $X_0$. 
Then the common image $g_B(S_0) = g_B(\ckS_0)$ covers the level $0$ puzzle except for $\ov{P}_0^B$ 
and the subset $D_0^B$ between the shortcut $\ckga_0^B=\bd{P}_0(\be)\cap D_0^B$ and $g_B(\ckga_0^B)\subset D_0^B$.

Note that the condition $g_B^{mq}(-1) \notin \ov{P_0^0}$ in \ItemDBref{par_post_satelite_case} is equivalent to
$g_B^{mq}(-1) \in \ckS_0$. 

When studying parameter space we shall as for polynomials also be interested in the following extension of condition \ItemDBref{par_post_satelite_case} on $B\in\W^\Mone_0(p/q)$. 

\ENUMDBp
\setcounter{enumi}{1}
\item\label{par_wake_post_satelite_case}
There exist $m\geq 1$ minimal with $g_B^{mq}(-1) \in \ckS_0$.
\ENDENUM
Parameters $B$ satisfying \ItemDBpref{par_wake_post_satelite_case} 
(see \defref{para_p_q_dyadic_wakes} for the definition) belongs to a dyadic sub-wake of 
the satellite-copy $\Mbrot^\Mone_{p/q}$:
\REFPROP{par_wake_post_satelite_case_in_dyadic_wake}
A parameter $B\in\W^\Mone_0(p/q)$ satisfies \ItemDBpref{par_wake_post_satelite_case} if and only if 
$$
B\in \W^\Mone_{mq-1}(p/q,r,m) := \W^\Mone_{mq-1}(p/q)\cap\W^\Mone(p/q,r,m)
$$
where $m\geq1$ is from \ItemDBpref{par_wake_post_satelite_case} and $r$ is odd with $0< r < 2^m$.
\ENDPROP
\PROOF
Let $B\in\W^\Mone_0(p/q)$ satsify $g_B^{mq}(-1) = g_B^{mq-1}(v_B)\in \ckS_0$ for some minimal $m\geq 1$. 
Then $v_B\in V_{mq-1}^B$ and thus also $B\in \W_{mq-1}^\Mone$. 
Moreover $g_B^{mq}(v_B)\notin \ov{P}_0^B$, hence $g_B^{(m-1)q}(v_B)$ belongs to the $1/2$ dyadic wake 
$\W_B(p/q,1,1)$ and thus $B\in\W_B(p/q,r,m)$ for some odd $r$ with $0< r < 2^m$ 
by induction and minimality of $m$. 
Hence also $B\in\W^\Mone(p/q,r,m)$.
\ENDPROOF

Recall that for $m\geq 1$ and $r$ odd $0<r<2^m$ the derooted dyadic decoration $\LL^\Mone_*(p/q,r,m)$ 
is the set of parameters 
$$
\LL^\Mone_*(p/q,r,m) = \LL^\Mone_{p/q}\cap\W^\Mone(p/q,r,m).
$$
Let $\LL^\Mone(p/q,r,m)$ denote the limb with root, i.e. $\LL^\Mone_*(p/q,r,m)$ union 
the root point of $\W^\Mone_*(p/q,r,m)$ (see also \defref{para_p_q_dyadic_wakes} and trailing comments).
This gives the following decomposition of the limb $\LL^\Mone_{p/q}$, corresponding to the decomposition of 
limbs of the Mandelbrot set.
\REFCOR{MoneLimb_stratification}
The limb $\LL^\Mone_{p/q}$ has a natural stratification as
$$
\LL^\Mone_{p/q} = \Mbrot_{p/q}\cup\bigcup_{\frac{r}{2^m}} \LL^\Mone_*(p/q,r,m).
$$
\ENDCOR
\PROOF
This follows immediately from the dichotomy, \ItemDBref{par_satelite_case}, \ItemDBref{par_post_satelite_case} above.
\ENDPROOF

As an immediated corollary of \propref{base_para_holomorphic_motion} we obtain 
\REFPROP{Monebeta_c_nest_moves_holomorphically}
The boundaries of the puzzle pieces in both the $\be_B$-nest 
$\NN(\be_B) = \{P_n(\be_B)\}_{n\geq 0}$ 
and the $\be'_B$-nest $\NN(\be'_B) = \{P_n(\be'_B)\}_{n\geq 0}$ move holomorphically with 
$B\in\W^\Mone_0(p/q)$.
\ENDPROP
It was proven in \cite{PR2} that the $\be$ and $\be'$-nest are convergent to $\be$ and $\be'$ respectively.

Similarly to the polynomial case, for any $B\in\W^\Mone_0(p/q)$ 
the level $0$ twin puzzle pieces $P_0(\be')$ and $P_0(\be) = \ckP_0(\be')$ 
each contain $q$ level $1$ puzzle pieces, 
which are mapped homeomorphically onto $P_0^0=\ckP_0^0$ and $\ckP_0^k$, $0<k<q$ 
except for the slight variation that due to the short cuts $\ov{P}_1^{q-1} \subset P_0^{q-1}$ and 
similarly for the twins, but $g_B(P_1^{q-1}) = g_B(\ckP_1^{q-1})\cap D_0^B$ differs in an 
inessential way from $\ckP_0^{q-1}\cap D_0^B$.
And every other level $0$-puzzle piece contains a unique level $1$ puzzle piece
$P_1^k\subset P_0^k$ respectively $\ckP_1^k\subset\ckP_1^k$ which is mapped properly onto 
$P_0^{k+1}$. 

By construction there are $q$ puzzle pieces of level $n$ adjacent to $\al_B$ for every $n$. 
And thus $q$ sequences of nested puzzle pieces $\NN^{\al,k}:= \{P_n^{\al,k}\}_{n\geq 0}$, adjacent to $\al_B$ 
and defined by $P_0^{\al,k}=P_0^k$. 
Moreover $g_B$ maps $P_{n+1}^{\al,k}$ properly onto $P_n^{\al,(k+1)\!\!\!\mod\! q}$ for every $n$ and $k$ 
and the degree is $1$ unless $k=0$ so that $P_{n+1}^{\al,k}=P_{n+1}^0$. 
It follows immediately that either all $q$ nests are convergent to $\al$ or none is convergent.

\REFLEM{par_Cantor_little_Julia}
Let $m\geq 1$ and suppose that $B\in\W^\Mone_{mq-1}(p/q)$ satisfies \ItemDBref{par_post_satelite_case}~with this $m$.
Let {\mapfromto f U {U'}} be a quadratic like map as in \lemref{par_q_renormalization} with $k=mq$, 
then 
\ENUM
\item\label{par_Kprime_Cantor_containment}
the filled Julia set $K'_B\subset \ov{P_{mq}^{\al,0}}\cup\ov{\ckP_{mq}^{\al,0}}\subset \ov{U}$, 
\item\label{par_Standard_Cantor_Construct}
the restriction {\mapfromto {f_B} {\ov{P}_{mq}^{\al,0}}{\ov{P}^0_{m(q-1)}\subset\ov{\ckS_0}}} 
is a holomorphic diffeomorphism,
\item\label{par_vanishing_sizes}
$\diam(\ov{P}_{(n+m)q}) \to 0$ as $n\to\infty$ uniformly over all connected components $P_{(n+m)q}$ of 
$f_B^{-n}(P_{mq}^{\al,0}\cup\ckP_{mq}^{\al,0})$. 
\item\label{par_almost_compactly_contained}
If $P_{(n+1+m)q}\subset P_{(n+m)q}$ are nested puzzle pieces with 
$f_B^n(P_{(n+m)q}), f_B^{n+1}(P_{(n+1+m)q})\in\{P_{mq}^{\al,0}, \ckP_{mq}^{\al,0}\}$, 
then $\partial P_{(n+1+m)q} \cap \partial P_{(n+m)q} \cap K_B' \subset f_B^{-(n+1)}(\al_B)$.
\item\label{par_good_neighbourhoods}
In particular for any $z\in K_B'$ either $z$ is not prefixed to $\al_B$ under $g_B^q$ and 
the nest $\{P_n(z)\}_n$ is convergent to $z$ or $g_B^{lq}(z) = \al_B$ for some minimal $l$ and 
there are $q$ nests $\{P_n^{z,k}\}_n$, $0\le k < q$ convergent to $z$ and with 
$g_B^{lq}(P_{n+lq}^{z,k}) = P_n^{\al,k}$ for each $n$ and $k$. 
\ENDENUM
\ENDLEM

In order to create a fundamental system of nested neighbourhoods of $\al_B$ we denote by $P_n^\al$ 
the interior of $\cup_{k=0}^{q-1}\ov{P_n^{\al,k}}$, 
so that $P_n^\al$ is an open neighbourhood of $\al_B$ for all $n$. 
However no $P_n^\al$ is a puzzle piece. 
Let $r/2^m$, $r$ odd and $0< r < 2^m$ be a dyadic rational and let 
$B\in\ov{\W^\Mone(p/q,r,m)}$. 
Then for all $n\ge (m+1)q$ $g_B$ maps $P_n^\al$ biholomorphically onto $P_{n-1}^\al$ and 
$P_n^\al\subset\subset P_{n-q}^{al}$. 
Moreover the the union $\cup\partial P_n^\al$ of boundaries of $\partial P_n^\al$ move holomorphically over 
$\W^\Mone(p/q,r,m)$ and continuously over the closure.

Recall that $f$ similarly maps $\ov{P_{mq}^{\al,0}}$ diffeomorphically onto $\ov{P_{m(q-1)}}$.
\PROOF
The proof is completely analogous to the proof of \lemref{Cantor_little_Julia} 
above for polynomials and is left to the reader.
\ENDPROOF

By \propref{par_wake_post_satelite_case_in_dyadic_wake} the hypothesis 
\ItemDBpref{par_wake_post_satelite_case}~of \lemref{par_Cantor_little_Julia} is equivalent to $B\in\W^\Mone_{mq-1}(p/q,r,m)$ for some odd $r$ with $0<r<2^m$. 
 Fix such $r$ and define for $c\in\W^\Mone_{mq-1}(p/q,r,m)$ the set 
 $$
 \Ga'_B := K_B'\cup\bigcup_{n\geq 0} f_B^{-n}(\bd{P_{mq}^{\al,0}}\cup\bd{\ckP_{mq}^{\al,0}}) 
 = \ov{\bigcup_{n\geq 0} f_B^{-n}(\bd{P_{mq}^{\al,0}}\cup\bd{\ckP_{mq}^{\al,0}})}.
 $$
 
 \REFPROP{par_Kprime_nest_move_holomorphically_over_dyadic_wakes}
 Let $m\geq 1$, let $r$ be odd with $0<r<2^m$ and fix $B_*\in\W^\Mone_{mq-1}(p/q,r,m)$. 
 Then there exists a holomorphic motion
  $$
 \mapfromto{\psi^{B_*}_{r,m}}{\W^\Mone_{mq-1}(p/q,r,m)\times\Ga'_{B_*}}\C
 $$
 with base point $B_*$ such that $\psi^{B_*}_{r,m}(B,\Ga'_{B_*}) = \Ga'_B$ and 
 $f_B\circ\psi^{B_*}_{r,m}(B,z) = \psi^{B_*}_{r,m}(B, f_{B_*}(z))$ for every $B\in\W^\Mone_{mq-1}(p/q,r,m)$ 
 and every $z\in\Ga'_{B_*}$. 
 \ENDPROP
 \PROOF
 The proof is completely analogous to the proof of 
 \propref{Kprime_nest_move_holomorphically_over_dyadic_wakes} 
above for polynomials and is left to the reader.
 \ENDPROOF
 
Note that $z\in K_B$ with $g_B^l(z)=\al_B$ will be adjacent to $2q$ nests if $-1$ belongs to the orbit of $z$.

For $0<k<q$ let $S_1^k\subset P_0(\be')$ and $\ckS_1^k\subset \ckP_0(\be') = P_0(\be)$ 
corresponding to $X_1^k$ and $\ckX_1^k$ denote the level $1$ such puzzle pieces 
different from $P_1^{\al,q-1}$ and $\ckP_1^{\al,q-1}$ respectively 
indexed so that $g_B(S_1^k) = \ckP_0^k$ for $0< k < q-1$ and $S_1^{q-1}=P_1(\be')$. 
As with $X_1$ in the polynomial case let $S_1$ denote the interior of $\cup_k \ov{S_1^k}$.

\REFPROP{par_basic_trichotomy} 
Let $B\in \LL_{p/q}^\Mone$ and assume $g_B^{(m-1)q}(-1) \in \ov{P^0}$ and $g_B^{mq}(-1) \notin \ov{P^0}$ for some $m\geq 1$.
Then for any $z\in K_B$ the orbit falls in precisely one of the following three categories:
\ENUMi
\item\label{par_z_prebeta}
There exists $l\geq 0$ such that $g_B^l(z) = \be$.
\item\label{par_z_capture}
There exists $l\geq 0$ such that $g_B^l(z)\in K'$.
\item\label{par_z_level_zero_recurrent}
There exists a strictly increasing sequence $\{l_n\}_{n\geq 0}$ 
with $g_B^{l_n}(z) \in S_1$ for all $n$.
\ENDENUM
Moreover in both cases \Itemiref{par_z_prebeta} and \Itemiref{par_z_capture} any nest $\{P_n\}_n$ such that $z\in\ov{P}_n$ 
for every $n$ is convergent to $z$.
\ENDPROP
In the last statement there is a unique such $P_n = P_n(z)$ except if $g_B^k(z) = \al_B$ for some $k$, 
in which case there are precisely $q$ such nests. By the above these $q$ nests are either all convergent 
or all divergent.

\PROOF
Again the proof is analogous to the proof of \propref{basic_trichotomy} and is left to the reader.
\ENDPROOF

\REFPROP{par_simpel_convergnece_of_z_Nests}
Let $B\in \LL^\Mone_{p/q}$ satisfy \ItemDBref{par_post_satelite_case}
In the first two cases \Itemiref{par_z_prebeta} and \Itemiref{par_z_capture} of \propref{par_basic_trichotomy} 
any nest $\{P_n\}_n$ such that $z\in\ov{P}_n$ for every $n$ is convergent to $z$. 
Moreover if $z=v_B$ then there exists $N\geq l$ such that the restriction 
{\mapfromto {g_B^{l-1}} {P_N} {P_{N-l+1}(\be')}} is univalent in .\Itemiref{par_z_prebeta} 
and {\mapfromto {g_B^l} {P_N} {g_B^l(P_N)}} is univalent in \Itemiref{par_z_capture}.
\ENDPROP
Recall that there is a unique $P_n = P_n(z)$ with $z\in\ov{P_n}$ except if $g_B^k(z) = \al_B$ for some $k$, 
in which case there are precisely $q$ such nests if the orbit of $z$ avoids the critical point $-1$ 
and $2q$ such nests if not.
\PROOF
The proof is mostly analogous to the proof of \propref{simpel_convergnece_of_z_Nests} we point out the difference and leave the rest to the reader. 
The difference is due to the fact that the relation between $P_{n+1}(\be')$ 
and $P_n(\be)$ is only partly dynamical because of the short-cuts on the boundary. 
However the only way to $\be$ from $v_B$ is via $\be'$ and all pre-images under iteration of any 
of the puzzle pieces $P_n(\be')$ are dynamical, 
hence the $l-1$ in place of $l$ in the formula above in the case \Itemiref{par_z_prebeta}.
Other than this the proof is completely analogous to the proof of \propref{simpel_convergnece_of_z_Nests}.
\ENDPROOF

As in the polynomial case the above Proposition immediately gives the following Corollary for parameterspace.

\REFCOR{par_sub_post_satelite_case}
Let $B$ belong to a dyadic decoration $\LL^\Mone(p/q,r,m)$ for some $r$ is odd and $0< r < 2^m$. 
Then precisely one of the following three cases occur
\ENUMi
\item\label{par_prebeta}
There exists $l\geq mq$ such that $g_B^l(-1) = \beta_B$.
\item\label{par_capture}
There exists $l> mq$ such that $g_B^l(-1) \in K_B'$.
\item\label{par_level_zero_recurrent}
There exists a strictly increasing sequence $\{l_n\}_{n\geq 0}$ with $l_0 = mq-1$ 
with $g_B^{l_n}(-1) \in S_1$ for all $n$.
\ENDENUM
Moreover in both cases \Itemiref{par_prebeta} and \Itemiref{par_capture} any nest $\{P_n\}_n$ such that 
the critical value $v_B\in\ov{P}_n$ 
for every $n$ is convergent to $v_B$.
\ENDCOR 
In the last statement there is a unique such $P_n = P_n(v_B)$ except if $g_B^k(v_B) = \al_B$ for some $k$, 
in which case there are precisely $q$ such nests all of which are convergent by the discussion above.

Note that $P_0^{q-1}\Sm\ov{S}_1$ is a non degenerate annulus contained in any of the annuli 
$P_0^{q-1}\Sm\ov{S}_1^k$, $0<k<q$.

For the rest of this paragraph we fix $B\in\LL_{p/q}^\Mone$ 
and we let $c\in L_{p/q}$ be a parameter with $\siminfB = \siminfc$ 
as provided by \lemref{samebehaviour}. 

\REFPROP{dyadics_correspond}
Let $B\in\LL^\Mone_{p/q}$ and let $c\in L_{p/q}$ be a parameter with $\siminfB = \siminfc$. 
If $B\in\Mbrot^\Mone_{p/q}$ then $c\in\Mbrot_{p/q}$. 
And if $B\in\LL^\Mone(p/q,r,m)$ for some $r$ is odd and $0< r < 2^m$, 
then $c\in L(p/q,r,m)$ and moreover
\ENUMi
\item\label{both_prebeta}
$g_B^l(-1) = \beta_B$ if and only if  $Q_c^l(0) = \beta_c$.
\item\label{both_capture}
$g_B^l(-1) \in K_B'$ if and only if $Q_c^l(0) \in K'_c$
\item\label{both_level_zero_recurrent}
$g_B^l(-1) \in S_1$ if and only if $Q_c^l(0) \in X_1$.
\ENDENUM
\ENDPROP
\PROOF
Let {\mapfromto {\ckg_B} {\PP} {\PP}} denote the map of puzzle pieces induced by $g_B$ (defined in \defref{abstract_g}). 
And let {\mapfromto \chi {\YY} {\PP}} denote the dynamical correspondence between puzzles of \propref{samebehaviour}. 
Then nests are mapped to nests and in particular the critical value nest $\{Y_n^c\}_n$ is mapped 
to the critical value nest $\{P_n^B\}_n$. 
From this it follows that $c\in L(p/q,r,m)$ and \Itemiref{both_level_zero_recurrent} follows. 
Combining further with the descriptions of $K_c'$ in \lemref{q_renormalization} and $K_B'$ in \lemref{par_q_renormalization} yields \Itemiref{both_capture}. 
Finally the $\beta$-nest $\{Y_n(\beta_c)\}_n$ is easily seen to always be convergent. 
And the $\beta$-nest $\{P_n(\beta_B)\}_n$ was proven to always be convergent in \cite[Prop 5.10]{PR2}. 
So that also \Itemiref{both_prebeta} also follows from $\chi$ conjugating puzzle dynamics
\ENDPROOF

We are now ready to state and prove a parabolic analog of \thmref{p:yocc}: 

\REFTHM{p:yoccparabo} 
Let $B\in\LL^\Mone_{p/q}$ satisfy the hypotheses of \corref{par_sub_post_satelite_case} and its property 
\Itemiref{par_level_zero_recurrent}. 
Then there exists a non degenerate annulus $A_{n_0}^B = P_{n_0}\Sm\overline{P}_{n_0+1}$ 
between nested puzzle pieces of the parabolic Yoccoz puzzle for $g_B$ with $g_B^{n_0}(P_{n_0}) = P_0^{q-1}$, 
$g_B^{n_0}(P_{n_0+1}) = S_1^k$ for some $0<k<q$. 

And there exists a nested sequence of annuli 
$A_{n_i}^B=P_{n_i}^B\setminus\overline{ P_{n_i+1}^B}$, $i>0$ 
with $n_0 < n_1\nearrow\infty$ surrounding the critical value $v_B$ such that\,:
\begin{itemize}
\item the map $g_B^{n_i-n_0} : A_{n_i}^B \to A_{n_0}^B$ 
is a covering map of degree $2^{d_i}$, $d_i\ge 0$ for $i\ge 1$\,;
 \item 
 in particular also all the annuli $A^B_{n_i}$, $i>0$ are non degenerate\,;
 \item
\begin{itemize}
\item 
either the sum $\displaystyle \sum _{i\ge 1}\Mod(A_{n_i}^B) = \Mod(A_{n_0}^B)\sum _{i\ge 1} \frac{1}{2^{d_i}}$ is infinite,\\
$P_n^B$ is defined for all $n$ and 
$$\End(\{P_n^B\}_n) = v_B,$$
\item or there exists $k>0$ such that for all $n$ large enough the map\\ 
$$g_B^k:P_{n+k}^B\to P_n^B$$ 
is quadratic-like with connected filled-in Julia set. 
\end{itemize} 
 \end{itemize}
\ENDTHM
Note that the above sums are finite only when $\sim$ is renormalizable of period $k>q$. 
As with \thmref{p:yocc} the annulus $A^B_{n_0}$ will in general not surround the critical value $v_B$. 
\PROOF
There are two immediate proof strategies. 
Either redo the usual puzzle argument or as we shall do here combine  
\propref{samebehaviour} with \thmref{p:yocc}. 
And let $c\in L_{p/q}$ be a parameter with $\siminfB = \siminfc$. 
Let $\{Y_{n_i}\}_i$ and $\{Y_{n_i+1}\}_i$ be the sequences of Yoccoz puzzle pieces given by \thmref{p:yocc}. 
And for each $i\ge 0$ let $P_{n_i} = \chi(Y_{n_i})$, $P_{n_i+1} = \chi(Y_{n_i+1})$. 
Then by \propref{samebehaviour} the desired properties for the non-degenerate annuli $A_{n_i}^B$ 
follows from the similar properties of the annuli $A_{n_i}^c$ in \thmref{p:yocc}. 
Moreover the covering degree $d_i$ are the same so that 
$$
\sum _{i\ge 1}\Mod(A_{n_i}^B) = \Mod(A_{n_0}^B)\sum _{i\ge 1} \frac{1}{2^{d_i}} = 
\frac{\Mod(A_{n_0}^B)}{\Mod(A_{n_0}^c)} \sum _{i\ge 1}\Mod(A_{n_i}^c).
$$
\ENDPROOF
\REFCOR{similar_critical_value_ends}
Let $B\in\LL^\Mone_{p/q}$ satisfy the hypotheses of \corref{par_sub_post_satelite_case} 
and let $c\in L_{p/q}$ be a parameter with $\siminfB = \siminfc$. 
Then $c$ satisfies the hypotheses of \corref{sub_post_satelite_case} and for any nest $\{Y_n\}_n$ with 
$c\in\ov{Y_n}$ for all $n$, $v_B\in\ov{P_n}$ for all $n$ where $P_n = \chi(Y_n)$ and 
 \begin{itemize}
\item 
$\End(\{Y_n\}_n) = \{c\}$ if and only if $\End(\{P_n\}_n) = \{v_B\}$ 
\item
$\End(\{Y_n\}_n)$ is the filled Julia set of a quadratic-like restriction of $Q_c^k$\\
if and only if \\
$\End(\{P_n\}_n)$ is the filled Julia set of a quadratic like restriction of $g_B^k$.
\end{itemize}
\ENDCOR

\section{Parabolic Parameter-Puzzles}\label{s:PPP}
 
 \subsection{Parabolic Parameter Puzzle}\label{s:parameterpuzzle}
In the $p/q$-wake $\W^\Mone(p/q)$, 
we define  the parameter  puzzle pieces using three different points of view. 
As a first definition, we take the universal parabolic $p/q$ graph $ \GP^n$ 
and the parametrization $ \Upsilon$ (see  \defref{Parametergraph}) 
to define a parameter parabolic graph.  
This way, the complementary regions, the parameter puzzle pieces, of level $n$ 
parametrize a holomorphic motion of the level $n+1$ dynamical graph $ \GP^B_{n+1}$. 
In particular this point of view  allows to compare pieces and 
annuli in the dynamical plane and in the parameter plane.  
We characterize then the parameter puzzle pieces as the set of parameters such that 
the critical value stays in the holomorphic motion of the same puzzle piece. 
The third characterization is in terms of laminations. 
A parameter puzzle piece of level $n$ corresponds to the set of parameters 
sharing up to level $n+1$ the same lamination associated to a center.

\REFDEF{Parametergraph} For $n\geq 0$, the parameter parabolic graph is defined by  $$\GPP_n:=\overline\W^{\Mone}(p/q)\cap \overline{\Upsilon^{-1}(\GP^n)}.$$ 
\ENDDEF
Note that $\Upsilon^{-1}(\GP^n)\subset\C\Sm\Mone$. For this reason, we add the accumulation  of  $\Upsilon^{-1}(\GP^n)$
 which consists of landing points of rays coming from the graph.  Those rays have angles which are pre-images of $\theta$ and $\theta'$, therefore they land at Misiurewicz parameters (see \lemref{Misiurewicz}). 
 Any $B\in \GPP_n\cap \LL^\Mone_{p/q}$  is a Misiurewicz parameter. 
 It is common landing point of exactly $q$ parabolic parameter rays in $\GPP_n$ 
 and the corresponding parabolic dynamical rays co-land at $v_B$ (see \lemref{Misiurewicz}).
The graph $\GPP_n$ consists in two parts : the sides which are parts of rays with landing points and the top which are short cuts. 

\begin{definition} Denote by $\PPP_n$ the set of parameter parabolic puzzle pieces of level $n$, 
they are the connected components of $\overline\W^{\Mone}(p/q)\setminus \GPP_n$ intersecting $\Mone$.
We write $PP_n(B)$ for the one containing the parameter $B$. 
\end{definition}
Puzzle pieces are either disjoint or nested, in which case they have different levels.

\REM There is a unique  parameter puzzle piece of level $0$  that we denote by  $PP_0$. 
It is the short-cutted version of the wake :  $\W^{\Mone}_0(p/q)$.
\ENDREM
\proof There  is  a unique connected component of $\W^{\Mone}(p/q)\setminus \GPP_0$ intersecting $\Mone$ since 
   $\GPP_0$ contains in  $\overline {\W^{\Mone}(p/q)}$ only the rays of angle $\theta$ and $\theta'$ with its short cut :   $\hat \gamma(\theta,\theta')$. Therefore the graph    $\GPP_0$ intersects $\Mone$ only at 
the root $B_{p/q}\in\overline {\W^{\Mone}(p/q)}$ of $\W^{\Mone}(p/q)$.  

\begin{lemma}\label{PP0} Let $\Bc\in PP_0$. 
The graph $ \GP^{\Bc}(1)$ admits a holomorphic motion 
$$\Psi_0 = \Psi_0^\Bc :PP_0(\Bc)\times \GP^{\Bc}(1)\to  \widehat\C
\qquad\textrm{such that}\qquad
\Psi_0(B,  \GP^{\Bc}(1))= \GP^B(1).
$$ 
\end{lemma}
\proof  From \propref{base_para_holomorphic_motion} the graph   $\GP_0^{\Bc}$ admits a holomorphic motion
parametrized by $\W_0^{\Mone }(p/q)$. 

The graph $ \GP^B(1)$ is defined by $ \GP^B(1)=  \GP_0^B\cup \GP_1^B$ where  
$$
\GP_{1}^B := g_B^{-1}(\GP_0^B\Sm \ckga_0^B)\cup \{\ckga^B_{1}\cup\ga^B_{1}).
$$ 

By \corref{vlocation} we have that $B \in\W^{\Mone}(p/q)$ 
if and only if $v_B\in\W_{B}(p/q)$.  
Hence, $v_B\notin \GP^B_0$ so that we can lift $\Psi_0(B, .)$ 
on $\GP_0^B\Sm \ckga_{0}^B$ to get a holomorphic motion 
of $g_\Bc^{-1}(\GP_0^\Bc\Sm \ckga_{0}^\Bc)$ :  
$$\Psi_0(B,z) = g_B^{-1}(\Psi_0(B,g_\Bc(z)) ).$$
The holomorphic motion of the short cuts $\ckga_{1}^\Bc\cup \ga_{1}^\Bc$, 
follows immediately from \propref{base_para_holomorphic_motion}
 \endproof

\REFLEM{PPN} Fix $n\geq 0$ and any $\Bc$ in a level $n$ parameter puzzle piece 
$PP_n(\Bc)$. 
There exists a holomorphic motion 
$$\Psi_n = \Psi_n^\Bc :PP_n(\Bc)\times \GP_{n+1}^{\Bc} \to \widehat\C
$$ such that for any $B\in PP_n(\Bc)$, 
 $\Psi_n(B,  \GP_{n+1}^{\Bc})= \GP_{n+1}^B$. 

Moreover $\Psi_n$ extends to a holomorphic motion of the union $\GP^{\Bc}({n+1})$ 
of all graphs upto and including $n+1$:
$$\widetilde\Psi_n = \widetilde\Psi_n^\Bc :PP_n(\Bc)\times  \GP^{\Bc} ({n+1})\to \widehat\C$$
by setting $\widetilde\Psi_n =\Psi_k$ on $PP_n(\Bc)\times\GP_{k+1}^\Bc$ for $-1\le k\le n$. 
\ENDLEM
\proof The proof goes by induction, \lemref{PP0}  provides a proof for   $n=0$.  In the induction  we prove that  for $B\in PP_{n}(\Bc)$ the  graph is $\GP_{n+1}^B=h_B^{-1}(\GP_{n+1})$.  
  
We define $\Psi_{n+1}$ by lifting of $\Psi_{n}$. Indeed,  for $B\in PP_{n}(\Bc)$, the critical value $v_B$ never crosses   $\GP_{n+1}^B$ . 
Otherwise, if $v_B\in \GP_{n+1}^B$,  $h_B(v_B)\in  \GP_{n+1}$  and $\Upsilon (B)= h_B(v_B)$ would be on $ \GP_{n+1}$. 
Therefore we can define $g_B^{-1}(\GP_n^B\Sm \ckga_n^B)$ it  coincides with  $h_B^{-1}(\GP_{n+1}\Sm \ckga_n)$ by induction. 
Then 
 $$
\GP_{n+1}^B := g_B^{-1}(\GP_n^B\Sm \ckga_n^B)\cup \{\ckga^B_{n+1}\cup\ga^B_{n+1}).
$$ 
By hypothesis of induction, for $B\in PP_{n}(\Bc)$ the graph $\GP_{n+1}^B$ equals  $h_B^{-1}(\GP_{n+1})$ since  $\ckga_{n}^B = h_B^{-1}(\ckga_{n})$ and 
$\ga_{n}^B = h_B^{-1}(\ga_{n})$. 
The holomorphic motion follows  from these considerations. 
\endproof

\begin{figure}[htbp]
\centerline{\includegraphics[height=80mm]{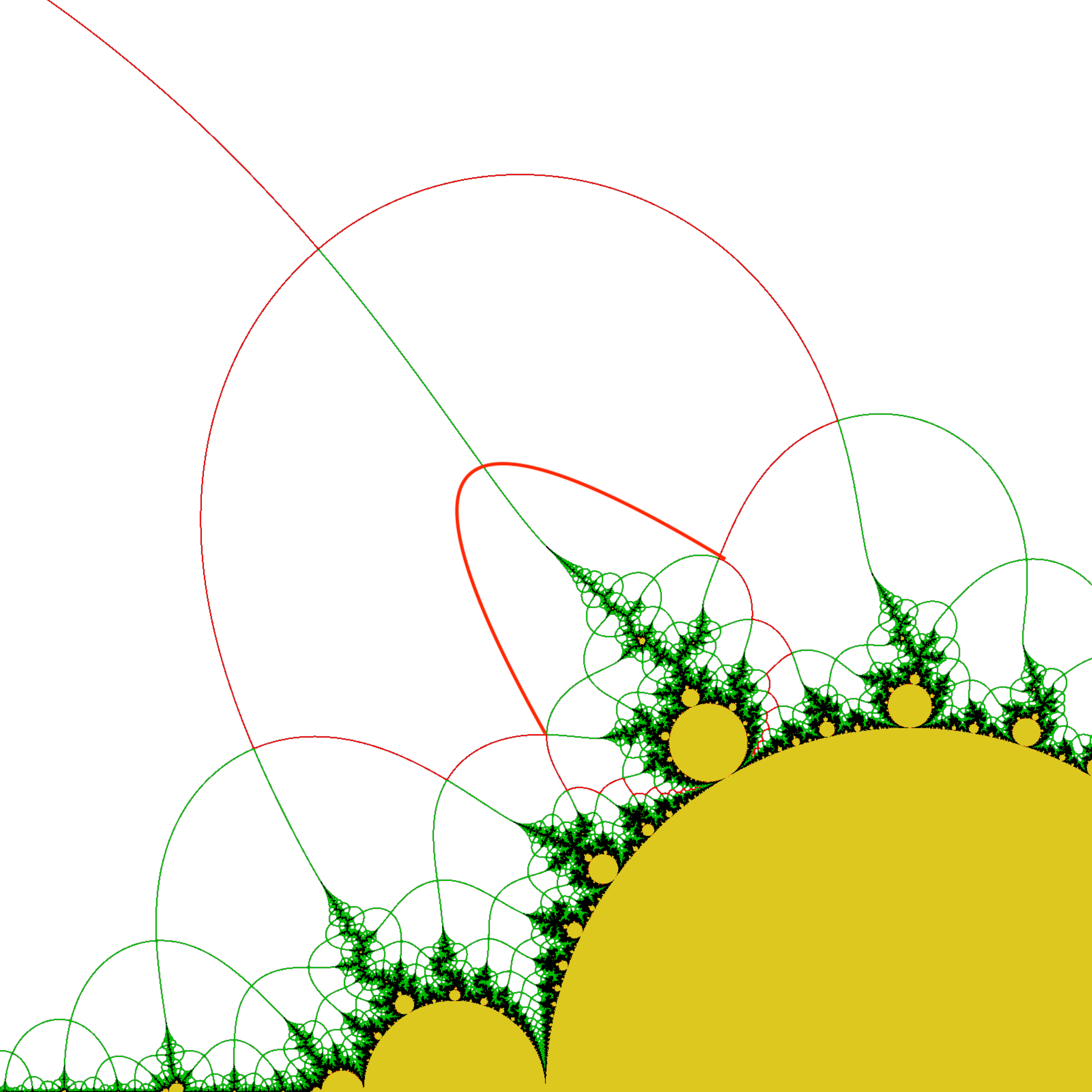}\includegraphics[height=90mm]{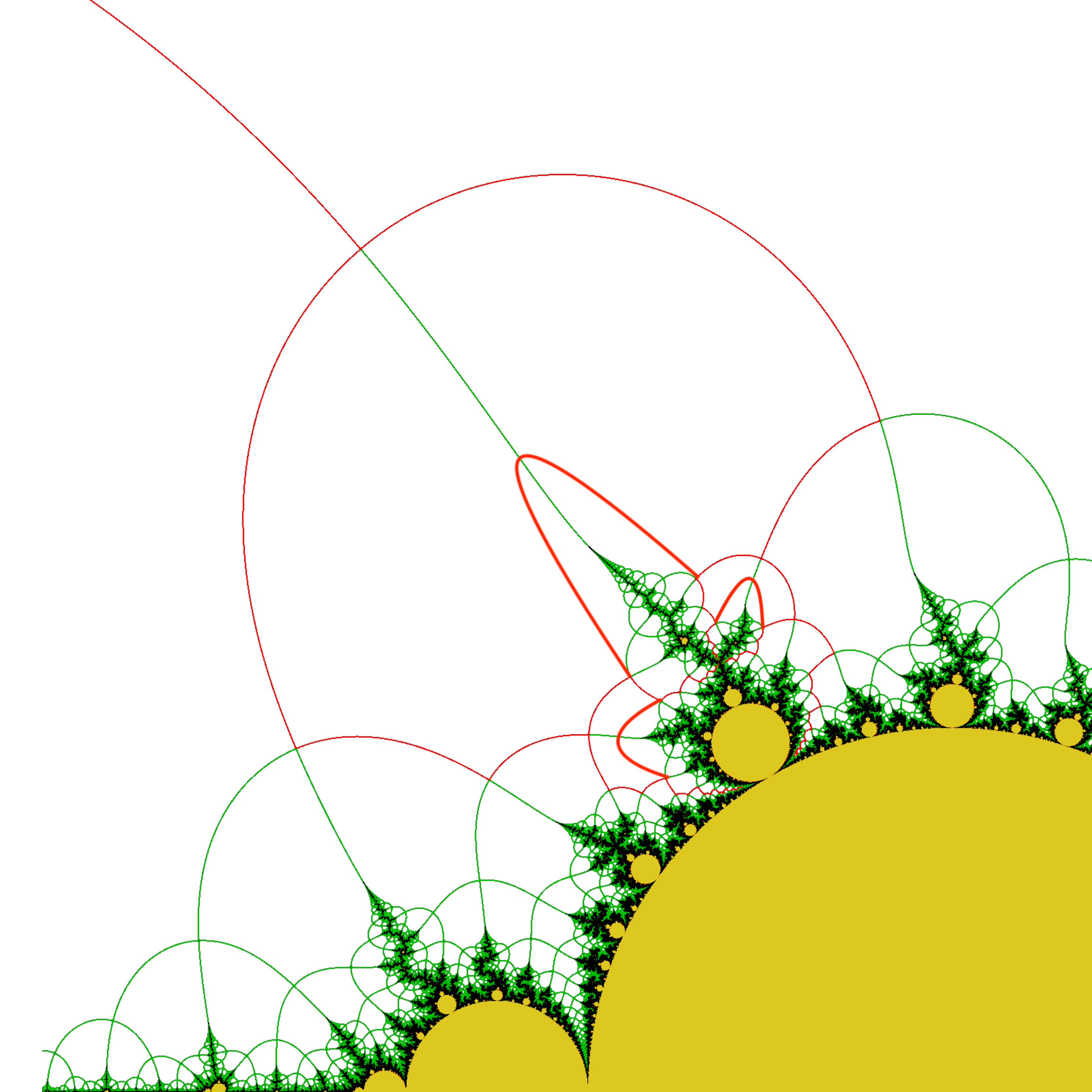}}
\caption{Parameter Parabolic Puzzle of first depths. }\label{parampuzzle}
\end{figure}

Let  $PP_n$ be a parameter puzzle piece, $\Bc\in PP_n$ and recall that $P_n^\Bc$ denotes 
the puzzle piece of level $n$ containing the critical value.  
We want to compare the situation in the parameter plane around $\Bc$ 
to the situation around the critical value $v_\Bc$ in the dynamical plane of $g_\Bc $ : 
compare the puzzle pieces and the annuli.  
Following the graph through the holomorphic motion, 
we have seen that the puzzle pieces up to level $n$ 
are homeomorphic so the situation is stable. 
Nevertheless, the critical value might cross the graph. 
Therefore the notion of puzzle piece containing the critical value is not continuous. 
For this reason, we give the name $\whP _n^B$ to the preferred puzzle piece, 
which is the holomorphic motion of this piece $P_n^\Bc$. More precisely, 
 \REFLEM{p:prefer}
For  $i\le n+1$ and  $B\in PP_n(\Bc)$, 
there is a unique parabolic puzzle piece $\whP_i^B$ bounded by the 
holomorphic motion $\Psi_n(B,\partial P_i^\Bc)$ of the critical value puzzle piece $P_i^\Bc$. 
\ENDLEM
 \proof  For $i\le n+1$ let $C_i^{\Bc}$ be  the boundary of the puzzle piece $P_i^{\Bc}$ containing the critical value for $g_\Bc$. It is a Jordan curve separating the critical value from $\GP_i^{\Bc}\setminus C_i^{\Bc}$. 
Following $C_i^{\Bc}$  in $PP_n(\Bc)$ through  the holomorphic motion  of the graph $\GP_i^{\Bc}$  defines a Jordan curve $C_i^B\subset \GP_{i}^B$ with $\GP_i^{B}\setminus C_i^{B}$ in a unique complementary component. Thus, we can define  the connected component of the complement of $C_i^B$  which is disjoint from $\GP_i^B$,   it   is our preferred puzzle piece denoted by $P_i^B$.
\endproof
  
\REFLEM{par_dyn_param_puzzle_p_corres} 
The puzzle piece $PP_{n+1}(\Bc)\subset PP_n(\Bc)$ is the set of parameters $B$ in   $PP_{n}(\Bc)$ such that  the prefered puzzle piece $\whP_{n+1}^B = P_{n+1}^B$. 
In particular $\whP_i^B = P_i^B$ for all $0\leq i \leq n$
\ENDLEM
\proof For parameters $B\in PP_{n+1}(\Bc)$, the critical value is clearly in $P_{n+1}^B$ since it never crosses the boundary of $P_{n+1}^B$. Then, being on the boundary of $PP_{n+1}(\Bc)$ and using the coordinates in the Blaschke model, one see that locally, if the parameter crosses the boundary  of $PP_{n+1}(\Bc)$ transversely,  
then the critical value follows the corresponding path in the dynamical plane. 
Hence it leaves the preferred puuzle piece $\whP_{n+1}^B$ 
and so has either to go into the puzzle piece adjacent to $\whP_{n+1}^B$ 
obtained by the holomorphic motion of the graph or to leave the level $n+1$ puzzle.
\endproof

\REFCOR{c:nondeg}
If $\overline{P_{n+1}^{\Bc}}\subset  P_n^{\Bc}$ then also $\overline{PP_{n+1}(\Bc)} \subset  PP_{n}(\Bc)$.
\ENDCOR 
\proof For a parameter $B$ in $PP_{n}(\Bc)$, 
the critical value belongs to $\whP_{n+1}^B\subset \whP_n^B = P_{n}^{B}$ 
and through the holomorphic motion we know that 
$\overline{\whP_{n+1}^{B}}\subset P_{n}^{B}$, 
so that if the critical value $v_B$ belongs to the annulus 
$P_{n}^{B}\setminus \overline{\whP_{n+1}^{B}}$, 
then parameter $B\in PP_n{\Bc}\Sm\overline{PP_{n+1}(\Bc)}$.
 \endproof
 
For $\Bc\in \Mone$, and $p/q$ such that $\Bc\in {\W^{\Mone}(p/q)}$, 
denote by  $ \sim^{\Bc}_\infty$ the lamination associated with 
the filled-in Julia set $K_\Bc$ ( see Section \ref{s:Tower}).
 
\begin{definition} Define $PP_n(\sim^{\Bc}_\infty) = PP(\sim_{n+1})$ 
to be the set of parameters $B$ in ${\W^{\Mone}(p/q)}$ 
such that $\sim^B_{n+1}=(\sim^{\Bc}_\infty)_{\vert_{n+1}} = \sim_{n+1}$.
 \end{definition} 

Recall from Section~\ref{s:parabopuzzle} 
that $V_n^\PP$ is the interior of the union of closures of level $n$ 
universal parabolic puzzle pieces and the reduced wakes $\W^\Mbrot_n(p/q)$ are 
$$
\W^\Mone_n(p/q) := \LL^\Mone_{p/q}\cup \{B\in\WMonepq | h_B(v_B)\in V_n^\PP\}, 
$$ 
so that 
\begin{lemma}
$PP_n(\sim^\Bc_\infty)\cap \W^\Mone_n(p/q) =PP_n({\Bc})$.
\end{lemma}

\proof
By definition $\Bc\in PP_n(\sim_\infty)$. Now  in $PP_n({\Bc})$ we have a holomorphic motion of the  parabolic rays in the graph $\GP_{n+1}^{\Bc}$  that gives the graph $\GP_{n+1}^B$.  Therefore, we keep the landing  relations for the parabolic rays in this graph in all the parameter puzzle piece. 
Hence, $PP_n({\Bc})\subset PP_n(\sim_\infty)$ 
and thus $PP_n(\sim_\infty)\cap  \W^\Mone_n(p/q)  \supset PP_n({\Bc})$.

For a parameter $B$ on the boundary of $PP_n({\Bc})$ the critical value $v_B$ is on the graph $\GP_{n}^{B}$, 
so either the second critical point $-1$ is on a pair of rays of the graph $\GP_{n+1}^{B}$ so that 
$\sim_{n+1}^\Bc\not= \sim_{n+1}^B$ 
or the critical value has escaped the level $n$ puzzle $\PP_n^B$. 
 \endproof

\REFPROP{para_enumerating}
For every $n$ there is 
\ENUM
\item\label{fertile_rea} a $1:1$ correspondence between the parameter puzzle pieces of level $n$ 
and the set of distinct fertile towers $\sim_{n+1}$ of level $n+1$. 
\item\label{terminal_rea} a $1:1$ correspondence 
between the set of level $n$ terminal towers and the set of points 
$$
\LL^*_{p/q}\cap(\GPP_n\Sm\GPP_{n+1}).
$$
\ENDENUM
\ENDPROP
\PROOF
The proof is by induction. For level $n=0$ there is only $1$ tower of level $n+1 =1$, it is fertile, 
the graph $\GP^B(1)$ moves holomorphically over $\W^\Mone_0(p/q)$ 
and induces the unique level $1$ tower $\sim_1$ and finally $\GPP_0$ does not intersect 
$\LL_*^\Mone(p/q)$.
Suppose the statement holds for $n\geq 0$, let $\sim_{n+1}$ be any fertile tower, 
let $PP_{n} = PP(\sim_{n+1})$ be the corresponding level $n$ parameter puzzle piece 
and let $\Bc$ be a parameter therein, i.e. $\sim_{n+1} = \sim_{n+1}^\Bc$. 
By \lemref{PPN} the graph $\GP^\Bc(n+1)$ moves holomorphically over $PP_{n}$ 
and hence by definition of puzzle pieces $\Upsilon$ 
defines a homeomorphism between the parameter sub-graph
$\GPP(n+1)\cap\ov{PP}_n$ and the dynamical sub-graph $\GP^\Bc(n+1)\cap \ov{P}_n^\Bc$ 
and hence induces a $1:1$ correspondence between the parameter puzzle pieces of level $n+1$ contained in $PP_n$ the level $n+1$ dynamical puzzles pieces for $g_\Bc$ contained in the critical value piece $P_n^\Bc$. 
And a $1:1$ correspondence between the points of $\GPP_{n+1}\cap PP_n\cap\LL^*_{p/q}$ and 
the points of $\GP_{n+1}^\Bc\cap P_n^\Bc\cap K_\Bc$. 
The (dynamical) puzzle pieces are in $1:1$ correspondence with the gaps of $\sim_{n+1}$ 
contained in the critical value gap $G_{n+1}'$ of $\sim_{n+1}$, i.e. the gap of $\sim_n$, 
which is the image of the critical gap of $\sim_{n+1}$. 
And the graph points are in $1:1$ correspondence with the set of level $n+1$ classes contained in $G_{n+1}'$. 
By definition of towers each gap $G'\subset G_{n+1}'$ of $\sim_{n+1}$ 
defines the unique level $n+2$  fertile tower extension $\sim_{n+2}$ of $\sim_{n+1}$ with critical value gap $G'$
and their totality enumerates all level $n+2$ fertile tower extensions of $\sim_{n+1}$.
Similarly each class $K'\subset G_{n+1}'$ defines the unique terminal tower extension of $\sim_{n+1}$ 
with critical value class $K'$ and their totality enumerates all level $n+2$ terminal tower extensions of $\sim_{n+1}$.
\ENDPROOF

\begin{figure}[h]\label{flaminations}
\begin{center}  
\includegraphics[ height=4 in]{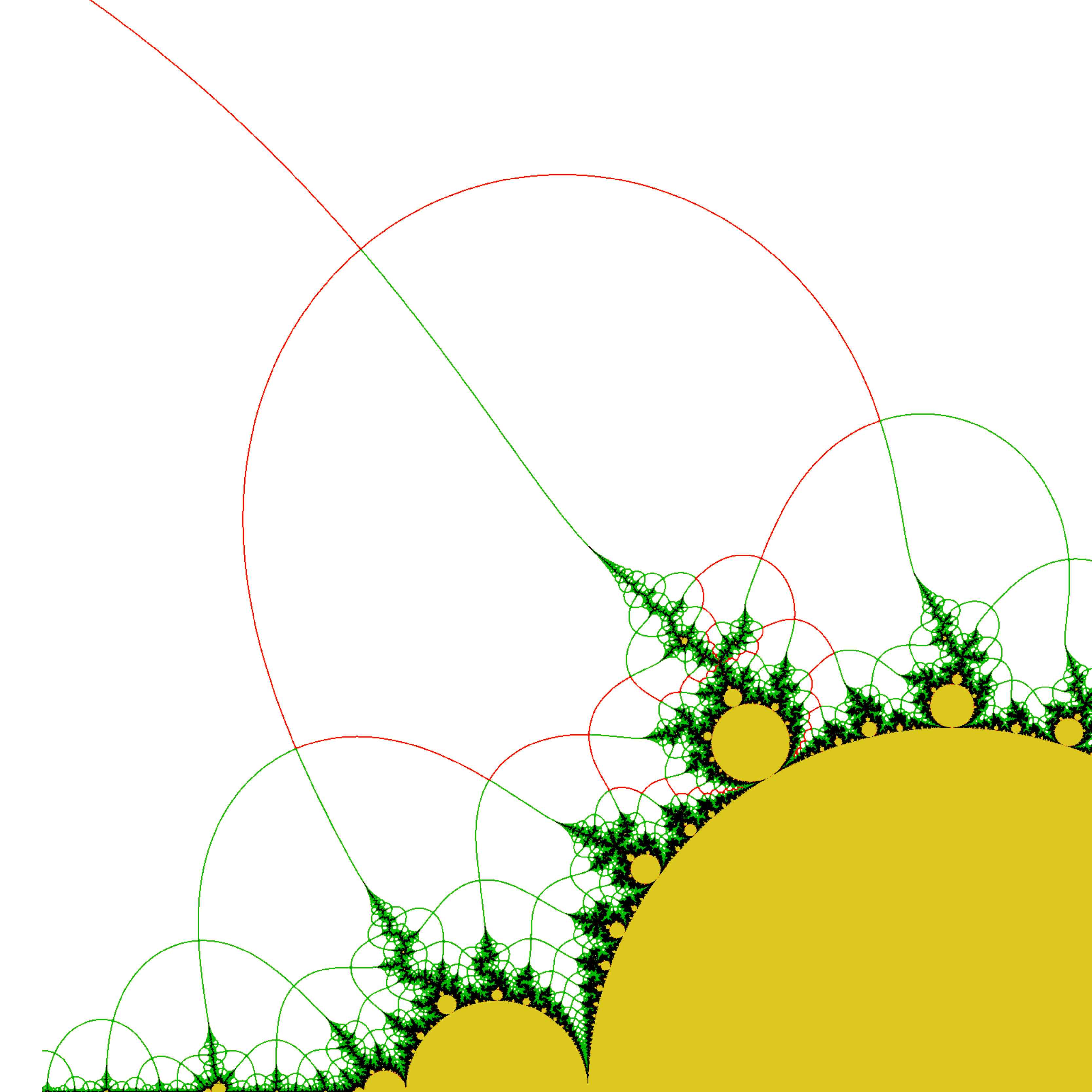}
\end{center}
\caption{Lamination for the puzzle pieces of first levels.}
\vskip -13.265em \hskip 23.165 em 
\includegraphics[height=0.5 in]{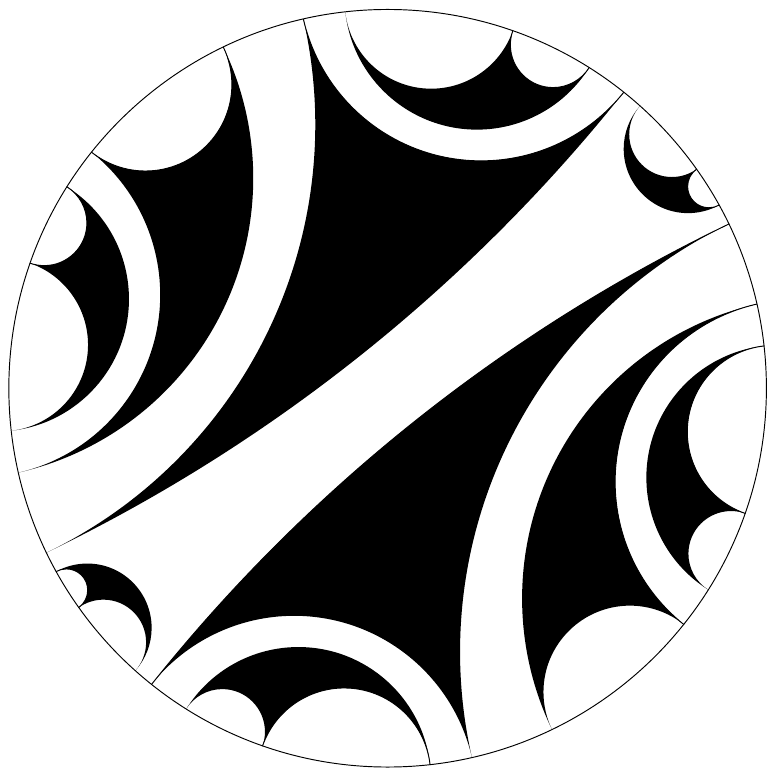} \vskip -13 em  \hskip 14em \includegraphics[height=1.5 in]{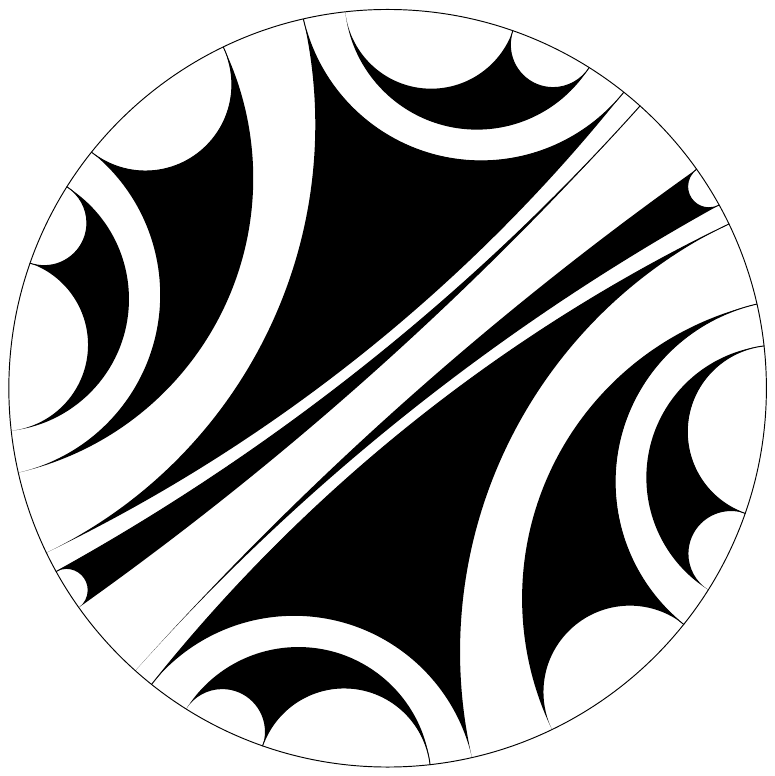}
 \vskip -2.4 em  \hskip 23.5em\includegraphics[height=0.35 in]{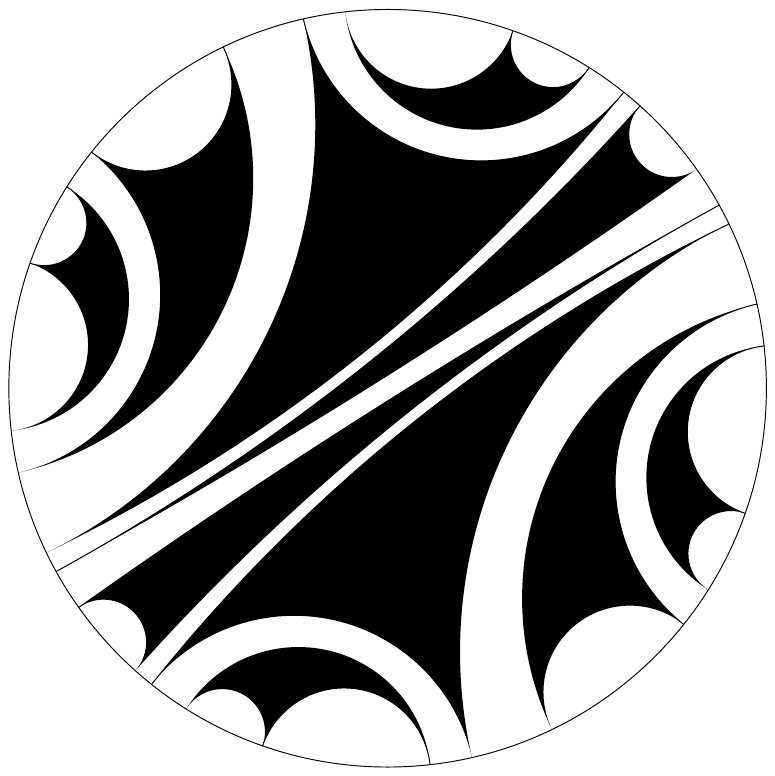}  
\end{figure}
\vskip 15 em

\REFDEF{terminal_tower_pieces}
In wiev of \itemref{terminal_rea} of the above Proposition we shall abuse notation and for every terminal tower 
$\sim$ write $PP(\sim)$ for the singleton consisting of the unique parameter $\Bc$ such that 
$\sim^\Bc = \sim$. For this parameter $g_\Bc^{n+1}(v_\Bc) = \al_\Bc$, where $n$ is the level 
of the critical value class for $\sim$.
\ENDDEF

\section{Parabolic Parameter Yoccoz Theorem, transfer to the parameter space}\label{Yoccoztheoremformone}

This section is devoted to proving that $\Mone$ is locally connected at any Yoccoz parameter, for a definition of such parameters see the item \Itemref{par_Yoccoz} below. 

Following Yoccoz approach to local connectivity of the Mandelbrot set we distinguish $3$ different types 
of parameters $B\in\Mone$:
\ENUMa
\item\label{par_alpha_non_repel}
Parameters $B\in\Mone$ such that the finite fixed point $\al_B$ is not repelling.
\item\label{par_g_B_renormalizeable}
Parameters $B\in\Mone$ such that some iterate $g_B^k$ is renormalizable around the second critical point $-1$ 
or equivalently around the second critical value $v_B$
\item\label{par_Yoccoz}
Parameters $B\in\Mone$ which is not in any of the two previous categories, also called Yoccoz parameters.
\ENDENUM

Local connectivity of $\Mone$ at a parameter $B$ of type \Itemref{par_alpha_non_repel} 
is most conveniently described in terms of the parameter $A=1-B^2\in\Dbar$. 
As a fundamental system of connected neighbourhoods of a parameter $A$ with $|A| = 1$ we may 
take a sequence of open intervals $J_n \subset\Sen$ shriking down to $A$ and with endpoints of irrational arguments, together with semi-disks in $\Delta_n\subset\D$, 
say bounded by the hyperbolic geodesic connecting the end-point of $J_n$ and together with the Limbs 
$\LL^\Mone_{p/q}$ with root in $J_n$. By \thmref{douadylandingthm}, \thmref{wake} and \corref{par_shringking_limbs} such sets form 
a fundamental system of connected neighbourhoods of $A$. 

We shall not here prove local connectivity of $\Mone$ at renormalizable parameters. 
It is not even known to be true in full genrality for the corresponding parameters in $\Mbrot$.
In fact our proof that $\Mone$ is homeomorphic to $\Mbrot$ works because it 
essentially does not rely on properties of renormalization copies beyond the second renormalization level.

In order to handle Yoccoz parameters we look to \secref{s:Transfer_to_Mone}.
Firstly we will consider one limb at a time, so we fix an irreducible rational $p/q$ and consider the limb
$\LL^\Mone_{p/q}$. 
Secondly we have the basic Dichotomy for such parameters :
\ENUMDB
\item%\label{par_satelite_case}
For all $n\in\N$ : $g_B^{nq}(-1)\in\ov{P_0^0}$.
\item%\label{par_post_satelite_case}
There exists $m\geq 1$ minimal such that $g_B^{mq}(-1) \notin \ov{P_0^0}$.
\ENDENUM
Where the first is equivalent to $g_B$ is $q$-renormalizable also denote immediate satelite type.
And the second is equivalent to $B\in\LL^\Mone(p/q,r,m)$ for some odd $r$ with $0< r < 2^m$ 
by \propref{par_wake_post_satelite_case_in_dyadic_wake}. 

Thridly by \corref{par_sub_post_satelite_case} 
the second condition $B\in\LL^\Mone(p/q,r,m)$ for some odd $r$ with $0< r < 2^m$ 
splits into three disjoint subsets or types of parameters:
\ENUMi
\item%\label{par_prebeta}
There exists $l\geq mq$ such that $g_B^l(-1) = \beta_B$.
\item%\label{par_capture}
There exists $l> mq$ such that $g_B^l(-1) \in K_B'$.
\item%\label{par_level_zero_recurrent}
There exists a strictly increasing sequence $\{l_n\}_{n\geq 0}$ with $l_0 = mq-1$ 
with $g_B^{l_n}(-1) \in S_1$ for all $n$.
\ENDENUM

All three cases will be handled by using holomorphic motions to define a homeomorphism from 
the boundaries of puzzle pieces surrounding $v_B$ in dynamical space to the boundaries of puzzle pieces surrounding $B$ in parameter space, where $B$ is any Yoccoz parameter in $\LL^\Mone(p/q,r,m)$.
In the two first types there are a sequence of boundaries puzzle pieces nesting down to $v_B$ 
which move holomorphically over a fixed domain where as in the third case 
the domain of holomorphic motion of puzzle piece boundaries shrink, when the level increases. 
Moreover as it appears in \thmref{p:yoccparabo} the third case splits-up further into two sub-cases 
renormalizable and not renormalizable. 

We shall apply variations of the following Proposition from the book of applications of holomorphic motions

For $B\in\W^\Mone_0(p/q)$ and $G^B\subset \GP^B$ a sub-graph consisting of the boundary of 
one or more puzzle pieces, 
not necessarily of the same level, 
denote by $\DD_{G^B}$ the connected component of $\Chat\Sm G^B$ 
containing the first critical point $1$. 
So that for any puzzle piece $P$ and $G=\partial P$ we have $\DD_{G^B} = \Chat\Sm\ov{P}$.

\REFPROP{par_transfer_to_param_space}
Let $\Bc\in\W^\Mone_0(p/q)$ and let $G^\Bc\subset \GP^\Bc$ 
be a sub-graph consisting of the boundary of one or more puzzle pieces $P_n$ with $v_\Bc\in\ov{P}_n$.
Suppose there exists a topological disk $U\subset\W^\Mone_0(p/q)$ with $\Bc\in U$ 
and a holomorphic motion {\mapfromto H {U\times G^\Bc}\Chat} with base point $\Bc$ and with 
$h_B(H(B,z)) = h_\Bc(z)$ for every $(B,z)\in U\times (G^\Bc\Sm J_\Bc)$. 
Let $G^B := H(B,G^\Bc)$ and suppose that 
the second critical value $v_B\in\DD_{G^B}$ on $U\Sm M$ for some connected compact set $M$. 
Then there exists a graph $G\subset M$ consisting of boundaries of parameter puzzle pieces 
$PP_n$ with $\Bc\in\ov{PP_n}$ such that $v_B\in G^B$ for all $B\in G$ and the map
$$
B\mapsto \zeta(B) := H_B^{-1}(v_B) : G \to G^\Bc, \qquad \textrm{where\qquad} 
H_B^{-1}(H(B,z)) = z
$$
is the restriction of a quasi-conformal homeomorphism $\zeta$, which is asymptotically conformal at $\Bc$.
Moreover for each $P_n$ with $v_\Bc\in\ov{P}_n$ and $\partial P_n\subset G^\Bc$ the pre-image 
$\zeta^{-1}(\partial P_n)\subset G$ is the boundary of a level $n$ parameter puzzle piece $PP_n$
with $\Bc\in\ov{\PP}_n$.
\ENDPROP
Note that the condition $h_B(H(B,z)) = h_\Bc(z)$ for every $(B,z)\in U\times (G^\Bc\Sm J_\Bc)$ 
means that $H$ is a holomorphic motion of puzzle piece boundaries, so that $G^B \subset \GP^B$ 
for every $B$. Thus if $v_\Bc\in P_n$ for some level $n$ and $\partial P_n\subset G^\Bc$, 
then the holomorphic motion $H$ coincides with the holomorphic motion $\Psi_{n-1}^\Bc$ of \lemref{PPN}
where ever both are defined and $H(B,\partial P_n)$ is the boundary of the prefered puzzle piece  
in the sence of \lemref{par_dyn_param_puzzle_p_corres}. 
And $B\in PP_n(\Bc)$ precisely when the preferred puzzle piece equals 
the critical value puzzle piece $\whP_n^B$. 
\PROOF
By Slodkowskis Theorem there exists a holomorphic motion extension 
{\mapfromto \wH {U\times\Chat} \Chat} of $H$. 
Define a quasi-regular map {\mapfromto \zeta U \Chat}, 
which is asymptotically conformal at $\Bc$ 
by $\zeta(B) := \wH_B^{-1}(v_B)$ and let $G = \zeta^{-1}(G^\Bc)$. 
Then by construction $G\subset M$ and 
$G$ consists of boundaries of parameter puzzle pieces and $\zeta$ has a non-zero degree over $G^\Bc$. 
Since $\Upsilon$ is univalent, the degree is $1$ so that the restriction 
{\mapfromto \zeta G {G^\Bc}} is a homeomorphism. 
Finally for each $P_n$ with $v_\Bc\in\ov{P}_n$ and $\partial P_n\subset G^\Bc$ the pre-image 
$\zeta^{-1}(\partial P_n)\subset G$ is the boundary of a level $n$ parameter puzzle piece $PP_n$
with $\Bc\in\ov{PP}_n$.
\ENDPROOF

\REFCOR{par_simple_param_nest_convergence}
Suppose for some parameter $\Bc\in\LL^\Mone_{p/q}$ and some nest $\NN = \{P_n\}_{n\geq 0}$ that 
$\{v_\Bc\}=\End(\NN)$ and that for some increasing sequence $\{n_k\}_{k\in\N}$ 
the graph $G^\Bc := \cup_k \partial P_{n_k}$ satisfies 
the hypotheses of \propref{par_transfer_to_param_space} then 
the corresponding parameter nest $\{PP_n\}_n$ with $\partial PP_n := \zeta(\partial P_n)$ 
is convergent with 
$$
\End(\{PP_n\}_n) = \{\Bc\}
$$
\ENDCOR

\REFCOR{par_type_i_and_ii_local_con}
Let $\Bc\in\LL^\Mone_{p/q}$ be a parameter satisfying 
\ItemDBref{par_post_satelite_case} of type 
\Itemref{par_prebeta} or \Itemref{par_capture} 
then the set $\Mone$ is locally connected at $\Bc$. 
If $v_\Bc$ is not prefixed to $\al_\Bc$ 
then intersection $\cap PP_n(\sim^\Bc)$ reduces to one point and thus 
$\{\Mone\cap PP_n(\sim^\Bc)\}_{n\geq 0}$ is a fundamental system 
of connected neighbourhoods of $\Bc$ in $\Mone$. 
And if $v_\Bc$ is prefixed to $\al_\Bc$ then $\Bc$ 
has a fundamental system of connected neighbourhoods consisting for each $n$ of the 
interior of the union of closures of the $q$ level $n$ parameter puzzle pieces with $\Bc$ on the boundary.
\ENDCOR

For \ItemDBref{par_post_satelite_case} type \Itemref{par_level_zero_recurrent} 
we need a refinement of \propref{par_transfer_to_param_space} above due to Shishikura.

Fix  any $\Bc\in \Mone$ of type iii) for the rest of the section.  
Recall that \thmref{p:yoccparabo}  provides a non degenerate annulus 
$A_{n_0}^\Bc = P_{n_0}\Sm\overline{P}_{n_0+1}$ 
between nested puzzle pieces of the parabolic Yoccoz puzzle 
for $g_\Bc$ with $g_\Bc^{n_0}(P_{n_0}) = P_0^{q-1}$, 
$g_\Bc^{n_0}(P_{n_0+1}) = S_1^k$ for some $0<k<q$ and a nested sequence of annuli $A_{n_i}^\Bc=P_{n_i}^\Bc\setminus\overline{ P_{n_i+1}^\Bc}$, 
$i>0$ with $n_0 < n_1\nearrow\infty$ surrounding the critical value $v_\Bc$ 
such that\,:
\begin{itemize}
\item the map $g_\Bc^{n_i-n_0} : A_{n_i}^\Bc \to A_{n_0}^\Bc$ 
is a covering map of degree $2^{d_i}$, $d_i\ge 0$ for $i\ge 1$\,;
\item 
in particular also all the annuli $A^\Bc_{n_i}$, $i>0$ are non degenerate\,;
\item
\begin{itemize}
\item
either the sum $\displaystyle \sum _{i\ge 1}\Mod(A_{n_i}^\Bc) = 
\Mod(A_{n_0}^\Bc)\sum _{i\ge 1} \frac{1}{2^{d_i}}$ is infinite,\\
\item
or there exists $k>0$ such that for all $n$ large enough the map\\ 
$$g_\Bc^k:P_{n+k}^\Bc\to P_n^\Bc$$ 
is quadratic-like with connected filled-in Julia set. 
\end{itemize} 
\end{itemize}
 
By \corref{c:nondeg}, the parameter $\Bc$ belongs to a complete sequence of puzzle pieces $PP_n(\sim)$  defining non degenerate annuli for the subsequence  $\AA_{n_i}(\sim)$  where   $\AA_{n}(\sim)$  denotes the  annulus  $PP_n(\sim)\setminus \overline {PP_{n+1}(\sim)}$. 

From \lemref{PPN}, for $n=n_0$, the  parameter puzzle piece  
$PP_{n_0}(\Bc)\subset  \W^\Mone_0(p/q)$,   parametrizes   a holomorphic motion  of the graph   $$\GP^\Bc_{n_0}\cup(\GP^\Bc_{n_0+1}\cap P^\Bc_{n_0}) \subset \GP^{\Bc}(n_0+1)$$  
as a restriction of  
 $$
 \widetilde\Psi_{n_0} = \widetilde\Psi_{n_0}^\Bc : 
 PP_{n_0}(\Bc)\times  \GP^{\Bc} ({n_0+1})\to \widehat\C
 $$

Then,  applying Slodkowsky's extension, we obtain a global holomorphic motion over $PP_n(\Bc) \subset \W^\Mone_0(p/q)$ of $\Chat$
(we are however only interested in the part inside $\ov{P}^\Bc_{n_0}$.)

Therefore by restriction, we get a holomorphic motion of the annulus $A^\Bc_{n_0} $, 
which gives an annulus that coincides with the annulus  $A^B_{n_0} = P^B_{n_0}\setminus \overline {P^B_{n_0+1}}$.

Shishikura's trick consists in lifting the holomorphic motion 
of the annulus to get a holomorphic motion of the annulus $A^\Bc_{n_i} $ 
defined in $\AA_{n_i} (\sim)$ with the same dilation. 
The lifting is possible since the map $g_\Bc^{n_i-n_0} : A_{n_i}^\Bc \to A_{n_0}^\Bc$ is a covering map of degree $2^{d_i}$, $d_i\ge 0$ for $i\ge 1$. 
Moreover, by  \lemref{par_dyn_param_puzzle_p_corres}, 
for parameters $B$  in  $PP_{n_i}(\Bc)$, 
the maps  $g_B^{n_i-n_0} : A_{n_i}^B \to A_{n_0}^B$ 
are all of the same type (covering map of degree $2^{d_i}$), 
since the critical value $v_B$ never passes though the boundary of $P_{n_i}^B$.

\begin{lemma}\label{modulus}
There exists a constant $K$ such that for any integer $n\in \{n_i\mid i\ge 0\}$
$$
\frac{\mod(A_n^\Bc)}{K}\le \mod \AA_N(\sim)\le K\mod(A_n^\Bc).
$$ 
\end{lemma}
\proof
The holomorphic motion of the boundary of $A^\Bc_{n_0}$ 
is defined in the whole para-puzzle piece $PP_{n_0}(\sim)$. 
By Slodkovskis Theorem we can extend it to a holomorphic motion 
of $\widehat \C$ still parametrized by $PP_{n_0}(\sim)$. 
Denote it by $H^0(B, z)$. 
Now, one can lift  $H^0(B, z)$ to 
a holomorphic motion $H^i(B, z)$ of $\overline{A_{n_i}^\Bc}$ 
using the unramified covering  $g_\Bc^{n_i-n_0}$. 
This holomorphic motion defines a quasi-conformal homeomorphism  
$H^i_{B}( z):=H^i(B, z)$, it has the same bound $K$ on the dilatation. 
Now, it follows from \cite{polylikemaps}Lemma IV.3, 
that the map $\zeta(B)=(H_{B}^i)^{-1}(v_B)$ is a quasi-conformal homeomorphism, 
it maps $\AA_n(\sim)$ to $A_n^\Bc$. 
The proof is exactly the same as in~\cite{PascaleinTanLeisbook}.
 \endproof

\REFCOR{par_local_con_Mone_type_iii_non_renorm}  
Let $\Bc\in\LL^\Mone_{p/q}$ be a parameter satisfying \ItemDBref{par_post_satelite_case} and of type \Itemref{par_level_zero_recurrent}. 
If $g_\Bc$ is not renormalizable or equivalently $\sim^\Bc$ is not renormalizable, 
then the intersection $\cap PP_n(\sim^\Bc)$ reduces to one point and thus 
$\{\Mone\cap PP_n(\sim^\Bc)\}_{n\geq 0}$ is a fundamental system of connected neighbourhoods 
of $\Bc$ in $\Mone$.
\ENDCOR
\proof  If $\sim^\Bc$ is not renormalizable, then the sum $\displaystyle \sum _{i\ge 1}\Mod(A_{n_i}^\Bc) = \Mod(A_{n_0}^\Bc)\sum _{i\ge 1} \frac{1}{2^{d_i}}$ is infinite, we deduce  from previous Lemma that   the sum 
$\displaystyle \sum _{i\ge 1}\Mod(\AA_n(\sim^\Bc))$ is infinite.  
Then the result follows from Gr\"otzsch inequality see~\cite{Ah}.\endproof
   
%%%%%%%%%%%%%%%%%%%%%%%%%%%%%%%%%%%%%%%%%%

\subsection{The renormalizable case}
We consider now a parameter  $\Bc\in\LL^\Mone_{p/q}$  satisfying \ItemDBref{par_post_satelite_case}, of type \Itemref{par_level_zero_recurrent} such that $g_\Bc$ is  renormalizable (equivalently $\sim^\Bc$ is  renormalizable). 
We use the Douady-Hubbard theory of polynomial like mapping 
to get that the intersection $\cap PP_n(\sim^\Bc)$ is a copy of the Mandelbrot set 
and we obtain in this way a straightening map that will serve 
to construct the bijection {\mapfromto {\Psi^1} \Mone \Mbrot}. 

\REFDEF{copy_of_M}\
A subset $\Mbrot_0$ of $\Mone$ is a {\it copy } of $\Mbrot$ if there exists a
homeomorphism $\chi$ and an integer $k>1$ (the {\it period}) such
that \begin{enumerate} \item$\Mbrot_0=\chi^{-1}(\Mbrot)$\,;\item
$\chi^{-1}(\partial \Mbrot)\subset\partial\Mone$ and\,;  \item  every
$B\in \Mbrot_0$  corresponds to a  renormalizable map  with $g^k$
topologically conjugated to $z^2+\chi(B)$ on
neighbourhoods of the filled Julia sets.
\end{enumerate}
\ENDDEF

\REFPROP{p:mandelbrot}
Suppose $\sim = \sim_\infty$ is a renormalizable tower of period $k\geq q$ 
and combinatorics $\sim_N$. 
Then for any $B\in PP_{N}(\sim) = PP(\sim_{N+1})$ the restriction 
$$
g_B^k : P_{N}^B \to P_{N-k}^B
$$
is quadratic like and the intersection 
$$
{\mathbf{M}}_{\sim}=\bigcap_{n\ge 0}\overline{PP_n(\sim)} = \bigcap_{n\ge 0}\overline{PP(\sim_{n+1})}
$$ 
is a copy of $\Mbrot$. 
\ENDPROP
Note that by definition $\Mbrot_\sim = \{B \mid\; \sim_\infty^B = \sim_\infty\} =: \Mbrot^\Mone(\sim_\infty)$.
\proof 
We develop here the case where the period is  $k\neq q$.
If the period is $k=q$, it corresponds to the satellite renormalizable case, 
the proof is similar except that one should consider enlarged puzzle pieces 
at the $\alpha$ fixed point for $P^B_n$.

The proof for $k>q$ is as follows. 
Let $c_0$ be the parameter with $Q_{c_0}^k(c_0) = c_0$ and $\sim_\infty^\cz = \sim_\infty$. 
Then {\mapfromto {Q_\cz^k} {Y^\cz_{N}}{Y_{N-k}^\cz}} is quadratic like and hybridly equivalent to $Q_0$. 
By \propref{samebehaviour} it follows that $g_B^k : P_{N}^B \to P_{N-k}^B$ is quadratic like 
for any $B\in PP_{N}(\sim) = PP(\sim_{N+1})$. 
Moreover the filled-in Julia set of this restriction is connected if $g_B^{mk}(v_B)\in P_{N}^B$ for all $m$. 

\noindent To simplify notation write $PP_n=PP_n(\sim)$ for all $n$. 
Consider the mapping $\mathbf{g}: \W'\to \W$ defined by
 $\W=\{(B,z) \mid B \in PP_N,\ z\in  P_{N-k}^B\}$, 
$\W'=\{(B,z) \mid B\in PP_N,\ z\in P_{N}^B\}$ and
$\mathbf{g}(B,z)=(B,g^k(z))$. 
It is an analytic family of quadratic-like maps in the sense of
Douady and Hubbard \cite[p.304]{polylikemaps} since it satisfy the
following three properties\,:

 \begin{itemize}
 \item the map $\mathbf{g}: \W' \to \W$ is holomorphic and
proper\,;

\item the holomorphic motion of the disk $P_N^B$, resp.
$P_{N-k}^B$, is a homeomorphism  between $\W$', resp.
$\W$, and $PP_N\times \D$ which  is fibered over $PP_N$ (since $B\in PP_N$)\,;

\item the projection  $\overline{ \W'} \cap \W \to  PP_N$ ({\it i.e.}
the first coordinate) is proper, since  
$\overline{ \W' }\cap \W=\{(B,z) \mid B \in PP_N,\ z\in\overline{P_N^B}\}$.
\end{itemize}

Let $\Mbrot_{\mathbf{g}}=\{B \in PP_N \mid K(g_B^k) \hbox{ is
connected}\}$ denote the connectedness locus of $\mathbf{g}$,
where $\displaystyle{K(g_B^k)=\bigcap_{i \ge 0}
(g_B^k)^{-i}(P^B_N)}$ denote its filled Julia set. 
Then $\Mbrot_{\mathbf{g}}$ coincides with $\Mbrot_{\sim}$. 
Indeed, for $B \in \Mbrot_{\sim}$, the critical point and its orbit under $g_B$ 
never cross the graphs.
Therefore the critical point of $g_B^k|_{P_n^B}$ does not escape the piece $P_N^B$ 
(by iteration by $g^k$). 
Hence $K(g_B^k)$ is connected and $B\in\Mbrot_{\mathbf{g}}$. 
Conversely, for $n\geq 0$ and $B \in PP_{N+nk}\setminus\ov{PP}_{N+(n+1)k}$, 
the common critical value $v_B$ for $g_B$ and $g^k_B$ restricted to $P^B_N$ 
belongs to the annulus $P^B_{N+nk}\setminus\overline P^B_{N+(n+1)k}$. 
Thus $g_B^{(n+1)k}(v_B)$ is not in $P_N^B$, that is the critical point of
$g_B^k$ escapes the domain. 
Hence the filled Julia set is not connected and so $B\notin \Mbrot_{\mathbf{g}}$.

Moreover, by \corref{c:nondeg} and \propref{p:yoccparabo},
there exists a sequence $n_i$ such that $\overline{PP}_{n_i+1}\subset PP_{n_i}$. 
Then $\Mbrot_{\sim}$ is also the intersection of the (open) pieces\,:
$\displaystyle \Mbrot_{\sim}=\bigcap_{n\ge0} PP_n$ and therefore  is
compactly contained in any of the parameter puzzle pieces $PP_m$, $m\ge N$.

 The theory of Mandelbrot-like families  of Douady
 and Hubbard  (see \cite{polylikemaps}, Theorem II.2, Propositions II.14 and IV.21)
 gives a continuous map $\chi : PP_N \to \C$ such that
 the maps  $g_B^k$ and  $z^2+\chi(B)$ are quasi-conformally conjugate on a
 neighbourhood of the filled Julia sets, for every $B \in PP_N$.

 Moreover, since $\Mbrot_{\sim}=\Mbrot_{\mathbf{g}}$ is compactly contained in $PP_N$, 
 the map $\chi$ induces a homeomorphism between $\Mbrot_{\mathbf{g}}$ 
 and the Mandelbrot set~$\Mbrot$  
 if we are in the following situation (see \cite{polylikemaps})\,:
 for a closed disk $\Delta\subset PP_n$
 containing $\Mbrot_{\mathbf{g}}$ in its interior,  the quantity $g_B^k(x_B)-x_B = v_B - x_B$,
 (where $x_B$ denotes the unique critical point of $g_B^k|_{P_n^B}$) turns exactly
 once around $0$ when $B$ describes $\partial \Delta$. 
 We verify this property in the following.

Let $n = N+k$ so that $\Mbrot_\sim\subset PP_n =: \Delta \subset\subset PP_N$. 
For $B\in PP_N$ both the boundary $\bd{P}_N^B$ and the critical point $x_B\in P_n^B$ 
move holomorphically with $B$ and thus $g_B^k$ is a locally diffeomorphic covering map of degree $2$ 
 from $\bd{P}_n^B$ onto  $\bd{P}_N^B$. Thus also $\bd{P}_n^B$ move holomorphically with $B$ over 
 the parameter disk $PP_N$.
 We compute the degree of $\gamma(B)=v_B-x_B$ on $\bd{\Delta}$. 
 Fix $\Bc\in PP_n$. 
 In order to do so we make use of the above holomorphic motions to transfer the problem to a 
 problem of winding number of a curve in the dynamical plane of $g_\Bc$
Let {\mapfromto H {PP_N\times(\bd{P}_n^\Bc\cup\{x_\Bc\})} \C}
 be the holomorphic motion with base point $\Bc$ just described and extend it to a global 
 holomorphic motion {\mapfromto H {PP_N\times \C} \C} using Slodkovskis theorem. 
 Let {\mapfromto \zeta {\bd{PP}_n} {\bd{P}_n^\Bc}} be the homeomorphism given by
 $H(B,\zeta(B)) = v_B$. 
 
Assume that $\Delta = PP_n$ is a round disk with center $\Bc$ (if not  use a conformal representation)\,; 
then the map {\mapfromto G {[0,1]\times\bd{\Delta}} \C} given by
$$G(t,B)=H(\Bc+t(B-\Bc),\zeta(B))- H(\Bc+t(B-\Bc),x_\Bc)$$ is a homotopy between 
$\zeta(B)-x_\Bc$ and $v_B-x_B$. And the degree of the first curve is simply the 
winding number $1$ of $\bd{P}_n^\Bc$ around $x_\Bc$. 
Hence $\Mbrot_\sim$ is a copy of the Mandelbrot set $\Mbrot$.
\endproof

\REFCOR{tower_surjectivity}
For every $c\in L^*_{p/q}$ there exists a $B\in\LL^*_{p/q}$ with $\sim_\infty^c\; =\; \sim_\infty^B$
\ENDCOR
\PROOF
Let $c\in L^*_{p/q}$ and let $\sim\; =\; \sim_\infty^c$. If $\sim$ is a terminal tower then there exists a unique 
$B\in \LL^*_{p/q}$ with $\sim_\infty^B\; =\; \sim_\infty^c$ by \Itemref{terminal_rea} of \propref{para_enumerating}. 
If $c\in\Mbrot^\Mbrot_{p/q}$, i.e. $Q_c$ is $q$ renormalizable take any $B\in\Mbrot^\Mone_{p/q}$ 
i.e. the unique $B$ such that $f_c$ and $f_B$ are hybrid equivalent. 

Finally if $c$ is any other parameter then $Y_n^c$ is defined for all $n$ and the parameter nest 
$\{YY_n(c)\}_n$ is a system of nested neighbourhoods of $c$. 
And for every $n$ there is $n' >n$ such that $\ov{Y}_{n'}^c\subset Y_n^c$ so 
that also $\ov{YY}_{n'}(c) \subset YY_n(c)$. 
By the theorems of this section the similar statement $\ov{PP}_{n'}(\sim) \subset PP_n(\sim)$ 
also holds so that for $B\in\cap\ov{PP}_n(\sim)$ we have $\sim_\infty^B\; =\; \sim\; =\; \sim_\infty^c$.
\ENDPROOF

%%%%%%%%%%%%%%%%%%%%%%%%%%%%%%%%%%%%%%%%%%
\section{Proof of the Main Theorem}\label{s:mainproof}
\subsection{The map $\Phi^1$ is a homeomorphism.}
In \cite[Proof of Theorem 1.1]{PR2} we have constructed a projection 
{\mapfromto {\Psi^1} \Mone \Mbrot} with the 
following properties, recall that $A = 1 - B^2$:
\ENUM
\item\label{ChypComp}
For $B\in\ov{\Hz}$ define $\Psi^1(B) := c$ where $c$ is the unique parameter such that 
the fixed point $\al_c$ for $Q_c$ has multiplier $A = 1- B^2$. 
\item\label{p_q_limb_condition}
For $B\notin\ov{\Hz}$ let $p/q$ be the irreducible rational such that $B\in\LL^\Mone_{p/q}$ and 
thus $\sim_\infty^B$ is well defined.
\ENUMa
\item\label{renomalizeable_case}
If $\sim_\infty^B$ is renormalizable of period $k$, then $g_B$ is $k$ renormalizable and 
$B\in M^\Mone = M^\Mone(\sim_\infty^B)$ a copy of $\Mbrot$ in $\Mone$. 
Let $M^\Mbrot := M^\Mbrot(\sim_\infty^B)$ be the copy of $\Mbrot$ in $\Mbrot$ 
such that $c\in M^\Mbrot$ if and only if $\sim_\infty^c = \sim_\infty^B$. 
Let {\mapfromto \chi {M^\Mone}{M^\Mbrot}} 
be the homeomorphism induced by straightening. Define $\Psi^1(B) = \chi(B)$.
\item\label{non_renormalizeable_case}
If $\sim_\infty^B$ is not renormalizable, i.e. a Yoccoz parameter let $c\in\Mbrot$ be the unique parameter such that 
$\sim_\infty^B = \sim_\infty^c$ and define $\Psi^1(B) = c$.
\ENDENUM
\ENDENUM

Define a map between parameter puzzles {\mapfromto \Xi \PPP \YY} by 
$\Xi(PP(\sim_n)) := YY(\sim_n)$, 
where the parameter puzzles $YY(\sim_n)$ for $\Mbrot$ are defined similarly to those for $\Mone$. 
Then by construction of $\Psi^1$ for any finite tower $\sim_n$ 
$$
\Psi^1(PP(\sim_n)\cap\Mone) = YY(\sim_n)\cap\Mbrot = \Xi(PP(\sim_n))\cap\Mbrot.
$$

By construction $\Psi^1$ is dynamic and unique:

For $B\in\ov{\Hz}$ let $c = \Psi^1(B)\in\Card$, 
we shall first construct using 
Ha{\"i}ssinsky's surgery a homeomorphism 
{\mapfromto {\rho_c} \Chat \Chat}, 
which is conformal a.e.~on the filled-in Julia set $K(Q_c)$ 
and which conjugates dynamics of $Q_c$ to that of some $g$ of the form 
$g(z) = z + B' + 1/z$ on their filled-in Julia sets.
\[\xymatrix{K(Q_c)\ar[r]^{Q_c}\ar[d]_{\rho_c} &K(Q_c)\ar[d]^{\rho_c} \\
K(g)\ar[r]_g&K(g) } \]
And then secondly see that $B = B'$.

If $A\in\D$ or is a root of unity, 
then clearly the critical point is not recurrent to the beta fixed point 
and Ha{\"i}ssinsky's theorem applies. 
And if $A\in\Sen$, but not a root of unity, then Dudko and Lyubich \cite{Dudko-Lyubich} 
have recently shown the existence of a "Mother hedgehog" for $Q_c$, 
a compact set $H_c$ containing the critical point $0$ and the fixed point $\al_c$ 
and such that the restriction {\mapfromto {Q_c}{H_c}{H_c}} is a homeomorphism. 
It follows that $\beta_c\notin H_c$ and hence that $0$ 
is not recurrent to $\beta_c$ in this case \cite{Dudko2}.
So that also in this case Ha{\"i}ssinsky's theorem applies. 

If $A\in\D$ of then the finite fixed point of $g$ also has multiplier $A$ 
since $\rho_c$ is conformal on the interior of $K(Q_c)$. 
And if $A\in\Sen$ it follows from Naishul's Theorem \cite{Naishul} that the 
the finite fixed point of $g$ also has multiplier $A$. 
Thus in either case $B' = B$. And $\Psi^1$ is uniquely defined on $\ov{\Hz}$. 

If $B\notin\ov{\Hz}$ and $c= \Psi^1(B)$, then we distinguish two cases. 
If $g_B$ is not renormalizable then the Julia 
sets are locally connected and the puzzle bijection $\Xi_B$ induces a conjugacy between 
the dynamics on the Julia sets. Moreover this conjugacy extends to a global homeomorphism, 
conformal a.e.~on $K(Q_c)$
since the Julia sets are locally connected and $K(Q_c)$ has measure $0$. 
In this case injectivity of $\Psi^1$ follows from uniqueness of the combinatorial invariant.
And if $g_B$ is renormalizable, 
then it is conjugate to $Q_c$ on the little Julia sets by straightening. 
And this conjugacy extends to a conjugacy on the Julia sets 
through the puzzle bijection $\Xi_B$. 
Thus the maps have the same combinatorial-analytic invariants. 
Existence of a global homeomorphism {\mapfromto {\rho_c} \Chat \Chat}, 
which is conformal a.e.~on the filled-in Julia set $K(Q_c)$ 
and which conjugates dynamics of $Q_c$ to that of some $g_B$ 
follows from Ha{\"i}ssinsky's theorem, because the critical point $0$ 
is not recurrent to $\beta_c$. 
This gives uniqueness also in this last case.

We proceed to show that $\Psi^1$ is a homeomorphism. 

Injectivity of $\Psi^1$ is an immediate consequence of 
\corref{par_type_i_and_ii_local_con} and \corref{par_local_con_Mone_type_iii_non_renorm}. 

For the surjectivity we need only to consider the case $c\in L^*_{p/q}$ for some irreducible rational $p/q$. 
Let $\sim\; =\; \sim^c$.

If $\sim$ is renormalizable let $M^\Mone = \cap PP_n(\sim)$, 
$M^\Mbrot = \cap YY_n(\sim)$ and {\mapfromto \chi {M^\Mone}{M^\Mbrot}} 
be as in \Itemref{renomalizeable_case}, then $\Psi^1(\chi^{-1}(c)) = c$ and we are done. 

If $\sim$ is not renormalizable then by \corref{tower_surjectivity} 
there exists $\Bc\in \LL^*_{p/q}$ with $\sim^\Bc\; =\; \sim^c$ and 
$\Psi^1(\Bc) = c$ by construction of $\Psi^1$.
Thus $\Psi^1$ is a bijection. 

We shall prove continuity of $\Psi^1$. The continuity of $\Phi^1 = (\Psi^1)^{-1}$ then follows 
since $\Psi^1$ is a continuous bijection between compact sets in a metric space. 

For the continuity of $\Psi^1$ fix $\Bc\in\Mone\Sm\ov{\Hz}$, $\cc = \Psi^1(\Bc)$ and 
$\sim_n = \sim_n^\Bc = \sim_n^\cc$, $n\geq 0$.

Let us start by noting that $\Psi^1$ is continuous at $\Bc$, 
when $\Bc$ is any Yoccoz parameter, 
by construction and by local connectivity of both $\Mone$ and $\Mbrot$ 
at all Yoccoz parameters. 
Indeed the parameter nest $\{PP_n(\Bc)\}_n = \{PP(\sim_n)\}_n$ 
form a fundamental system of neighbourhoods of $\Bc$ 
and similarly $\{YY_n(\cc)\}_n = \{YY(\sim_n)\}_n$ 
form a fundamental system of neighbourhoods of $\cc$. 

Thus we only need to prove that 
$\Psi^1$ is continuous at the boundary of any top-level renormalization copy 
$M^\Mone(\sim_\infty^\Bc)$ in $\Mone$. 
So let $\Bc\in M^\Mone\subset \LL^\Mone_{p/q}$ 
be a boundary point, $M^\Mbrot = M^\Mbrot(\sim_\infty^\Bc)$ 
and {\mapfromto \chi {M^\Mone}{M^\Mbrot}} be as in \Itemref{renomalizeable_case}.

We must show that $\Psi^1$ is continuous at $\Bc$. 
By construction $\Psi^1$ coincides with $\chi$ and so is continuous on $M^\Mone$. 
Hence we only need to show that if $\{B_n\}_n\subset \LL^\Mone_{p/q}\Sm M^\Mone$ 
is a sequence converging to $\Bc$, 
then the sequence $c_n := \Psi^1(B_n)$ 
converges to $\cc := \Psi^1(\Bc) = \chi(\Bc)$.

To this end we invoke the shrinking of dyadic decorations theorem. 
Recall that any renormalization copy comes equipped with dyadic limbs, 
which are the extremities of 
$\Mbrot$ or $\Mone$ beyond a renormalization copy (see e.g. \defref{para_p_q_dyadic_wakes}).

\RREFTHM{{\cite[Theorem 4]{PRcarrots}, \cite{Dudko}}}{shrinking_carrots}
For any copy $M$ of $\Mbrot$ in $\Mbrot$ or in $\Mone$
$$
\lim_{s\to\infty} \diam(L_M(r,s)) = 0
$$
where $\diam(\cdot)$ denotes Euclidean diameter and $L_M(r,s)$ denotes the $r/2^s$ dyadic limb of $M$, 
i.e. if $M= M^\Mbrot(\theta,\theta')$ then $L_M(r,s) = L^\Mbrot(\theta, \theta', r, s)$ and 
if $M= M^\Mone(\veps,\veps')$ then $L_M(r,s) = L^\Mone(\veps, \veps', r, s)$.
\ENDTHM

For $\{B_n\}_n$ a sequence converging to $\Bc$ as above let $L^\Mone_M(r_n,s_n)$ 
denote the dyadic limb of $M^\Mone$ containing $B_n$ and let $B_n'\in M^\Mone$ denote the root of that limb.
Then by construction $c_n=\Psi^1(B_n) \in L^\Mbrot_M(r_n, s_n)$ 
and $c_n' = \Psi^1(B_n')$ is the root of that limb.

Passing to a subsequence if necessary we can assume that 
either $B'_n$ and $(r_n, s_n)$ are eventually constant or $s_n$ diverges to infinity, 
since any two distinct dyadic limbs of $M^\Mone$ are strongly separated. 

If the sequence $s_n$ diverges to $\infty$, then both $|B_n-B_n'|\to 0$ and $|c_n-c_n'|\to 0$ as $n\to\infty$ 
by \thmref{shrinking_carrots}. 
Thus $B_n'\to \Bc$ since $B_n\to \Bc$ and hence $c_n'\to \cc$ as $n\to\infty$ since $\chi$ is continuous. 
And combining with $|c'_n-c_n|\to 0$ yields the desired $\Psi^1(B_n) = c_n\to \cc =\Psi^1(\Bc)$ as $n\to\infty$.

Next suppose $B_n' = B'$ and so $(r_n,s_n) = (r, s)$ for $n\geq N_0$. 
Moreover $c_n' = c' = \Psi^1(B')$ for $n\geq N_0$. 
Then $\Bc = B'$ since $\Bc, B'\in M^\Mone$  and $B'$ is the only intersection of $M^\Mone$ 
and $L^\Mone_M(r',s')$. And similarly $c' = \cc$. 

If the renormalization period $k>q$ then $M^\Mone = \cap PP_n(\sim)$ and 
$M^\Mbrot= \cap YY_n(\sim)$, where $\sim = \sim_\infty^\Bc = \sim_\infty^\cc$. 
Thus for every $N_1$ there exists $N_2$ such that $B_n\in PP_{N_1}(\sim)$ 
and hence $c_n\in YY_{N_1}(\sim)$
for every $n\geq N_2$. Hence $c_n\to c' = \cc$ as $n\to\infty$.

Finally if the renormalization period is $q$, so that $M^\Mone = \Mbrot^\Mone_{p/q}$ and 
$M^\Mbrot = \Mbrot^\Mbrot_{p/q}$.
Then $g_\Bc^{qs}(v_\Bc) = \al_\Bc$ and similarly $Q_\cc^{qs}(\cc) = \al_\cc$. 
Instead of developing augmented puzzles using more rays in the base puzzle, 
see e.g.~\cite{PRcarrots} and its illustrations we shall give here an ad hoc argument. 

Recall that for $\Bc\in\ov{\W^\Mone_{sq-1}(p/q,r,s)}$ 
we have defined a fundamental system of neighbourhoods $\{P_n^\al\}_n$ of $\al_\Bc$, 
such that $P_n^\al\cap K_\Bc$ is connected, where $P_n^\al$ 
is the interior of $\cup_{k=0}^{q-1}\ov{P_n^{\al,k}}$. 
Moreover the union of boundaries $\cup_n\partial P_n^\al$ move continuously with 
$B$ in the closure of $\W^\Mone_{sq-1}(p/q,r,s)$. 
Similarly for the quadratic polynomials $Q_c$, 
the union of boundaries $\cup_n\partial Y_n^\al$ move continuously with 
$c\in\ov{\W^\Mbrot_{sq-1}(p/q,r,s)}$. 
By hypothesis $g_{B_n}^{sq}(v_{B_n}) \to g_{\Bc}^{sq}(v_{\Bc}) =\al_\Bc$ as $n\to \infty$. 
Thus given $N_1$ there exists $N_2$ such that $g_{B_n}^{sq}(v_{B_n})\in P_{N_1}^{\al}\Sm P_0^0$ 
for every $n\geq N_2$. 
By construction  $g_{B_n}^{sq}(v_{B_n})\in P_{N_1}^{\al}\Sm P_0^0$ if and only if 
$Q_{c_n}^{sq}(c_n) \in Y_{N_1}^\al\Sm Y_0^0$ for every $n$ and $N_1$.
Thus $Q_{c_n}^{sq}(c_n) \in Y_{N_1}^\al\Sm Y_0^0$ for every $n\geq N_2$, i.e. 
also $Q_{c_n}^{sq}(c_n) \to Q_\cc^{sq}(\cc) = \al_\cc$. 
Finally $Q_\cc^{sq}$ is a local diffeomorphism from a neighbourhood of $\cc$ 
to a neighbourhood of $\al_\cc$ and the map $(c,z) \mapsto Q_c^{sq}(z)$ is holomorphic, 
so that $Q_{c_n}^{sq}(c_n) \to Q_\cc^{sq}(\cc) = \al_\cc$ implies $c_n\to\cc$ 
as $n\to\infty$. This completes the proof that $\Psi^1$ and $\Phi^1$ are homeomorphisms.

\subsection{The homeomorphism $\Phi^1$ is nowhere H{\"o}lder on $\partial\Mbrot$}
Let $c_0\in L_{p/q}$ be the $m/2^n$-dyadic tip of $\Mbrot$, where $0<m<2^n$, $m$ odd, i.e. 
the landing point of the parameter ray $\RR_\theta^\Mbrot$ of argument $\theta = m/2^n$. 
Let $\dN_{c_0}(\beta) = \dN(\beta) := \{\dY_n(\beta)\}_{n\geq 0}$ 
denote the nest around $\beta = \beta_{c_0}$ in the dynamical plane of $Q_{c_0}$. 
And let $\pN(c_0) := \{\pY_n(c_0)\}_{n\geq 0}$ denote parameter nest around $c_0$. 
Recall that the set
$$
\Ga_{c_0} := \{\beta\}\cup\bigcup_{n\geq 0} \bd\dY_n(\beta)
$$
moves holomorphically and equivariantly with $c$ in $\pY_0(c_0) := \pY_0$.
 
Let {\mapfromto H {\pY_0\times\ov{\dY_0(\beta)}} \C} 
be a holomorphic motion extending this motion, 
so that for for all $n\geq 0$ and all $c\in\pY_0$ : $H(z,\beta) = \beta_c$ and 
{\mapfromto H {\bd\dY_n(\beta)} {\bd\dY_n(\beta_c)}} is a homeomorphism with 
$Q_c\circ H(c,z)= H(c,Q_{c_0}(z))$, where both sides are defined.

For $n\geq 0$ denote by $z_n$ the unique point in $\bd{\dY_n(\beta)}\cap K_{c_0}$, in particular $z_0 = \al'$. 
Equivalently let {\mapfromto {\psi  = \psi_{c_0}} \C \C} denote the linearizer of $Q_{c_0}$ at $\beta$, 
normalized by $\psi(1) = \al'$. 
Let $\rho = Q_{c_0}'(\beta)$ denote the multiplier of $\beta$. 
Then $z_k = \psi(\rho^{-k})$ and 
asymptotically $(z_k-\beta)\rho^{-k} \to A$ as $k\to\infty$ 
for some non-zero complex number $A$. 
Let $z_n(c) = H(c, z_n)$ denote the motion of $z_n$ under the holomorphic motion.
Let $l$ be minimal such that $Q_{c_0}^k(c_0)\notin \dY_l(\beta)$ for $0\leq k < n$. 
Then for every $c\in\pY_{n+l}(c_0)$ the restriction 
{\mapfromto {Q_c^n} {\ov{\dY_{n+l}(c)}}{\ov{\dY_l(\beta_c)}}} is a holomorphic diffeomorphism. 
Define a quasi-conformal homeomorphism, asymptotically conformal at $c_0$, 
{\mapfromto \xi {\pY_{n+l}(c_0)}{\dY_l(\beta)}} by 
$$
\xi(c) := H_c^{-1}(Q_c^n(c)),\qquad\textrm{where}\qquad H_c^{-1}\circ H(c,z) = z.
$$
Note that a priori this map is a proper quasi-regular map, 
but the degree on the boundary is $1$, 
so that indeed it is a q.c.-homeomorphism. 
Let $c_k$, $k\geq l$ be the sequence of parameters $c_k = \xi^{-1}(z_k)\subset\pY_{n+l}(c_0)$ 
so that for $k\geq l$ 
$Q_{c_k}^n(c_k) = z_k(c_k)$ and thus $Q_{c_k}^{n+k}(c_k) = \al'(c_k)$. 
In particular $c_k\in L_{p/q}$ are Misiurewicz parameters and belong to $\bd\Mbrot$. 
Moreover by \lemref{assymptotic_conformality} below 
$$
\frac{(c_k-c_0)}{\rho^k} \longrightarrow \frac{A}{a} \qquad \textrm{as}\qquad k\to\infty,
$$
where $a\not=0$  is the difference of the derivatives at $c_0$ 
of the functions $Q_c^n(c)$ and $H(c, \beta)$.

Let $B_0 = \Phi^1(c_0)$, so that $g_{B_0}^n(B_0-2) = \beta_{B_0}$ and 
$g_{B_0}^{n-1}(B_0-2) = \beta'_{B_0}$. 
Recall that $\beta_B \equiv \infty$ and $\beta'_B \equiv 0$. 
For each $B$ near $B_0$ let 
{\mapfromto {\psi_{B}}\C \Chat} be the repelling Fatou-parameter for 
$g_B$ normalized by $\psi_B(0) = \al'_B$. 
Then $\psi_B$ depends holomorphically on $(B,z)$.
Define similarly as above a sequence of iterated pre-images of $\al_{B_0}$, 
$\{\whz_k = \psi_{B_0}(-k)\}_{k\geq 0}$ converging directly to $\beta_{B_0}$, 
i.e.~ $\whz_0 = \al_{B_0}'$, $g_{B_0}(\whz_{k+1}) = \whz_k = \dP_k(\beta_{B_0})\cap K_{B_0}$ for $k\geq 0$ and $\whz_k\to\beta_{B_0}$ as $k\to\infty$. 
Moreover let $\whz_k(B) = \psi_B(-k)$ be the motion of $\whz_k$ under the equivariant holomorphic motion on the parameter piece $\pPP_0(B_0) = \pPP_0$. 

Continuing similarly as above let $\{B_k\}_{k\geq l}\subset\pP_{n+l}(B_0)$ be 
the sequence of parameters such that $g_{B_k}^n(B_k-2) = \whz_k(B_k)$ for $k\geq l$. 
Then $\Phi^1(c_k) = B_k$ for $k\geq l$ and $B_k\to B_0$ as $k\to\infty$. 
And since $g_B(1/z) = g_B(z)$ the preimage of $\whz_k(B)$ near $0$ is $1/\whz_{k+1}(B)$,  
hence $g_{B_k}^{n-1}(B_k-2) = 1/\whz_{k+1}(B_k)$. 
Moreover since $\psi_B(z) = Bz + B^2\log(-z) + O(1)$ at infinity 
we have $\whz_k/k = \psi_{B_0}(-k)/k\simeq -B_0$ as $k\to\infty$. 
And thus $k/\whz_{k+1} \to -1/B_0$ as $k\to\infty$. 

Applicating \lemref{assymptotic_conformality} similarly to above 
we obtain 
$$
(B_k-B_0)\cdot k \longrightarrow \frac{-1}{bB_0}
\qquad\textrm{as}\qquad k\to\infty
$$
where $b$ is the derivative of $B\mapsto g_{B}^{n-1}(B-2)$ at $B = B_0$.

Finally we obtain that for any exponent $\kappa>0$ 
$$
\left|\frac{\Phi^1(c_k) - \Phi^1(c_0)}{{(c_k - c_0)}^\kappa}\right| 
= \left|\frac{B_k - B_0}{{(c_k - c_0)}^\kappa}\right| 
\geq \frac{1}{2}\frac{|a|^\kappa }{|bB_0||A|^\kappa} \frac{|\rho|^{k\kappa}}{k} 
\underset{k\to\infty}\longrightarrow\infty.
$$
Hence $\Phi^1$ is not H{\"o}lder-$\kappa$ for any $\kappa > 0$ at $c_0$. 
Since the dyadic tips are dense in the boundary of $\Mbrot$ the Theorem is proved.
 
 For completenes let us state precisely the asymptotic conformality result used above. 
 A variant can be found in \cite[Lemma in prof of Prop. 20.]{polylikemaps}. 
 Let {\mapfromto H {\D\times\D} \C} be a holomorphic motion with base point $\la_0=0$ and  
 let {\mapfromto f \D \C} be a holomorphic map. Then the expression 
$H_\la^{-1}(f(\la)) := \phi(\la)$ defines a quasi-regular map with $\phi(0) = 0$ 
in a neighbourhood of $0$. Let {\mapfromto \zeta \D \C} be the holomorphic map 
(motion of $0$) $\zeta(\la) := H(\la,0)$ and set $\si := f'(0) - \zeta'(0)$.
\REFLEM{assymptotic_conformality}
If $\si\neq 0$ then $\phi$ is quasi-conformal on a neighbourhood of $0$. 
Moreover $\phi$ is asymptotically conformal at $0$ with derivative $1/\si$ 
in the sense that 
$$
\lim_{\la\to 0} \frac{\phi(\la)}{\la} = \frac{1}{\si}.
$$
\ENDLEM
This proves Theorem B from, which the property that $\Phi^1$ admits no quasi-conformal 
extension to any neighbourhood of any boundary point of $\Mbrot$ easily follows, 
since any $K$-quasi-conformal homeomorphism is $1/K$-H{\"o}lder.

Aknowledgement:  The authors would like to thank: Arnaud Ch\'eritat for producing illustrations for figures 3, 9, 13, 16, 24 and C.~T.~McMullen for helpful suggestions 
C.~L.~P would like to thank Roskilde University, the Danish Council for Independent Research $|$ Natural Sciences for support via grants
DFF-1026-00267 and DFF-4181-00502. P.~Roesch would like to thank Lone at first for her patience, the University of Lille, Roskilde University and University of Toulouse, 
as well as the ANR for the grant ABC. 
Her research was also partial supported by ANR project LAMBDA,  ANR-13-BS01-0002.   
This work was also supported by the ANR LabEx CIMI (grant ANR-11-LABX-0040) within the French State Programme "Investissements d'Avenir."

Addresses:

Carsten Lunde Petersen, 
Department of Mathematical Sciences, 
University of Copenhagen, 
Universitetsparken 5, 
DK-2100 Copenhagen {\O}.
e-mail: lunde@math.ku.dk

Pascale Roesch,  
Institut de Math\'ematiques de Toulouse, 
Universit\'e Paul Sabatier,  Toulouse 3,   
31062 Toulouse Cedex 09, 
France. 
e-mail: roesch@math.univ-toulouse.fr

\end{document}